%% file: main.tex
\documentclass[ijoc,nonblindrev]{informs3_nojournal}

\OneAndAHalfSpacedXII %

\usepackage{natbib}
 \bibpunct[, ]{(}{)}{,}{a}{}{,}%
 \def\bibfont{\small}%

\TheoremsNumberedThrough     %

\EquationsNumberedThrough    %

\graphicspath{{Images/}}
\usepackage{bm}
\usepackage{bbm}
\usepackage{algorithm}

\usepackage{booktabs}
\usepackage{diagbox} %
\usepackage{subfig}
\usepackage{graphicx}
\usepackage{enumerate} 
\usepackage{algcompatible}

\newcommand{\calQ}{\mathcal{Q}}

\newcommand{\bbR}{\mathbb{R}}

\newcommand{\sA}{\mathcal{A}}
\newcommand{\sR}{\mathcal{R}}
\newcommand{\pset}{\mathcal{P}_{\delta,n,\rho}}

\newcommand{\calP}{\mathcal{P}}
\newcommand{\costF}{\text{Cost}(F)}
\newcommand{\espbullet}{\epsilon^{\bullet}(\{x_t,p_t,\theta_t\}_{t=1}^T)}
\newcommand{\espcirc}{\epsilon^{\circ}(\{x_t,p_t,\theta_t\}_{t=1}^T)}
\newcommand{\espi}{\epsilon^{i}(\{x_t,p_t,\theta_t\}_{t=1}^T)}

\newcommand{\hvt}{\hat{\vartheta}}
\newcommand{\kbullet}{\kappa^{\bullet}}
\newcommand{\kcirc}{\kappa^{\circ}}

\newcommand{\bF}{\mathbf{F}}
\newcommand{\bG}{\mathbf{G}}
\newcommand{\bu}{\mathbf{u}}
\newcommand{\sAhat}{\sA_{\hat{i}}}
\newcommand{\sAfeas}{\mathcal{A}_{\text{feas}}}

\begin{document}

\RUNAUTHOR{Im and Grigas}

\RUNTITLE{Stochastic First-Order Algorithms for Constrained DRO}

\TITLE{Stochastic First-Order Algorithms for Constrained Distributionally Robust Optimization}

\ARTICLEAUTHORS{%
\AUTHOR{Hyungki Im}
\AFF{Department of Industrial Engineering and Operations Research, University of California, Berkeley, Berkeley, CA 94720, \EMAIL{hyungki.im@berkeley.edu}} %
\AUTHOR{Paul Grigas}
\AFF{Department of Industrial Engineering and Operations Research, University of California, Berkeley, Berkeley, CA 94720, \EMAIL{pgrigas@berkeley.edu}}
} %

\ABSTRACT{%
We consider distributionally robust optimization (DRO) problems, reformulated as distributionally robust feasibility (DRF) problems, with multiple expectation constraints.
We propose a generic stochastic first-order meta-algorithm, where the decision variables and uncertain distribution parameters are each updated separately by applying stochastic first-order methods.
We then specialize our results to the case of using two specific versions of stochastic mirror descent (SMD): {\em (i)} a novel approximate version of SMD to update the decision variables, and {\em (ii)} the bandit mirror descent method to update the distribution parameters in the case of $\chi^2$-divergence sets.
For this specialization, we demonstrate that the total number of iterations is independent of the dimensions of the decision variables and distribution parameters. Moreover, the cost per iteration to update both sets of variables is nearly independent of the dimension of the distribution parameters, allowing for high dimensional ambiguity sets.
Furthermore, we show that the total number of iterations of our algorithm has a logarithmic dependence on the number of constraints. 
Experiments on logistic regression with fairness constraints, personalized parameter selection in a social network, and the multi-item newsvendor problem verify the theoretical results and show the usefulness of the algorithm, in particular when the dimension of the distribution parameters is large.
}%

\KEYWORDS{distributionally robust optimization, stochastic first-order methods, saddle point problems} 

\maketitle

\section{Introduction}
\label{submission}

\input{Introduction.tex}

\section{Stochastic First-Order Meta-Algorithm for DRF}\label{sec:Generic}
\input{Meta-Algorithm}

\section{$\epsilon$-Stochastic and Bandit Mirror Descent for DRF}\label{sec:SOFO}

\input{FirstOrderMethods.tex}

\section{Experiments}\label{sec:experiment}
\input{Experiments}

\section{Conclusion}
\input{Conclusion}

\ACKNOWLEDGMENT{%
This work was supported, in part, by NSF AI Institute for Advances in Optimization Award 2112533.
}%

\bibliographystyle{informs2014} %
\bibliography{ref.bib} %

\newpage
\begin{APPENDICES}
\input{Appendix}

\end{APPENDICES}

\end{document}

%% file: Introduction.tex
Recent increases in computational power, data availability, and improved modeling techniques have necessitated solving more and more real-world optimization problems under uncertainty. Models that perform robustly against changes in an uncertain environment are particularly appealing in both theory and practice. Distributionally robust optimization (DRO) \citep{wiesemann2014distributionally}, a popular paradigm that is robust against distribution shift, has gained great interest in both the operations research \citep{bertsimas2018data,taskesen2021sequential} and machine learning \citep{staib2019distributionally,duchi2021learning} communities.
DRO is a framework that seeks a robust solution that minimizes the worst-case expectation of the objective function in an ambiguity set of distributions $\calP$. 
Most of the existing work considers problems with a single distributionally robust objective or constraint, and DRO with multiple expectation constraints has not been studied as extensively.

There are important problems that are formulated with multiple expectation constraints, including supervised learning problems with fairness constraints \citep{zafar2017fairness, taskesen2020distributionally,akhtar2021conservative}, the Neyman-Pearson classification model \citep{tong2016survey}. For example, fairness in machine learning has been widely studied recently, since, without such constraints or other mechanisms, it has been observed that machine learning models can learn historical biases toward sensitive information such as gender and race \citep{raghavan2020mitigating}. Of course, this can be problematic if such models are used in important decision-making contexts, such as hiring decisions. At the same time, DRO has gained interest in the operations research and machine learning communities as a way to mitigate the risks of overfitting and distribution shift, and to generally increase model reliability \citep{shafieezadeh2015distributionally,mohajerin2018data,gotoh2021calibration}.

\par 
This paper considers a distributionally robust feasibility (DRF) problem with multiple expectation constraints. Namely, the task is to find $x \in X$ such that
 \begin{equation} 
    \underset{p^i \in \calP^i}{\sup} \left\{f^i(x,p^i):= \mathbb{E}_{z^i \sim p^i}[F^i (x,z^i)]\right\} \leq 0, 
 \label{form:DRO feas original}
 \end{equation}
for all $i \in \{1, \cdots,m \}$.
The domain of the decision variable $X \subseteq \bbR^d$ is a closed and convex set, and the set $\calP^i$ is a convex and compact ambiguity set for the $i$-th constraint. Although we consider feasibility problems, our methodology can solve a corresponding class of distributionally robust optimization problems by using binary search on the optimal value and repeatedly solving a small number of instances of \eqref{form:DRO feas original}.
In the model \eqref{form:DRO feas original}, the random variable $z^i$ follows the distribution $p^i$, which is only known to live in the ambiguity set of distributions $\calP^i$. Due to this ambiguity, the robust optimization methodology looks for a solution that satisfies the constraint for {\em all} possible distributions $p^i \in \calP^i$. Therefore, there is a trade-off between the robustness of the model and the final quality of the solution $x$ when the corresponding decision is implemented. As several authors have pointed out (see the discussion in Section \ref{sec:related}), a careful choice of the ambiguity set $\calP^i$, depending on the application context, can help balance this trade-off and enhance the effectiveness of model \eqref{form:DRO feas original}.

We define $Z^i$ as the range of possible values of $z^i$ and assume that $Z^i$ is a finite set of size $n$. As a result, $p^i \in \Delta_n$, where $\Delta_n: = \{p \in \mathbb{R}_+^n|\sum_{r=1}^n p^{(r)}=1 \}$ denotes a standard unit simplex, and $\calP^i \subseteq \Delta_n$. 
Throughout this paper, we assume that $F^i(x,z^i)$ is a $G$-Lipschitz convex function of $x$, for all $z^i \in Z^i$ and $i \in \{1, \dots,m\}$. Then, \eqref{form:DRO feas original} becomes a convex feasibility problem, albeit with a special structure. 
Under the assumption that $Z^i$ is a finite set, we can reformulate DRF to a robust feasibility problem (RF), which is the corresponding feasibility problem for robust optimization (RO). We discuss in detail this connection between DRF and RF in Section \ref{subsec:SP problem}. Based on this connection, algorithms for RF are clearly also applicable to solve the DRF problem we consider.
Recently, \citet{ho2018online} presented an algorithm based on first-order methods. 
Although their algorithm has cheap per-iteration costs and is scalable to the dimension of the decision variables, it lacks scalability to the dimension $n$ of the distribution parameters $p^i$.
Since the dimension of the distribution parameters is sometimes set to the number of data samples in DRF, the lack of scalability in $n$ makes the application of these algorithms to solve large-scale problems challenging.
In particular, the use of online first-order methods (FOMs) that use full gradient information is a crucial roadblock to the scalability of the algorithm of \citet{ho2018online}. 
\par
In this paper, we propose an online stochastic first-order meta-algorithm that is scalable to both the dimensions of $x$ and $p^i$. This meta-algorithm can be viewed as an extension of \citet{ho2018online} to a stochastic FOM version.
However, the use of usual stochastic FOMs, such as stochastic gradient descent (SGD) or stochastic mirror descent (SMD), turns out to be costly, as the complexity of obtaining a stochastic gradient for \eqref{form:DRO feas original} is linear in $n$. In particular, as discussed in Section \ref{sec:Generic}, we first need to obtain the index $i^*$ of the maximally violated constraint in \eqref{form:DRO feas original} in order to obtain an appropriate stochastic gradient for \eqref{form:DRO feas original}. In general, the complexity of obtaining $i^*$ is $O(n)$, which we discuss in Section \ref{subsec:epsilon SMD}. 
To overcome this challenge, we propose a novel variant of SMD called $\epsilon$-SMD. The $\epsilon$-SMD method uses a stochastic $\epsilon$-subgradient as its descent direction, and the complexity of obtaining the stochastic $\epsilon$-subgradient is much cheaper in this setting than other stochastic FOMs such as SGD and SMD. Furthermore, we present a customized version of the meta-algorithm that uses $\epsilon$-SMD as a stochastic FOM to update $x$ when the ambiguity set is the $\chi^2$-divergence set. 
Our contributions can be summarized as follows:
\begin{enumerate}
    \item We propose a stochastic first-order meta-algorithm for the DRF problem. Our meta-algorithm allows using arbitrary and possibly distinct ``basic'' stochastic FOMs to update the $x$ and $p^i$ variables, respectively. We show that our meta-algorithm for the DRF problem inherits high-probability convergence results from the basic stochastic FOMs. Also, we argue that high-probability convergence is a key ingredient that makes the meta-algorithm scalable to the number of constraints. 
    
    \item We propose an online stochastic FOM called $\epsilon$-stochastic mirror descent ($\epsilon$-SMD), which we apply as the basic stochastic FOM to update the decision variables $x_t$ in our meta-algorithm. We show that the per-iteration cost of $\epsilon$-SMD is independent of the dimension of the distribution parameters. Furthermore, we introduce an efficient feasibility testing algorithm that uses stochastic approximation and takes advantage of the information obtained while running $\epsilon$-SMD. We show that it is possible to determine the $\epsilon$-feasibility of a given solution independently of the dimensions of both the decision variable and the distribution parameters.
   
    \item We propose using the Bandit Mirror Descent (BMD) method of \citet{namkoong2016stochastic} to update the variables $p^i$ when the ambiguity set is a $\chi^2$-divergence set. BMD has a low cost per iteration to update the variables $p^i$. We extend the convergence analysis of BMD to a high-probability convergence guarantee and use this to show that the total iteration complexity of the customized meta-algorithm is almost independent of the dimension of the decision variables and distribution parameters.

    \item We examine the performance of our meta-algorithm using $\epsilon$-SMD and BMD in extensive numerical experiments on large instances for three problem classes:  logistic regression with fairness constraints, personalized parameter selection for a large-scale social network, and the multi-item newsvendor problem with a conditional value at risk (CVaR) constraint. We demonstrate performance improvements over the deterministic approach of \citet{ho2018online}, the state-of-the-art method for large-scale RF problems, when $n$ is large.
\end{enumerate}

\subsection{Related Work}\label{sec:related}
We review related works concerning DRO, RO, and associated algorithmic solution approaches.
In DRO, the ambiguity set $\calP$, the set of possible distributions of uncertain parameters, is a crucial design choice for the performance of the DRO model. 
Ideally, we want to choose an ambiguity set that balances accurate modeling (e.g., it encompasses the true distribution, if it exists), computational tractability, and an appropriate tuning of the level of conservativeness. 
Popular approaches include usßing moment-based constraints \citep{delage2010distributionally,hanasusanto2015distributionally,chen2019distributionally} or constraints based on divergence/distance functions \citep{duchi2021statistics,blanchet2019robust,blanchet2022optimal} to define the ambiguity set. In particular, there has been an extensive body of literature focusing on studying DRO models with ambiguity sets defined by $\phi$-divergences \citep{hu2013kullback, namkoong2016stochastic, duchi2021learning} and the Wasserstein distance \citep{kuhn2019wasserstein,gao2022wasserstein,gao2022finite}, largely due to their computational tractability and statistical properties. \cite{duchi2021learning}, among others, studied DRO with $\phi$-divergence constraints, demonstrating finite sample minimax upper and lower bounds while also showcasing its distributional robustness. Additionally, other divergence options, such as maximum mean discrepancy (MMD) \citep{staib2019distributionally, zhu2020kernel}, have also been explored. \cite{staib2019distributionally} studied DRO with MMD distance, demonstrating that MMD DRO is approximately equivalent to regularization by the Hilbert norm and established connections with kernel ridge regression.

Alternatively, rather than defining the ambiguity set as a constraint, it is frequently incorporated into the objective function as a penalization \citep{levy2020large,gotoh2021calibration,jin2021non, qi2021online}. \cite{gotoh2021calibration} investigated the out-of-sample performance of DRO solutions with $\phi$-divergence and demonstrated that calibrating the robust parameter, which determines the size of the ambiguity set, results in a greater variance reduction in out-of-sample loss with a minor trade-off in out-of-sample mean.
\par
There has been some progress in developing scalable algorithms for DRO under various ambiguity sets and settings. \citet{namkoong2016stochastic} proposed a primal dual algorithm to solve DRO with a $\phi$-divergence ambiguity set.
Also, due to the high cost of obtaining unbiased gradient estimators for the DRO objective, several works have focused on using biased estimators to solve DRO problems.
\citet{levy2020large} propose a primal algorithm that utilizes a mini-batch gradient estimator for DRO with a $\chi^2$-divergence set and CVaR-divergence set. Furthermore, they proposed an algorithm that uses a Multi-Level Monte Carlo (MLMC) estimator, in which the base sample size is proportional to $\ln(\frac{1}{\epsilon})$, where $\epsilon$ represents the tolerance rate in their settings. \cite{wang2021sinkhorn} recently introduced DRO with the Sinkhorn distance, which is a variant of the Wasserstein distance based on entropic regularization. They provide a dual reformulation of their problem and an efficient first-order algorithm utilizing biased gradients.
Our approach also uses a biased estimator for updating $x$, but it is important to note that the bias in our estimator stems from the presence of multiple constraints, whereas the bias in the aforementioned approaches results from the approximation error associated with the worst-case distribution for a given $x$.
\par 
As DRO has its roots in the robust optimization literature \citep{ben2013robust}, a strong connection can be established between DRO and RO, allowing for the exchange of ideas between the two sub-fields.
A traditional approach to solving robust convex optimization is to change the problem into an equivalent robust convex counterpart, which can be solved by state-of-the-art convex solvers \citep{ben2009robust,bertsimas2011theory,ben2015deriving}. However, the robust counterpart approach is often not scalable to the dimension of the decision variables. To remedy this scalability problem, \citet{ben2015oracle} propose an algorithm that alternatively updates the decision variables ($x$) and the uncertain parameters ($p$) using two oracles: a feasibility oracle and a pessimization oracle. Although their algorithm handles the scalability of the dimension of $x_t$, it still relies on the two oracles, which might be costly. \citet{ho2018online} introduced an algorithm that uses first-order methods.
Their algorithm has much cheaper iterations than the cost of the feasibility and pessimization oracles, is more scalable to the dimension of the decision variable $x$, and maintains the same convergence rates. In particular, \citet{ho2018online} bound the duality gap of their robust feasibility problem (in particular, its saddle point reformulation) by the sum of two weighted regret terms, and they update each weighted regret by using online convex optimization (OCO) methods \citep{shalev2011online,hazan2016introduction}. However, the use of deterministic gradients in their context makes it challenging for their algorithm to scale with respect to the dimension of the distribution parameters.

\subsection{Notation}
For any $n \in \mathbb{N}$, define $[n]:= \{1,\dots,n\}$. Superscripts are used to represent items corresponding to the $i$-th constraint.
We use $x^{(r)}$ to denote the $r-$th element of the vector $x \in \bbR^d$ and use $\mathbbm{1}_n$ to denote the $n$-dimensional vector of all ones. For a finite-dimensional real vector space $V$, $V^*$ denotes its dual space, and $\|\cdot\|, \|\cdot\|_*$ represents the primal norm and associated dual norm respectively.

%% file: Meta-Algorithm.tex
In this section, we present a stochastic first-order meta-algorithm for the distributionally robust feasibility (DRF) problem. We begin by reformulating \eqref{form:DRO feas original} to a convex-non-concave saddle point (SP) problem. All omitted proofs are included in the Appendix.
\subsection{Convex-Non-Concave Saddle Point Problem} \label{subsec:SP problem}
We begin by introducing the formal setting of our problem.
We assume that the random variable $z^i$ has finite support represented by $Z^i:= \{z_1^i, \cdots, z_n^i \}$ for all $i \in [m]$ and assume that there exists $M^i > 0$ such that $|F^i(x,z^i)| \leq M^i$, for all $x \in X$, $z^i \in Z^i$, and for each constraint $i \in [m]$. 
Also, we assume the random variables $z^i$ for all $i \in [m]$ are mutually independent. This assumption is without loss of generality and our results can be extended to dependent random variables $z^i$; we include a detailed discussion of the dependent case in the Appendix \ref{appn:dependent RV}.
We define $\calP:=\calP^1 \times \dots \times \calP^m$, which is the Cartesian product of sets $\calP^i$ for all $i \in [m]$. Define $F^i_r: X \rightarrow \bbR$ by $F^i_r(x):= F^i(x,z^i_r)$, for all $r \in [n]$. Furthermore, we define a vector-valued function $\bF^i(x): X \rightarrow \mathbb{R}^{n}$ by $\bF^i(x):=(F_1^i(x),\dots,F_{n}^i(x))$. Under our notation, we have $
    \mathbb{E}_{z^i \sim p^i}[F^i (x,z^i)] = (p^i)^\top\bF^i(x)$.
Then, we can write the DRF problem formally as:
 \begin{equation*}
    \begin{cases} 
    \text{feasible: exists $x \in X$ s.t. } \underset{i \in [m]}{\max} \  \underset{p^i \in \calP^i}{\sup} \ (p^i)^\top \bF^i(x)\leq 0 \\
    \text{infeasible: for all }  x \in X, \ \underset{i \in [m]}{\max} \  \underset{p^i \in \calP^i}{\sup} \ (p^i)^\top \bF^i(x)> 0.
    \end{cases}
 \end{equation*}
 Let us define a function $\phi(x,p):X \times \calP \rightarrow \mathbb{R}$ as $\phi(x,p):= \underset{i \in [m]}{\text{max}}  \ (p^i)^\top \bF^i(x)$. Note that $\phi(x,p)$ is a convex function of $x$ but not necessarily a concave function of $p$. After changing the order of the maximum and supremum above, the standard relaxation to an $\epsilon$-approximate feasibility problem is:
\begin{equation}
    \begin{cases}
    &\epsilon \text{-feasible:} \quad \underset{x \in X}{\text{inf}} \  \underset{p \in \calP}{\text{sup}} \  \phi(x,p) \leq \epsilon \\
    &\text{infeasible:} \quad \underset{x \in X}{\text{inf}} \    \underset{p \in \calP}{\text{sup}} \ \phi(x,p) >0 \label{form:DRO feas2}
\end{cases}    
\end{equation} 
We call $x \in X$ as a \textit{distributionally robust $\epsilon$-feasibility certificate} if it satisfies the $\epsilon-$feasibility condition of \eqref{form:DRO feas2}. Similarly, a realization of a distribution parameter $p \in \calP$ is called a \textit{distributionally robust infeasibility certificate} if $p$ satisfies $\phi(x,p)>0$ for all $x \in X$. 
As $\phi(x,p)$ is a convex-non-concave function, \eqref{form:DRO feas2} can be considered as a convex-non-concave saddle point problem (convex-non-concave SP).
We define the SP gap $\epsilon^{\phi}(\bar{x},\bar{p})$, which measures the accuracy of the given solution $(\bar{x},\bar{p}$) to \eqref{form:DRO feas2}, by 
\begin{align}
        \epsilon^{\phi}(\bar{x}, \bar{p}) := \underset{p \in P }{\text{sup}} \ \phi(\bar{x},p) - \underset{x \in X}{\text{inf}} \ \phi(x,\bar{p}). \label{eqn:SP gap}
\end{align}
\begin{theorem}[Theorem 3.1 of \citet{ho2018online}]\label{thm:Nam thm 3.1}
Let $\phi:X\times \calP \rightarrow \mathbb{R}$ be a convex function of $x$ but not necessarily a concave function of $p$ and $\epsilon>0$ be given. Suppose that we have $\bar{x} \in X, \bar{p} \in \calP$ and $\tau \in (0,1)$ such that $\epsilon^{\phi}(\bar{x},\bar{p}) \leq \tau \epsilon.$ Then if $\phi(\bar{x}, \bar{p}) \leq (1-\tau) \epsilon$, we have $\underset{p \in P}{\sup} \  \phi(\bar{x},p) \leq \epsilon$. Moreover, if $\phi(\bar{x},\bar{p})>(1-\tau)\epsilon$ and $\tau \leq \frac{1}{2}$, we have $\underset{x \in X}{\inf} \ \phi(x, \bar{p})>0$. 
\end{theorem}
Theorem \ref{thm:Nam thm 3.1} implies that it is sufficient to obtain a solution $(\bar{x},\bar{p}$) that satisfies $\epsilon^{\phi}(\bar{x},\bar{p}) \leq \frac{\epsilon}{2}$, for example, to solve \eqref{form:DRO feas2}.
Our meta-algorithm will obtain such a solution $(\bar{x},\bar{p}$) by iteratively generating $x_t \in X$ and $p_t \in \calP$ according to separate ``basic" stochastic FOMs, and then taking convex combinations of the generated points $\{x_t,p_t\}_{t=1}^T$. 
We define $\bar{x}_T = \sum_{t=1}^T \theta_t x_t$ and $\bar{p}_T = \sum_{t=1}^T \theta_t p_t$, given a vector of convex combination weights $\theta \in \Delta_T$.
\par 
To give an upper bound of $\epsilon^{\phi}(\bar{x}_T, \bar{p}_T)$, we define new error terms for given $\{x_t,p_t,\theta_t\}_{t=1}^T$:
\begin{alignat*}{3}
    &\epsilon^{\bullet}(\{x_t,p_t,\theta_t\}_{t=1}^T) ~:=~ &&\sum_{t=1}^T  \theta_t \phi(x_t,p_t) - \underset{x \in X}{\inf}  \sum_{t=1}^T  \theta_t \phi(x,p_t)\\
    &\epsilon^i(\{x_t,p_t^i,\theta_t\}_{t=1}^T) ~:=~ &&\underset{p^i \in P^i}{\sup} \sum_{t=1}^T \theta_t (p^i)^\top\bF^i(x_t)-\sum_{t=1}^T \theta_t (p_t^i)^\top\bF^i(x_t),\\
    &\epsilon^{\circ}(\{x_t,p_t,\theta_t\}_{t=1}^T) ~:=~ &&\underset{i \in [m]}{\max} \big\{ \espi \big\}.
\end{alignat*}
Then we have
\begin{align*}
    \epsilon^{\phi}(\bar{x}_T, \bar{p}_T) \leq \espbullet + \espcirc. %
\end{align*}
As a result, it is sufficient to solve the $\epsilon$-feasibility problem by obtaining $\{x_t,p_t,\theta_t\}_{t=1}^T$ that satisfies \begin{align}
    \epsilon^{\bullet}(\{x_t,p_t,\theta_t\}_{t=1}^T) + \epsilon^{\circ}(\{x_t,p_t,\theta_t\}_{t=1}^T)  \leq \frac{\epsilon}{2}. \label{eqn:term cond1}
\end{align}
\begin{remark}[Weighted Regret]
    The two terms $\espbullet$ and $\espi$ can be interpreted as weighted regret terms and we can use results from online convex optimization (OCO) to bound these terms respectively.
In particular, let us first consider $\epsilon^{\bullet}(\{x_t,p_t,\theta_t\}_{t=1}^T)$. Suppose $\theta \in \Delta_T$ and a sequence $\{p_t\}_{t=1}^T$ such that $p_t \in \calP$, for all $t \in [T]$, are given and let us define a function $\phi_t:X \rightarrow \mathbb{R}$ as $\phi_t(x):=\underset{i \in [m]}{\max} \  (p_t^i)^\top\bF^i(x)$. As $(p_t^i)^\top\bF^i(x)$ is a convex function of $x$, $\phi_t(x)$ is also a convex function of $x$. Then, we have 
\begin{align*}
    \epsilon^{\bullet}(\{x_t,p_t,\theta_t\}_{t=1}^T) &= \sum_{t=1}^T \ \theta_t \underset{i \in [m]}{\max} \ \phi(x_t,p_t) - \underset{x \in X}{\text{inf}} \ \sum_{t=1}^T \ \theta_t \underset{i \in [m]}{\max} \ \phi(x,p_t)\\
    &=\sum_{t=1}^T \theta_t \phi_t(x_t) - \underset{x \in X}{\text{inf}} \ \sum_{t=1}^T \theta_t \phi_t(x).
\end{align*}
We see that $\epsilon^{\bullet}(\{x_t,p_t,\theta_t\}_{t=1}^T)$ can be interpreted as a weighted regret of a sequence of convex functions $\{\phi_t(x)\}_{t=1}^T$.
\par  
For $\epsilon^{\circ}(\{x_t,p_t,\theta_t\}_{t=1}^T)$, let us assume a sequence $\{x_t\}_{t=1}^T$ such that $x_t \in X$, for all $t \in [T]$, is given. Define a function $\tilde{f}_t^i:\calP^i \rightarrow \mathbb{R}$ as $\tilde{f}_t^i(p^i) = -(p^i)^\top\bF^i(x_t)$, for all $i \in [m]$. Since $(p^i)^\top\bF^i(x_t)$ is a linear function of $p^i$, $\tilde{f}_t^i(p^i)$ is a concave function of $p^i$. Then, for each $i \in [m]$, we have 
\begin{align*}
\underset{p^i \in \calP^i}{\sup} \sum_{t=1}^T \theta_t (p^i)^\top\bF^i(x) - \sum_{t=1}^T \theta_t (p_t^i)^\top\bF^i(x) = \sum_{t=1}^T \theta_t \tilde{f}^i(p_t^i) - \underset{p^i \in \calP^i}{\inf} \ \sum_{t=1}^T \theta_t \tilde{f}_t^i(p^i).   
\end{align*}
Therefore, $\epsilon^{i}(\{x_t,p_t,\theta_t\}_{t=1}^T)$ can be interpreted as a weighted regret of a sequence of convex functions $\{\tilde{f}_t^i\}_{t=1}^T$, for all $i \in [m]$. 
\hfill \Halmos
\end{remark}

\subsection{Stochastic First-Order Meta-Algorithm for DRF} \label{subsec:general SFO}
We are now ready to introduce our stochastic first-order meta-algorithm for the DRF problem.
Similar to the work of \citet{ho2018online} in the deterministic case, we update the decision variables $(x_t)$ and the dimension parameters $(p_t)$ separately according to specific stochastic FOMs. Let us denote $\sA_x$ as the stochastic FOM used to update $x_t$ and denote $\sA_p^i$ as the stochastic FOM used to update $p^i$, for each $i \in [m]$. We use $\sA_x$ to update $x_{t+1}$ based on $\{x_s,p_s\}_{s=1}^{t}$ and $\sA_p^i$ to update $p_{t+1}^i$ based on $\{x_s,p_s \}_{s=1}^{t}$ for all $i \in [m]$. We represent these updates as
\begin{equation}
    \begin{aligned}
    x_{t+1} &= \sA_x(\{x_s,p_s\}_{s=1}^{t}) \in X, \\
    p_{t+1}^i &= \sA_i(\{x_s,p_s\}_{s=1}^{t}) \in \calP^i, \ \forall i \in [m]. \label{eqn:update process}
    \end{aligned}
\end{equation}

\begin{assumption}[High-Probability Convergence up to Precision $\bar{\epsilon}$]\label{assm:hpc}
For given sequences $\{p_t\}_{t=1}^T$ and $\{x_t\}_{t=1}^T$, and given $\bar{\epsilon} \geq 0$, algorithms $\sA_x$ and $\sA_p^i$ satisfy high-probability convergence up to precision $\bar{\epsilon}$ if there exists an increasing and sublinear function $w(\Omega): \mathbb{R}_+ \rightarrow \mathbb{R}_+$ and decreasing functions $\sR_x(T)$ and $\sR_p^i(T)$ such that
\begin{align*}
    &\mathbb{P} \big( \espbullet > w(\Omega) \sR_x(T) + \bar{\epsilon} \big) \leq 2\exp (-\Omega), \\
    &\mathbb{P} \big( \espi > w(\Omega) \sR_p^i(T)  \big) \leq 2\exp (-\Omega),
\end{align*}
for all $i \in [m]$.
\end{assumption}

\par
In Assumption \ref{assm:hpc}, both $\sR_x(T)$ and $\sR_p(T)$ are expected to converge to zero, and $\bar{\epsilon}$ is related to the desired level of tolerance. In particular, the additional error term $\bar{\epsilon}$ appears in the $\epsilon$-SMD algorithm introduced in Section \ref{subsec:epsilon SMD}, and it arises from the utilization of biased gradients within $\epsilon$-SMD. However, $\bar{\epsilon}$ is typically set to 0 in other stochastic FOMs with unbiased gradients. 
Some exemplary algorithms that satisfy the assumption of high probability convergence are stochastic gradient descent (SGD) and stochastic mirror descent (SMD). Both methods have high-probability convergence results with $R_x(T) \sim O(\frac{1}{\sqrt{T}})$, $w(\Omega) \sim O(\sqrt{\Omega})$, and $\bar{\epsilon}=0$. Under Assumption \ref{assm:hpc}, we can set $T$ large enough so that $\espbullet$ and $\espi$ are sufficiently small with high probability and, correspondingly, by Theorem \ref{thm:Nam thm 3.1} and \eqref{eqn:term cond1}, we can solve \eqref{form:DRO feas2} with high probability. 
Although it is plausible to set $T$ sufficiently large according to our theoretical results, it is usually more efficient to use the SP gap to terminate early in situations where the SP gap can be cheaply calculated periodically. Indeed, we can terminate early as soon as we find $\bar{x}$ and $\bar{p}$ that satisfy
$ \epsilon^{\phi}(\bar{x}, \bar{p}) \leq \frac{\epsilon}{2}$, where the SP gap $\epsilon^{\phi}(\bar{x}, \bar{p})$ is defined in \eqref{eqn:SP gap}. We can often calculate the SP gap in a reasonable time since $\phi(x,p)$ is a piecewise linear function of $p$ and the complexity of calculating $\epsilon^{\phi}(\bar{x},\bar{p})$ is often $O(n)$. We present the details of the SP gap calculation in Appendix \ref{appn:SP gap}. 
\par 
\begin{algorithm}[t]
    \caption{Exact Feasibility Testing Algorithm}
    \label{alg:basic FT}
    \begin{algorithmic}[1]
        \STATE \textbf{Input}: $\epsilon$, $\bar{x} \in X, \ \bar{p} \in\calP$.
        \STATE \textbf{Output}: $\epsilon$-\textit{feasibility} or \textit{infeasibility} certificate of the $\epsilon$-feasibility problem \eqref{form:DRO feas2}.
        \STATE Compute $\phi(\bar{x}, \bar{p}) =  \underset{i \in [m]}{\text{max}}  \ (\bar{p}^i)^\top\bF^i(\bar{x})$
        \IF{$\phi(\bar{x}, \bar{p}) > \frac{\epsilon}{2}$}
            \STATE \textbf{return} Infeasible
        \ELSE
            \STATE \textbf{return} $\epsilon$-feasible solution $\bar{x}$
        \ENDIF
    \end{algorithmic}
\end{algorithm}
In addition, we use the feasibility testing algorithm, denoted $\mathcal{A}_{\text{feas}}$, which yields the values of $x$ and $p$ if they constitute an $\epsilon$-feasible solution for the DRF problem; otherwise, it indicates infeasibility. Clearly, there are several algorithmic options available for $\sAfeas$. According to Theorem \ref{thm:Nam thm 3.1}, given a pair $(\bar{x}, \bar{p})$ satisfying condition $\epsilon^\phi(\bar{x}, \bar{p}) \leq \frac{1}{2} \epsilon$, we can select $\mathcal{A}_{\text{feas}}$ as an algorithm that produces an $\epsilon$-feasible solution in cases where $\phi(\bar{x}, \bar{p}) \leq \frac{1}{2} \epsilon$; otherwise, it returns infeasibility.
Algorithm \ref{alg:basic FT} represents this exact feasibility testing algorithm.
However, this can be improved in terms of computational efficiency with specific choices of $\sA_x$ and $\sA_p^i$. We discuss this further in Section \ref{subsec:efficient feas}.

\par 
Algorithm \ref{alg:SOFO} represents the stochastic first-order meta-algorithm for the DRF problem, which includes the optional early termination criterion and uses $\sAfeas$ for the feasibility testing algorithm at the very end (if we do not stop early). If we do allow the option of early termination, then the SP gap is calculated periodically at every $T_s$-th iteration.
\begin{algorithm}[t]
\caption{Stochastic Distributionally Robust Feasibility (DRF) Solver}
\label{alg:SOFO}
\begin{algorithmic}[1]
\STATE \textbf{Input}: $\sA_x$, $\sA_p^i$,$\sA_{\text{feas}}$, $\epsilon$, $T$ , $\theta \in \Delta_{T}$, optional SP gap calculation frequency $T_s$.
\STATE \textbf{Initialization}: $x_1 \in X$ and $p_1 \in \mathcal{P}$.
\STATE \textbf{Output}: $\epsilon-$\textit{feasibility} or \textit{infeasibility} certificate of the $\epsilon$-feasibility problem \eqref{form:DRO feas2}.
\FOR{$t=1,\dots,T-1$}
    \FOR{$i=1,\dots,m$}
        \STATE Update $p_{t+1}^i = \sA_p^i(\{x_s,p_s\}_{s=1}^{t}) \in \calP^i$.
    \ENDFOR
    \STATE Update $x_{t+1} = \sA_x(\{x_s,p_s\}_{s=1}^{t}) \in X$.
    \STATE \textbf{SP Gap Early Termination (Optional):}
    \IF {$t \equiv 0$ (mod $T_s$)}
    \STATE Compute $\bar{x} = \sum_{s=1}^{t+1} \frac{\theta_s}{\boldsymbol\theta} x_s$, $\bar{p} = \sum_{s=1}^{t+1} \frac{\theta_s}{\boldsymbol\theta} p_s$, where $\boldsymbol\theta = \sum_{k=1}^{t+1} \theta_k$, and compute $\epsilon^{\phi}(\bar{x},\bar{p})$.
    \IF { $\epsilon^{\phi}(\bar{x}, \bar{p}) \leq \frac{\epsilon}{2}$}
             \IF{$\phi(\bar{x}, \bar{p}) > \frac{\epsilon}{2}$}
            \STATE \textbf{return} Infeasible
        \ELSE
            \STATE \textbf{return} $\epsilon$-feasible solution $\bar{x}$
        \ENDIF
    \ENDIF
    \ENDIF
\ENDFOR
\STATE Compute $\bar{x}_T = \sum_{s=1}^{T} \theta_s x_s$ and $\bar{p}_T = \sum_{s=1}^{T} \theta_s p_s$. \\
\textbf{return} $\mathcal{A}_{\text{feas}}(\epsilon,\bar{x}_T, \bar{p}_T)$

\end{algorithmic}
\end{algorithm}
The following theorem shows the convergence guarantee of Algorithm \ref{alg:SOFO}.
\begin{theorem}\label{thm:high prob m>1}
Let $\{x_t,p_t\}_{t=1}^T$ be a sequence generated by $\sA_x$ and $\sA_p^i$, for all $i \in [m]$, according to the update rules \eqref{eqn:update process}, and suppose that Assumption \ref{assm:hpc} holds. Then, for any $\Omega \geq 0$, Algorithm \ref{alg:SOFO} satisfies:
\begin{align*}
    \mathbb{P}\Big\{  \espbullet + \espcirc &\leq w(\Omega)\left(\sR_x(T) +  \underset{i \in [m]}{\max} \ \sR_p^i(T)\right) + \bar{\epsilon} \Big\} \\
    &\geq 1 - 4m \exp (- \Omega). 
\end{align*} 
\end{theorem}

Theorem \ref{thm:high prob m>1} translates the high-probability convergence guarantees of the basic algorithms (Assumption \ref{assm:hpc}) to Algorithm \ref{alg:SOFO} and implies if $\Omega \geq \ln (\frac{4m}{\nu_0})$, then 
\begin{align}
    \mathbb{P}\Big\{  &\epsilon^{\bullet}(\{x_t,p_t,\theta_t\}_{t=1}^T) + \espcirc \leq  w(\Omega)(\sR_x(T) + \underset{i \in [m]}{\max} \ \sR_p^i(T)) +\bar{\epsilon} \Big\} \geq 1- \nu_0. \label{eqn:Omega bound}
\end{align}
Equation \eqref{eqn:Omega bound} implies that $\Omega$ has low dependency on the number of constraints $m$ and the parameter $\nu_0.$

%% file: FirstOrderMethods.tex
In this section, we introduce two specific basic stochastic FOMs to improve the scalability and applicability of Algorithm \ref{alg:SOFO}. Namely, we introduce a new algorithm, $\epsilon$-Stochastic Mirror Descent ($\epsilon$-SMD), to update the $x_t$ variables (Section \ref{subsec:epsilon SMD}). $\epsilon$-SMD not only improves the per-iteration cost associated with updating $x_t$, but also offers a more efficient approach to design the feasibility test $\mathcal{A}_{\text{feas}}$ (Section \ref{subsec:efficient feas}). We also introduce Bandit Mirror Descent (BMD) \citep{namkoong2016stochastic} to update the $p_t$ variables (Section \ref{subsec:BMD}), and we develop a new high probability convergence guarantee for this method, which is needed in our setting to apply Theroem \ref{thm:high prob m>1}.
\par 
\subsection{$\epsilon$-Stochastic Mirror Descent}\label{subsec:epsilon SMD}
As mentioned earlier, $\espbullet$ can be interpreted as a weighted regret for $\{x_t\}_{t=1}^T$, suggesting the use of online convex optimization methods to perform the updates. Note that these methods typically use unbiased stochastic gradients. However, we point out that using an online unbiased stochastic FOM for $\sA_x$ could cause a substantial computational burden when $m$ and $n$ are large.
Note that in \eqref{eqn:update process}, the update of $x_t$ and $p_t$ does not occur simultaneously. As a result, we consider a sequence $\{p_t\}_{t=1}^T$ as given when we discuss the update of $x_t$ and vice versa. As we only discuss the update of $x_t$ in this section, we simplify our notation as follows throughout this section:  we define $f_t^i: X \rightarrow \mathbb{R}$  as $f_t^i(x) := (p_t^i)^\top\bF^i(x)$ and $\phi_t:X \rightarrow \mathbb{R}$ as $\phi_t(x) := \underset{i \in [m]}{\max} \ f_t^i(x)$.
\par 
The main inefficiency of using an unbiased stochastic FOM comes from the structure of the function $\phi_t(x) = \underset{i \in [m]}{\max} \ f_t^i(x) $.
In order to obtain an unbiased estimator of $\phi_t(x)$ at $x_t$ and associated first-order information, we need to calculate an index $i_t^* \in \underset{i \in [m]}{\text{argmax}} \ f_t^i(x_t)$.
However, as $f_t^i(x_t) = (p_t^i)^\top\bF^i(x_t)$, for all $i \in [m]$, the complexity of obtaining $i_t^*$ is $O(mn\costF)$, where $\costF$ is the average cost to compute $F_r^i(x)$, for all $r \in [n]$, $i \in [m]$ and $x \in X$. Thus, calculating $i_t^*$ dominates the other steps of the stochastic FOMs and causes Algorithm \ref{alg:SOFO} to lose scalability in $n$. 
\par 
To overcome this issue, we use an approximate maximizing index $\hat{i}_t$ instead of an exact maximizing index $i_t^*$.
Algorithm \ref{alg:i hat} describes a procedure for obtaining this $\hat{i}_t$. We use $K$ samples of indices $s_1,\cdots,s_K$ from distribution $p_t^i$ to obtain $\hat{f}_t^i(x_t):= \frac{1}{K}\sum\limits_{k=1}^K F_{s_k}^i(x_t)$, which approximates $f_t^i(x_t)$ for all $i \in [m]$. We call $K$ the ``per-iteration sample size'' to distinguish it from a batch size for stochastic gradients in a stochastic FOM. 

 \begin{algorithm}[t]
\caption{Procedure for Obtaining $\hat{i}$ ($\sA_{\hat{i}}(x,p,K)$)}
\begin{algorithmic}[1]
\STATE \textbf{input}: per-iteration sample size $K$, $x$ and $p$
\STATE \textbf{output}: $\hat{i}$
\FOR{$i=1,\dots,m$}
        \STATE Sample i.i.d $K$ indices $s_{1},\cdots, s_K \in [n]$ from  distribution $p^i$ with replacement.
\ENDFOR
\STATE Compute $\hat{i} \in \underset{i \in [m]}{\text{argmax}} \  \{  \hat{f}^i(x):= \frac{1}{K} \sum\limits_{k=1}^K F_{s_k}^i(x)\}$
\STATE \textbf{return} $\hat{i}$
\end{algorithmic}
\label{alg:i hat}
\end{algorithm}

Lemma \ref{lem:main lemma} below shows that if the sample size $K$ is large enough to satisfy the condition in Definition \ref{def:approx}, then $\hat{i}_t$ is a good approximation of $i_t^*$ in the sense that the approximate subgradient set $\partial \hat{\phi}_{t}(x):= \partial \left(\underset{i \in [m]}{\max} \ \hat{f}_{t}^{i}(x)\right)$ is close to the actual subgradient set $\partial \phi_t(x)$. To explain this formally, we introduce the concept of $\bar\epsilon$-subgradients (see, e.g., \citet{borwein1982note}, \citet{kiwiel2004convergence} and \citet{bonettini2016scaling} for more details).
 \begin{definition}[$\bar\epsilon$-subgradient]
For any given convex function $f: X \rightarrow \mathbb{R}$ and $\bar{\epsilon}>0$, where $X \subseteq \mathbb{R}^n$ is a convex set, a vector $g \in \mathbb{R}^n$ is called an $\bar{\epsilon}$-subgradient of $f$ at $\bar{x}$ if 
\begin{align*}
    f(x) \geq f(\bar{x}) + g^\top(x-\bar{x}) - \bar{\epsilon}, \quad \forall x \in X.
\end{align*}
We denote the set of all $\bar{\epsilon}$-subgradients of $f$ at $\bar{x}$ by $\partial_{\bar{\epsilon}}f(\bar{x})$. \hfill \Halmos
\end{definition}
\begin{definition}[Approximation Condition]\label{def:approx}
    For given $\epsilon >0$, the per-iteration sample size $K$ satisfies the approximation condition, with parameters $C_K \in [0,\frac{1}{2})$ and $\nu_1 \in (0,1)$, if it holds that
\begin{align}
    |\hat{f}_{t}^{i}(x)-f_t^{i}(x)|\leq \frac{C_K \epsilon}{2}, \ \forall x \in X, i \in [m], t \in [T]. \label{eqn:K cond}
\end{align}
with probability at least $1- \nu_1$. \hfill \Halmos
\end{definition}

Notice that both $C_K$ and $\nu_1$ are fixed, user-specified, parameters that are independent of any other parameters in our context. We now introduce the key lemma for $\epsilon$-SMD.
\begin{lemma}\label{lem:main lemma}
For given $\epsilon >0$, suppose that the per-iteration sample size $K$ satisfies the approximation condition with parameters $(C_K, \nu_1)$. Let $\bar{\epsilon}$ be defined as $\bar{\epsilon} = C_K\epsilon$ and $x \in X$ be given. Then for any $h_t \in \partial f_t^{\hat{i}_t}(x)$, $h_t$ is an $\bar{\epsilon}$-subgradient of $\phi_t(x)$, i.e., $h_t \in \partial_{\bar{\epsilon}} \phi_t(x)$, with probability at least $1- \nu_1$.
\end{lemma}

\begin{remark}[Value of Per-iteration Sample Size K]
    As $\hat{f}_{t}^{i}(x)$ is an unbiased estimator of $f_t^{i}(x)$, it is possible to choose $K$ that satisfies the approximation condition for small $\nu_1$ by the law of large numbers. However, $K$ being large is not desirable since the complexity of obtaining $\hat{i}$ is $O(mK \costF)$. Intuitively, choosing $K$ on the order of $n$ is overly conservative to satisfy the approximation condition. In fact, we can satisfy this condition with $K \ll n$ since, in Appendix \ref{appn:lower bound}, we show that $K$ is almost independent of $n$ when $n$ is large. For example, using Bennet's inequality \citep{wainwright2019high}, we can choose $K = O\left(\frac{\log (1/\nu_1)}{C_K\epsilon \log(1+C_K\epsilon)}\right)$ to satisfy the approximation condition for any $C_K$ and $\nu_1$. This value is independent of $n$ and much smaller than $n$ if $n$ is large. \hfill \Halmos
\end{remark}

\par 
Our main basic stochastic FOM to update $x_t$, $\epsilon$-SMD, is a variant of SMD. So, we assume that the set $X$ follows the following mirror descent setup \citep{juditsky2011first}.
\begin{assumption}[Mirror Descent Setup]\label{assm:mds}
Let $\psi_x: X \rightarrow \mathbb{R}$ be a given differentiable and 1-strongly convex function and $V_X$ be the Euclidean space containing $X$. We define the Bregman distance $B_{\psi_x}: X \times X \rightarrow \mathbb{R}$ as $B_{\psi_x}(x,y):= \psi_x(x) - \psi_x(y) - \nabla \psi_x(y)^\top (x-y)$ and the $\psi_x$-diameter $D_x:= \underset{x,y \in X}{\max} \ B_{\psi_x}(x,y) < \infty$. Then for any given $w \in X$, $g \in V_X^*$ and $\alpha>0$, it is easy-to-compute $\underset{x \in X}{\mathrm{argmin}} \ \big\{ g^\top x + \frac{1}{\alpha}B_{\psi_x}(x,w)  \big\}$.
\end{assumption}
\par 
\begin{algorithm}[t]
\caption{$\epsilon$-Stochastic Mirror Descent ($\epsilon$-SMD) Update Rule}
\begin{algorithmic}[1]
\STATE \textbf{input}: step-size $\alpha_{x,t}$, sample size $K$, $x_{t}$ and $p_t$, $\sAhat$
\STATE \textbf{output}: $x_{t+1}$
\STATE Compute $\hat{i}_t = \sAhat\big(x_t,p_t,K \big)$
\STATE Sample an index $I_{x,t} \in [n]$ from distribution $p_t^{\hat{i}_t}$. Then compute $g_{x,t}$ such that  $g_{x,t} \in \partial F_{I_{x,t}}^{\hat{i}_t}(x_t)$. \label{line:alg2 stoc sub}
\STATE $x_{t+1} \leftarrow \underset{x \in X}{\text{argmin}} \ \big\{ g_{x,t}^{\top}(x- x_t) + \frac{1}{\alpha_{x,t}}B_{\psi_x}(x,x_t) \big\}$ \label{line:SMD}
\STATE \textbf{return} $x_{t+1}$
\end{algorithmic}
\label{alg:Stoc Ax}
\end{algorithm}

We introduce our main online stochastic FOM for updating the $x_t$ variables, which we call $\epsilon$-SMD, in Algorithm \ref{alg:Stoc Ax}.
The main difference from unbiased stochastic FOMs such as SGD or SMD is that $\epsilon$-SMD uses a biased subgradient (which approximates an unbiased subgradient) of the objective. 
A desirable property of SMD is its high-probability convergence guarantee \citep{nemirovski2009robust}.
We now present such a high-probability convergence result for $\epsilon$-SMD. In addition to the notation defined in Assumption \ref{assm:mds}, which will be used in Theorem \ref{thm:x high prob conv}, recall that $G$ denotes the uniform Lipschitz constant of functions $F^i(\cdot,z^i)$ for all $z^i \in Z^i$ and $i \in \{1, \dots,m\}$.

\begin{theorem} \label{thm:x high prob conv}
Let $w(\Omega)$ be defined as $w(\Omega)= 3\sqrt{\Omega}$ for $\Omega \geq 1$ and $\epsilon>0$ be given.
Given a sequence $\{p_t\}_{t=1}^T$ such that $p_t \in \calP$ for all $t \in [T]$, let $\{x_t\}_{t=1}^T$ be the sequence generated by Algorithm \ref{alg:Stoc Ax} with diminishing step size $\alpha_{x,t} = \frac{c_x}{\sqrt{t}}$, where the constant $c_x$ is defined as $c_x= \frac{1}{G}\sqrt{\frac{D_x}{\Omega}}$. Moreover, let $\theta \in \Delta_T$ be the sequence derived by normalizing $\{\alpha_{x,t}\}_{t=1}^T$. Given that the per-iteration sample size sequence $K$ satisfies Equation \eqref{eqn:K cond} in the approximation condition (Definition \ref{def:approx}), for any $\Omega \geq 1$, we have
 \begin{align}
     \mathbb{P}\Big \{&\epsilon^{\bullet}(\{x_t,p_t,\theta_t\}_{t=1}^T) \leq \frac{ w(\Omega)\cdot G \cdot \ln(T)\cdot \sqrt{2 D_x }}{\sqrt{T}}+ C_K \epsilon ~|~ \eqref{eqn:K cond} \text{ holds}\Big\}
     \geq 1 - 2\exp (- \Omega).\label{eqn:nu-naive}
 \end{align}
Otherwise, if Equation \eqref{eqn:K cond} holds with probability at least $1 - \nu_1$, we have
 \begin{align}
     \mathbb{P}\Big \{&\epsilon^{\bullet}(\{x_t,p_t,\theta_t\}_{t=1}^T) \leq \frac{ w(\Omega)\cdot G \cdot \ln(T)\cdot \sqrt{2 D_x }}{\sqrt{T}}+ C_K \epsilon \Big\}
     \geq 1 - 2\exp (- \Omega)-\nu_1. \label{eqn:nu}
 \end{align}
\end{theorem}
Note that we choose the constant $c_x$ of the diminishing step-size in the above theorem to minimize the upper bound of $\espbullet$ in the statement of the Theorem. Also, the added $\nu_1$ term in equation \eqref{eqn:nu} reflects the probability of satisfying equation \eqref{eqn:K cond} under the approximation condition.
In Equation \eqref{eqn:nu-naive}, we observe that $\epsilon$-SMD satisfies Assumption \ref{assm:hpc} with $\sR_x(T)=\frac{G \cdot \ln(T)\cdot \sqrt{2 D_x }}{\sqrt{T}}$, $\bar{\epsilon}=C_K\epsilon$, and $w(\Omega) = 3\sqrt{\Omega}$.
Theorem \ref{thm:x high prob conv} tells us that, relative to SMD, $\epsilon$-SMD achieves a low cost per iteration of $O(mK\costF)$ at the cost of an increase in tolerance error by $C_K \epsilon$. 
Furthermore, we can see that the upper bound of $\espbullet$ does not depend directly on the dimensions $n$ and $d$, and often has a very low dependence on these parameters. For example, when $x \in \Delta_d$, and $\psi_x$ is the entropy function, we have $D_x = \log(d)$ which indicates the low dependence of the upper bound of $\espbullet$ on $d$. We emphasize that $\epsilon$-SMD is independent of the choice of the convex ambiguity set. In addition to the reduced per iteration cost offered by $\epsilon$-SMD, the information of $\{\hat{f}_t^i(x_t)\}_{t=1}^T$ for all $i \in [m]$, which we obtain from Algorithm \ref{alg:i hat}, can be reused to develop an efficient feasibility testing algorithm, which is elaborated on in Section \ref{subsec:efficient feas}.

\begin{remark}[Practical Performance with Diminishing Step Size]
The upper bound of $\espbullet$ in Theorem \ref{thm:x high prob conv} contains a logarithmic term $\ln (T)$ that can be removed if we use a constant step size. However, in practice, using a diminishing step size typically leads to a higher-quality solution with fewer iterations compared to a constant step size. As a result, to optimize the advantages of early termination, we choose the diminishing step size, despite its theoretical limitations. \hfill \Halmos
\end{remark}

\subsection{Efficient Feasibility Test}\label{subsec:efficient feas}
In Algorithm \ref{alg:basic FT}, we determine the $\epsilon$-feasibility or infeasiblity of a given solution $(\bar{x},\bar{p})$, satisfying $\epsilon^\phi(\bar{x},\bar{p}) \leq \frac{1}{2}\epsilon$, by comparing $\phi(\bar{x},\bar{p})$ with $\frac{1}{2}\epsilon$. Given that $\epsilon^\phi(\bar{x},\bar{p}) \leq \frac{1}{2}\epsilon$, which can satisfied by taking $T$ sufficiently large in Algorithm \ref{alg:SOFO}, the complexity of Algorithm \ref{alg:basic FT} is dominated by computing $\phi(\bar{x},\bar{p}) = \underset{i \in [m]}{\max}\ \left(\bar{p}^i\right)^\top \bF ^i(\bar{x})$, which is $O(mn\costF)$. Since $n$ may be very large in settings with high dimensional distribution parameter vectors, the linear dependence of this complexity on $n$ is problematic.
To confront this challenge, we introduce an \textit{efficient feasibility test}.
This test, when used in conjunction with $\epsilon$-SMD as the updating algorithm for $x$, streamlines the feasibility testing process and completely removes the linear dependence on $n$. By cleverly reusing the information $\hat{f}_t^i(x_t)$ that was obtained during the application of $\epsilon$-SMD, the new efficient feasibility test has computational complexity $O(m\tilde{T})$, where $\tilde{T} \approx T$ represents a new value for the total number of iterations in Algorithm \ref{alg:SOFO}.
Importantly, we show that $\tilde{T}$ is within a constant factor of the value of $T$ required when we instead use the exact feasibility test Algorithm \ref{alg:basic FT}.
As we will show in Section \ref{subsec:Meta-chi}, $T$ and $\tilde{T}$ are both independent of $n$.

Towards the goal of developing the efficient feasibility test, we introduce a result below that provides an alternative way of producing a certificate of the DRF problem.
\begin{proposition}\label{prop:Nam 3.2} \textnormal{\textbf{(Extension of Corollary 3.2 of \cite{ho2018online})}}
Let $\{ x_t, p_t, \theta_t \}_{t=1}^T$ be a sequence with $x_t \in X, p_t \in \mathcal{P}$, for all $t \in [T],$ and $\theta \in \Delta_T$ and $\epsilon>0$ be given. Suppose that there exists $\nu_1, \kappa^{\bullet}, \kappa^{\circ} \in (0,1)$ satisfying $\espcirc \leq \epsilon \kappa^{\circ}$, $\espbullet \leq \epsilon \kappa^{\bullet}$, and $\kappa^{\bullet} + \kappa^{\circ} \leq 1$, with probability at least $1- \nu_1$. Let $\tau \in [\kappa^{\circ}, 1- \kappa^{\bullet}]$ and $\vartheta(\{x_t,p_t,\theta_t\}_{t=1}^T):=\underset{i \in [m]}{\max}\sum_{t=1}^T \theta_t (p_t^i)^\top\bF^i(x_t)$. With probability at least $1-\nu_1$,
if $\vartheta(\{x_t,p_t,\theta_t\}_{t=1}^T) \leq (1-\tau)\epsilon$, then $\bar{x}_T$ is an $\epsilon$-feasible certificate, and if $\vartheta(\{x_t,p_t,\theta_t\}_{t=1}^T) > (1-\tau)\epsilon$, then problem \eqref{form:DRO feas2} is infeasible. 
\end{proposition}
Following Proposition \ref{prop:Nam 3.2}, if we choose $T$ that satisfies 
\begin{align}
    w(\Omega)\left(\sR_x(T) + \underset{i \in [m]}{\max} ~ \sR_p^i(T)\right) +C_K\epsilon \leq \epsilon, \label{eqn:T2}
\end{align}
 and set $\kappa^\bullet$ as $\kappa^\bullet = \frac{w(\Omega) \sR_x(T)}{\epsilon} + C_K$ and set $\kappa^\circ = \underset{i \in [m]}{\max} ~  \sR_p^i(T)/ \epsilon$, then we can set $\tau = 1-\kappa^\bullet$ and produce a certificate for problem \eqref{form:DRO feas2} by comparing $\vartheta(\{x_t,p_t,\theta_t\}_{t=1}^T)$ with $\kappa^\bullet \epsilon$.
Unfortunately, the complexity of computing $\vartheta(\{x_t,p_t,\theta_t\}_{t=1}^T)$ is $O(mnT\costF)$, which actually exceeds the complexity of determining $\phi(\bar{x}_T,\bar{p}_T)$, $O(mn\costF)$, as required by Algorithm \ref{alg:basic FT}. 
Nevertheless, using a stochastic approximation of $\vartheta(\{x_t,p_t,\theta_t\}_{t=1}^T)$, which we obtain almost for free by reusing the information of $\hat{f}_t^i(x_t)$, allows us to develop the efficient feasibility test. Indeed, notice that $\hat{f}_t^i(x_t)$ is the approximate value of $f_t^i(x_t)$ that satisfies \eqref{eqn:K cond}.
In this case, we assume $C_K \in (0,\frac{1}{3})$ and adjust the total number of iterations $\tilde{T} \geq T$ to satisfy
\begin{align}
    w(\Omega)\left(\sR_x(\tilde{T}) + \underset{i \in [m]}{\max} \ \sR_p^i(\tilde{T}) \right) + C_K\epsilon  \leq (1-2C_K)\epsilon. \label{eqn:bound sum of kappa}
\end{align}
For clarity, we simplify our notation and denote $\vartheta:=\vartheta(\{x_t,p_t,\theta_t\}_{t=1}^{\tilde{T}})$ and $\vartheta^i := \sum_{t=1}^{{\tilde{T}}} \theta_t (p_t^i)^\top\bF^i(x_t)$ to simply express $\vartheta$ as $\vartheta  = \underset{i \in [m]}{\max} \ \vartheta^i$. We approximate $\vartheta^i$ with $\hvt^i$, which we define as $\hvt^i := \sum_{t=1}^{{\tilde{T}}} \theta_t \hat{f}_t^i(x_t)$, for all $i \in [m]$ and further define $\hvt$ as $\hvt:= \underset{i \in [m]}{\max} \ \hvt^i$.
\begin{lemma}\label{lemma:vartheta}
    Let $\epsilon>0$ be given and $\theta \in \Delta_{\tilde{T}}$ and $\{x_t,p_t\}_{t=1}^{\tilde{T}}$ be a sequence that is generated by Algorithm \ref{alg:SOFO} with $\epsilon$-SMD as $\sA_x$. If the per-iteration sample size $K$ satisfies the approximation condition with parameters $C_K\in (0,\frac{1}{3})$ and $\nu_1 \in (0,1)$, then we have
\begin{align}
    |\vartheta - \hvt| \leq C_K\epsilon. \label{eqn:vartheta and hvt diff}
\end{align}
with probability at least $1- \nu_1$.
\end{lemma}

Lemma \ref{lemma:vartheta} suggests that $\hvt$ is sufficiently close to $\vartheta$ to serve as a substitute for $\vartheta$. Additionally, the complexity of calculating $\hvt$ is $O(m\tilde{T})$, as it reuses the values of $\hat{f}_t^i(x_t)$ that were previously calculated. Because we use $\hvt$ instead of $\vartheta$, we need a different guideline to determine $\epsilon$-feasibility other than comparing $\vartheta$ and $\kbullet \epsilon$.
Instead, we compare $\hvt$ with $(\kbullet+ C_K)\epsilon$, and Algorithm \ref{alg:eff FT} presents the efficient feasibility testing algorithm using this stochastic approximation.
Notice that Algorithm \ref{alg:eff FT} can be used after a sufficient number of iterations $\tilde{T}$ that satisfies \eqref{eqn:bound sum of kappa}.
The following theorem implies that after $\tilde{T}$ iterations, Algorithm \ref{alg:eff FT} gives us a correct solution to \eqref{form:DRO feas2} with high probability.

\begin{algorithm}[t]
    \caption{Efficient Feasibility Testing Algorithm Using Stochastic Approximation}
    \label{alg:eff FT}
    \begin{algorithmic}[1]
        \STATE \textbf{Input}: $\epsilon$, $\theta \in \Delta_T$, $\kbullet$, $\bar{x}_T$, $C_K$, $\{\hat{f}_s^i(x_t)\}_{t=1}^T$, for all $i \in [m]$.
        \STATE \textbf{Output}: $\epsilon-$\textit{feasibility} or \textit{infeasibility} certificate of the $\epsilon$-feasibility problem \eqref{form:DRO feas2}.
        \STATE Compute $\hvt=\underset{i \in [m]}{\max} \ \sum_{t=1}^T \theta_t \hat{f}_t^i(x_t)$
        \IF{$\hvt > (\kbullet + C_K)\epsilon $}
            \STATE \textbf{return} Infeasible
        \ELSE
            \STATE \textbf{return} $\epsilon$-feasible solution $\bar{x}_T$
        \ENDIF
    \end{algorithmic}
\end{algorithm}

\begin{theorem}\label{thm:feas}
Let $\epsilon>0$ be given, and $\{x_t,p_t\}_{t=1}^{\tilde{T}}$ be a sequence generated by Algorithm \ref{alg:SOFO} with $\epsilon$-SMD as $\sA_x$ and the per-iteration sample size $K$ satisfies the approximation condition with parameter $C_K \in (0,\frac{1}{3})$ and $ \nu_1 \in (0,1)$.
Suppose $\tilde{T}$ satisfies \eqref{eqn:bound sum of kappa}, and let $\kbullet = \sR_x(\tilde{T})/\epsilon + C_K$.
 Then, Algorithm \ref{alg:eff FT} returns either an $\epsilon$-feasible solution or an infeasibility certificate of the DRF problem \eqref{form:DRO feas2}, with probability at least $1-\nu_1$.
\end{theorem}

The remaining task is to contrast $\tilde{T}$, which represents the total number of iterations required when using $\hvt$ to determine the $\epsilon$-feasibility of the problem, with $T$, which is the total number of iterations when using Algorithm \ref{alg:basic FT} for the same purpose. Notice that $\tilde{T}$ needs to satisfy \eqref{eqn:bound sum of kappa}, while $T$ needs to satisfy 
\begin{align}
    w(\Omega)\left(\sR_x(T) + \underset{i \in [m]}{\max} ~ \sR_p^i(T)\right) +C_K\epsilon \leq \frac{\epsilon}{2}.\label{eqn:T Theorem1}
\end{align}
In fact, considering $\sR_x(T) \sim O(\frac{1}{\sqrt{T}})$, the relationship between $T$ and $\tilde{T}$ is $\frac{T}{\tilde{T}} \sim O(\big(\frac{2(1 - 3C_K)}{1-2C_K}\big)^2)$. In particular, for all $C_K \in (0, 1/3)$, we have that $T$ and $\tilde{T}$ will be within a constant factor of each other. In fact, if $C_K \in (0, 1/4)$, then $\tilde{T} \leq T$, although they are still within a constant factor of each other. In any case, we see that using the efficient feasibility test Algorithm \ref{alg:eff FT} reduces the complexity of the $\sAfeas$ subroutine of Algorithm \ref{alg:SOFO}, and sometimes requires less itearations in total. This reduction is highly beneficial when early termination (SP gap termination) in Algorithm \ref{alg:SOFO} is unattainable and we ultimately need to rely on $\sAfeas$ to solve the problem.

\subsection{Bandit Mirror Descent}\label{subsec:BMD}
Let us now introduce Bandit Mirror Descent (BMD) \citep{namkoong2016stochastic}, which provides an efficient update scheme for the variables $p_t$ in the case where the ambiguity set $\calP^i$ are constructed from $\chi^2$-divergence sets. The $\chi^2$-divergence set is defined as
$\mathcal{P}_{n,\rho}:= \{p \in \mathbb{R}_+^n ~|~ \mathbbm{1}_n^{\top}p = 1, D_{\chi^2}(p\| \frac{\mathbbm{1}_n}{n} )\leq \frac{\rho}{n}\}$ with parameter $\rho \in (0,1)$, where $D_{\chi^2}(p\|q)$ represents the $\chi^2$-divergence between distributions $p$ and $q$.
In particular, if the distributions $p$ and $q$ have finite support, we have $ D_{\chi^2}(p\|q)= \sum_{r=1}^n  \frac{q^{(r)}}{2} (\frac{p^{(r)}}{q^{(r)}}-1)^2$. We consider $\pset :=\{p \in \mathbb{R}^n ~|~ p\geq \frac{\delta}{n}\mathbbm{1}_n, D_{\chi^2}(p\| \frac{\mathbbm{1}_n}{n} )\leq \frac{\rho}{n}\}$ with parameters $\rho,\delta \in (0,1)$ to be our ambiguity set, which is a slight variation of the set $\calP_{n,\rho}$. This modification of the ambiguity set is related to BMD as we discuss later in this section. Notice that $p \in \pset$ is bounded below with parameter $\delta \in (0,1)$ and $p\in \pset$ no longer needs to satisfy $\mathbbm{1}_n^\top p =1$. Instead, we include a normalization step when obtaining a stochastic gradients of $(p^i)^\top \bF^i(x)$ at $x=x_t$ or $p^i = p_t^i$.
\par

\par 
BMD is a variant of SMD, so it follows the mirror descent setup specified in Assumption \ref{assm:mds} with 
$\psi_p(p)=\frac{1}{2}\|p\|_2^2$ and $\psi_p$-diameter $D_p = \frac{4\rho}{n^2}$. The derivation for $D_p$ can be found in Appendix \ref{appn:Dp bound}.
Algorithm \ref{alg:p under chi} shows the Bandit Mirror Descent under $\chi^2$-divergence set. 
The main reason for the choice of $\pset$ instead of $\mathcal{P}_{n,\rho}$ is because under $\mathcal{P}_{n,\rho}$, the individual probabilities $p_t^{i^{(r)}}$ could be small, causing the gradient estimator $g_{p,t}^i$ to have an excess variance. \citet{namkoong2016stochastic} showed that the convergence of the algorithm still holds even if we use $\pset$ instead of $\mathcal{P}_{n,\rho}$. We further extend their convergence result to a high-probability convergence result. Before presenting the upcoming proposition, let us recall the definition of $M^i$. $M^i$ is a constant that uniformly bounds $|F^i(x,z^i)|$ for all $x \in X$ and $z^i \in Z^i$ for $i \in [m]$. With this in mind, we present the following proposition.
\begin{proposition} \label{prop:Rp with entropy}
Let $C_g = 1 + \sqrt{\frac{2\rho}{n}}$ and $w(\Omega) = 3\sqrt{\Omega}$ for $\Omega \geq 1$.
Given a sequence $\{x_t\}_{t=1}^T$ where $x_t \in X$  for all $t \in [T]$, let $\{p_t^i\}_{t=1}^T$ be the sequence generated from Algorithm \ref{alg:p under chi} using a diminishing step-size
$\alpha_{p,t}=\frac{c_p}{\sqrt{t}}$, where the constant $c_p$ is defined as 
\begin{align}
   c_p= \frac{2\delta}{C_g\cdot M^i\cdot n^2}\cdot \sqrt{\frac{\rho}{\Omega}}, \label{eqn:p step size}
\end{align}
and $\theta \in \Delta_T$ be the normalization of $\{\alpha_{p,t}\}_{t=1}^T$.
For any $\Omega \geq 1$, we have
\begin{align*}
    \mathbb{P}\Big\{  &\epsilon^{i}(\{x_t,p_t^i,\theta_t\}_{t=1}^T) \leq \frac{w(\Omega) \cdot 2C_g \cdot M^i\cdot\ln (T)\cdot\sqrt{2\rho}}{\delta\sqrt{T}} \Big\}\leq 1 - 2\exp (- \Omega), \ \forall i \in [m].
\end{align*}
\end{proposition}
Note that the upper bound of $\espi$ does not directly depend on $n$, and $\sR_p(T) = \frac{2C_g \cdot M^i\cdot\ln (T)\cdot\sqrt{2\rho}}{\delta\sqrt{T}}$ satisfies Assumption \ref{assm:hpc} with $w(\Omega) = 3\sqrt{\Omega}$. Also, following the work of \citet{namkoong2016stochastic}, we propose an efficient way of updating the $p_t$ variables in BMD in Appendix \ref{appn:efficient update}. In particular, we show that the complexity of the update of the $p_t$ variables at each iteration of BMD is $O(\log(n))$. These results show the scalability of BMD to update $p_t$ variables.
\begin{algorithm}[t]
\begin{algorithmic}[1]
\STATE \textbf{input} $x_t$ and $p_t^i$, stepsize $\alpha_{p,t}^i$, $\psi_p(p) = \frac{1}{2}\|p\|_2^2$
\STATE \textbf{output} $p_{t+1}^i$ 
\STATE Sample $I_{p,t}^i \sim \frac{p_{t}^i}{\mathbbm{1}_n^{\top}p_{t}^i}$
\STATE Compute gradient estimator $g_{p,t}^i$ such that $g_{p,t}^{i^{(r)}} = \frac{(\mathbbm{1}_n^{\top}p_t^i)F_r^i(x_t)}{p_t^{i^{(r)}}} \mathbbm{1}(r = I_{p,t}^i)$ for all $r \in [n]$
\STATE Update $w_{t+1}^i = p_t^i + \alpha_{p,t}^i g_{p,t}^i$
\STATE Perform projection: $p_{t+1}^i = \underset{p^i\in \mathcal{P}_{\rho,n,\delta}}{\text{argmin }} B_{\psi_p}(p^i,w_{t+1}^i)$

\STATE \textbf{return} $p_{t+1}^i$

\end{algorithmic}
\caption{Bandit Mirror Descent Update Rule (\citet{namkoong2016stochastic})}
\label{alg:p under chi}
\end{algorithm}
\par 
\subsection{Meta-Algorithm under $\chi^2$-divergence set}\label{subsec:Meta-chi}

Now we are ready to consider Algorithm \ref{alg:SOFO}, in the special case where we use Algorithm \ref{alg:Stoc Ax} to update the $x_t$ variables and use Algorithm \ref{alg:p under chi} to update the $p_t^i$ variables. The following result shows its convergence and follows by combining Theorem \ref{thm:x high prob conv} and Proposition \ref{prop:Rp with entropy} together.
\begin{theorem} \label{thm:alg1 conv}
Let $C_g = 1 + \sqrt{\frac{2\rho}{n}}$ and $w(\Omega) = 3\sqrt{\Omega}$ for $\Omega \geq 1$, and $\epsilon>0$ be given.
Suppose that $K$ satisfies the approximation condition (Definition \ref{def:approx}) with $C_K \in (0,\frac{1}{2})$ and $\nu_1 \in (0,1)$.
Consider the sequence $\{x_t, p_t\}_{t=1}^T$ generated by Algorithm \ref{alg:SOFO}, where the $x_t$ variables are generated using $\epsilon$-SMD with a diminishing step-size of $\frac{c_x}{\sqrt{t}}$ with $c_x = \frac{1}{C_g\cdot G}\sqrt{\frac{D_x}{\Omega}}$, and the $p_t$ variables are generated using BMD with a diminishing step-size of $\frac{c_p}{\sqrt{t}}$ with $c_p$ as specified in \eqref{eqn:p step size}.
Let $\theta \in \Delta_T$ be the normalization of $\{\alpha_{x,t}\}_{t=1}^T$, which is also the normalization of $\{\alpha_{p,t}\}_{t=1}^T$.
 For any $\Omega \geq 1$, let
\begin{align}
    &\sR_x(T) = \frac{C_g\cdot G\cdot \ln(T)\cdot\sqrt{2 D_x}}{\sqrt{T}},\ \sR_p^i(T) = \frac{2C_g\cdot M^i\cdot \ln(T)\cdot\sqrt{2\rho }}{\delta\sqrt{T}}, \ \forall i \in [m],\label{eqn:R_x and R_p}
\end{align}
then we have 
\begin{align*}
    \mathbb{P}\Big\{  \espbullet + \espcirc &\leq  w(\Omega)\left(\sR_x(T) + \underset{i \in [m]}{\max} \ \sR_p^i(T)\right) +C_K\epsilon \Big\}\\
    &\geq
    1 - 4m \exp (- \Omega)- \nu_1.   
\end{align*}
\end{theorem}
\begin{remark}[Scalability of Algorithm \ref{alg:SOFO} under $\chi^2$-divergence set]

The total number of iterations $T$ of the version of Algorithm \ref{alg:SOFO} considered in Theorem \ref{thm:alg1 conv} must satisfy Equation \eqref{eqn:T Theorem1}, where $\sR_x(T)$ and $\sR_p^i(T)$ are as defined in \eqref{eqn:R_x and R_p}. Notably, both $\sR_x(T)$ and $\sR_p^i(T)$ do not have direct dependence on the dimensions $d$ and $n$. Therefore, the total number of iterations $T$ will also not depend directly on the dimensions.
In Section \ref{subsec:efficient feas},
we show that this independence further extends to $\tilde{T}$, the total iterations of Algorithm \ref{alg:SOFO} when using efficient feasibility testing.
Therefore, both $T$ and $\tilde{T}$ do not have direct dependence on the dimensions $d$ and $n$. 
Moreover, as we discussed in \eqref{eqn:Omega bound}, it is adequate to let $\Omega \sim O(\ln(\frac{m}{\nu_0}))$. This, in turn, makes $T$ and $\tilde{T}$ scale logarithmically with the number of constraints $m$.
Lastly, as outlined in Section \ref{subsec:epsilon SMD} and Appendix \ref{appn:efficient update}, the complexity of updating $x$ and $p$ is almost independent of $n$. Overall, these results imply that our algorithm is scalable to $n$ and $d$ under the $\chi^2$-divergence set.
\hfill \Halmos
\end{remark}

%% file: Experiments.tex
In this section, we conduct numerical experiments to compare the performance of our stochastic online first-order meta-algorithm (SOFO-based approach) with the deterministic online first-order framework (OFO-based approach) of \citet{ho2018online}. 
Note that we do not compare our approach with other algorithms, such as the pessimization oracle-based algorithm proposed by \cite{mutapcic2009cutting} and the dual-subgradient meta-algorithm proposed by \cite{ben2015oracle}, because these algorithms tend to struggle with the large-dimensional problem sizes considered in our study. Additionally, \cite{ho2018online} have already demonstrated the superior performance of their OFO-based approach as compared to these algorithms, especially as the dimensions get larger. As we consider problem sizes comparable to and even larger than \cite{ho2018online} (especially larger with regard to the dimenson $n$ of the distribution/uncertain parameters), we opt to only compare with the OFO-based approach of \cite{ho2018online}.

\subsection{Logistic Regression with Fairness Constraints}
We first consider logistic regression with fairness constraints that reduce the correlation between the predicted labels of the classifier and sensitive attributes such as race or gender. In particular, we consider the following formulation adopted from \citet{zafar2017fairness}.
    \begin{align*}
    &\underset{\theta}{\min} - \frac{1}{n}\sum_{i=1}^n y_i \log (\sigma(\theta^\top x_i)) + (1-y_i) \log(1-\sigma(\theta^\top x_i))\\
    &\text{s.t. } -c \leq \frac{1}{n}\sum_{i=1}^n (z_i - \bar{z})\theta^\top x_i \leq c.
    \end{align*}    
In this formulation, $x_i$ denotes the feature vector for $i$-th sample and $y_i \in \{0,1\}$ denotes the corresponding response. The function $\sigma(\cdot)$ is the standard sigmoid function and $z_i$ is a sensitive attribute with sample mean $\bar{z}$. As shown in \citet{zafar2017fairness}, $\frac{1}{n}\sum_{i=1}^n (z_i - \bar{z})\theta^\top x_i$ approximates the covariance between the sensitive attribute and the decision boundary of the classifier. We want $c$ to be small so that the classifier is uncorrelated to the sensitive attribute, and $c$ is often chosen to satisfy certain fairness criteria. 
\par 
We use the adult income dataset \citep{asuncion2007uci}, which contains 45,000 examples with 14 features, including age, gender, and marital status. Our goal is to classify whether the income level is greater than \$5K, and we consider gender as a sensitive attribute.
In addition to the provided features, we add new features generated by considering all polynomial combinations of continuous variables with degrees less than equal to either 3 or 4.
This results in the dimension of the features, and correspondingly $\theta$, being either $d=174$ or $d=300$, respectively.
For simplicity of comparison, instead of solving the optimization problem, we solve the corresponding feasibility problem with the right-hand side of the objective function constraint set to $0.5$; as the optimization problem can be reduced to several feasibility problems, the results will be similar. 
We present the detailed parameter setting in Appendix \ref{appn:FL param}.
We simulate 20 times for each experiment and setting of parametrs that we vary. The experiments are performed on a Linux machine with a 2.3-GHz processor and 64GB memory using Python v3.7. While computing the SP gap, if a closed-form solution is not available, we used CVXPY v1.1.7 \citep{diamond2016cvxpy} and Mosek v10.0.34 \citep{mosek} with standard parameter settings.
\par
\begin{figure}[t]
    \centering
    \includegraphics[width=10cm]{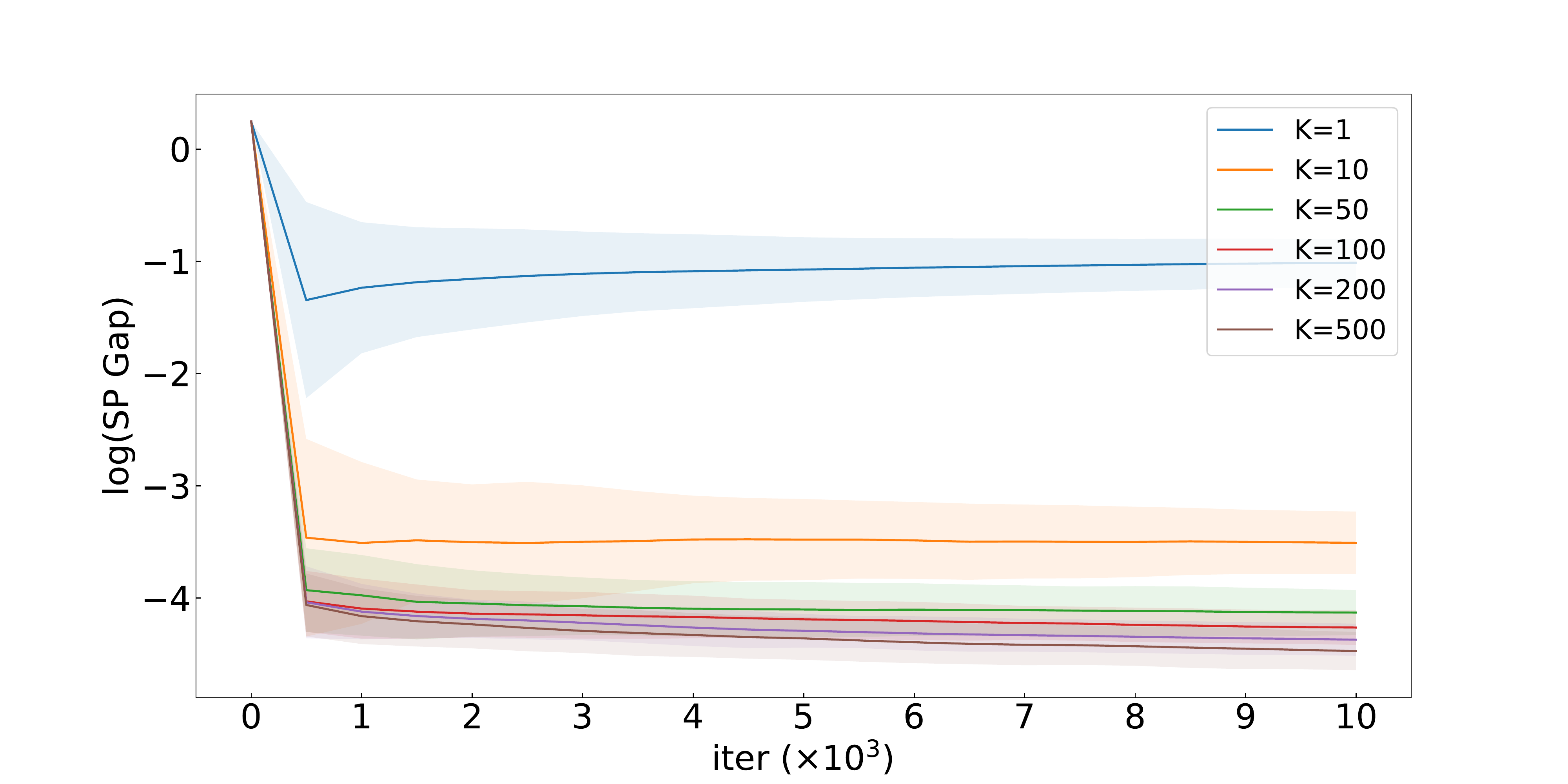}
    \caption{Average log(SP Gap) versus the number of iterations of SOFO, across different values of $K$  $(d=300, n=45000)$.}
    \label{fig:K Iter}
\end{figure}

 
To find an appropriate practical per-iteration samples size $K$, we tested different values and report the results in Table \ref{table:approximation of i}. Specifically, Table \ref{table:approximation of i} shows the average percentage of iterations that satisfy either $\hat{i}_t=i_t^*$ or $|f_t^{\hat{i}_t}(x_t) - f_t^{i_t^*}(x_t)| \leq C_K \epsilon$ out of the first 10,000 iterations (when $n=45000$). We can see that even $K=1$ has around 38\% accuracy, as we only have three constraints. Furthermore, the accuracy increases as $K$ increases. Figure \ref{fig:K Iter} compares the convergence rate of the SOFO-based approach across different values of $K$ when $d=300$. As our theory suggests, $K$ determines the convergence rate and also the precision of the solution that we can obtain. For the following experiments, we choose $K=200$, which was large enough for SOFO to achieve the desired solution for our setting of $\epsilon$ ($\epsilon = 0.02$).
\begin{table}[t]
\begin{center}
\caption{Average percentage of iterations with either $\hat{i}_t=i_t^*$ or $|f_t^{\hat{i}_t}(x_t) - f_t^{i_t^*}(x_t)| \leq C_K \epsilon$ out of first 10000 iterations ($n = 45000$).}

\begin{tabular}{ |c| c | c | c | c | c | c| }
\toprule
  \backslashbox{$d$}{$K$} & 1 & 10 & 50 & 100 & 200 & 500 \\ 
   \hline
   174& 38.4 & 43.4 &49.5&53.2&56.4&60.3\\
   \hline
 300& 38.4 & 43.5 & 49.7  & 53.2 & 56.5 & 60.4 \\
 \bottomrule
\end{tabular}
\label{table:approximation of i}
\end{center}
\end{table}

\par 
\begin{figure}[t]
    \centering
    \includegraphics[width=10cm]{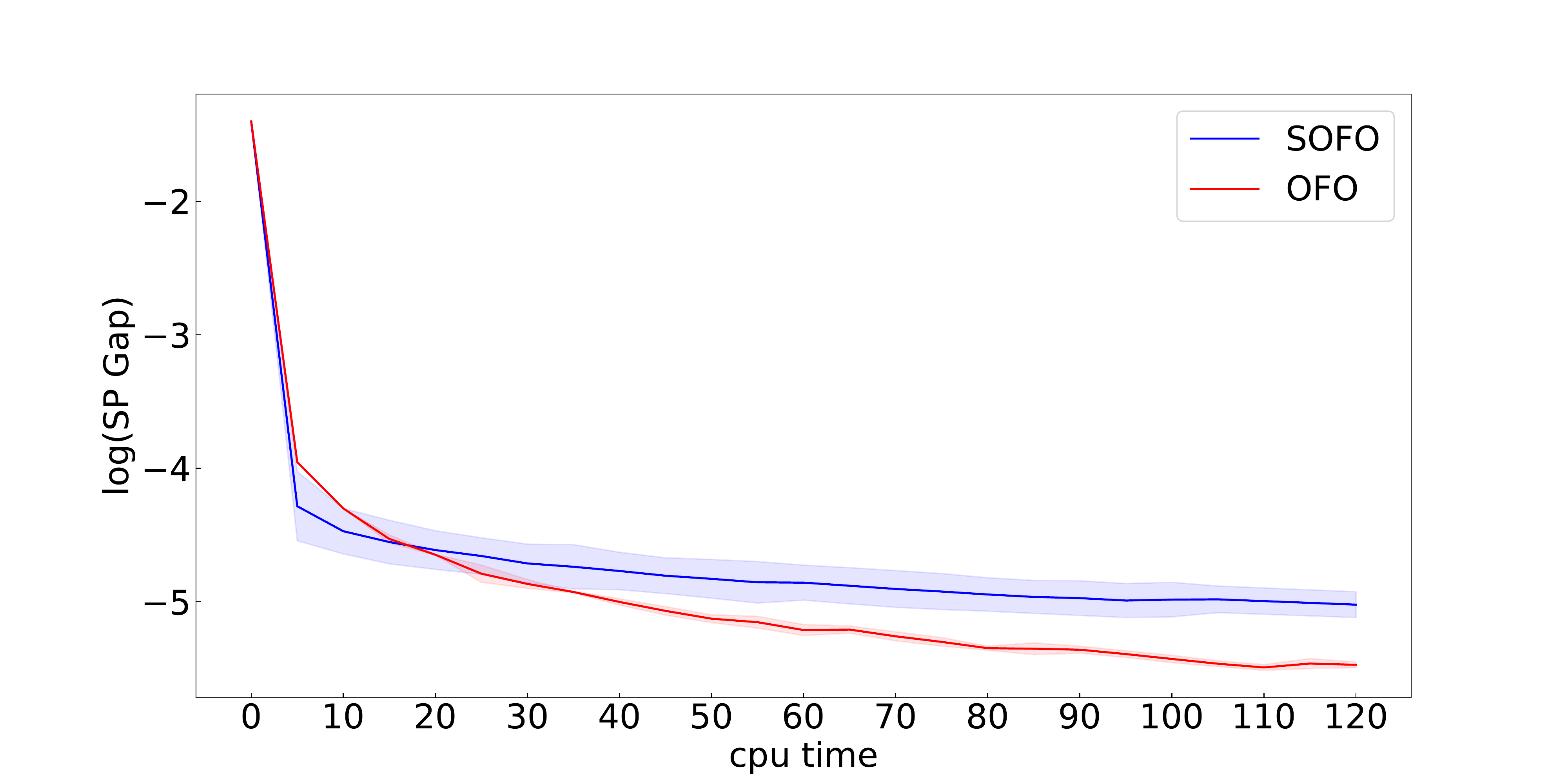}
    \caption{Average log(SP Gap) versus CPU time (in seconds) for SOFO and OFO-based approaches ($d=300, n=45000$). The shaded area represents the 95\% confidence intervals.}%
    \label{fig:K time test}%
    \end{figure}
Figure \ref{fig:K time test} plots the convergence speed between the SOFO-based and OFO-based approaches. The SOFO-based approach obtains the solution with a small SP gap ($\epsilon^\phi(\bar{x},\bar{p}) = 0.007$) faster than the OFO-based approach, while OFO eventually achieves a better solution. This result highlights the limitations of our work in the regime where a very accurate solution, i.e., a small value of $\epsilon$ is desired. This result also follows the commonly recognized trend that stochastic gradient methods perform better early on, but full gradient methods perform better after enough computation time.

\begin{figure}[t]
    \centering
    \subfloat[$d = 174$]{{\includegraphics[width=7.5cm]{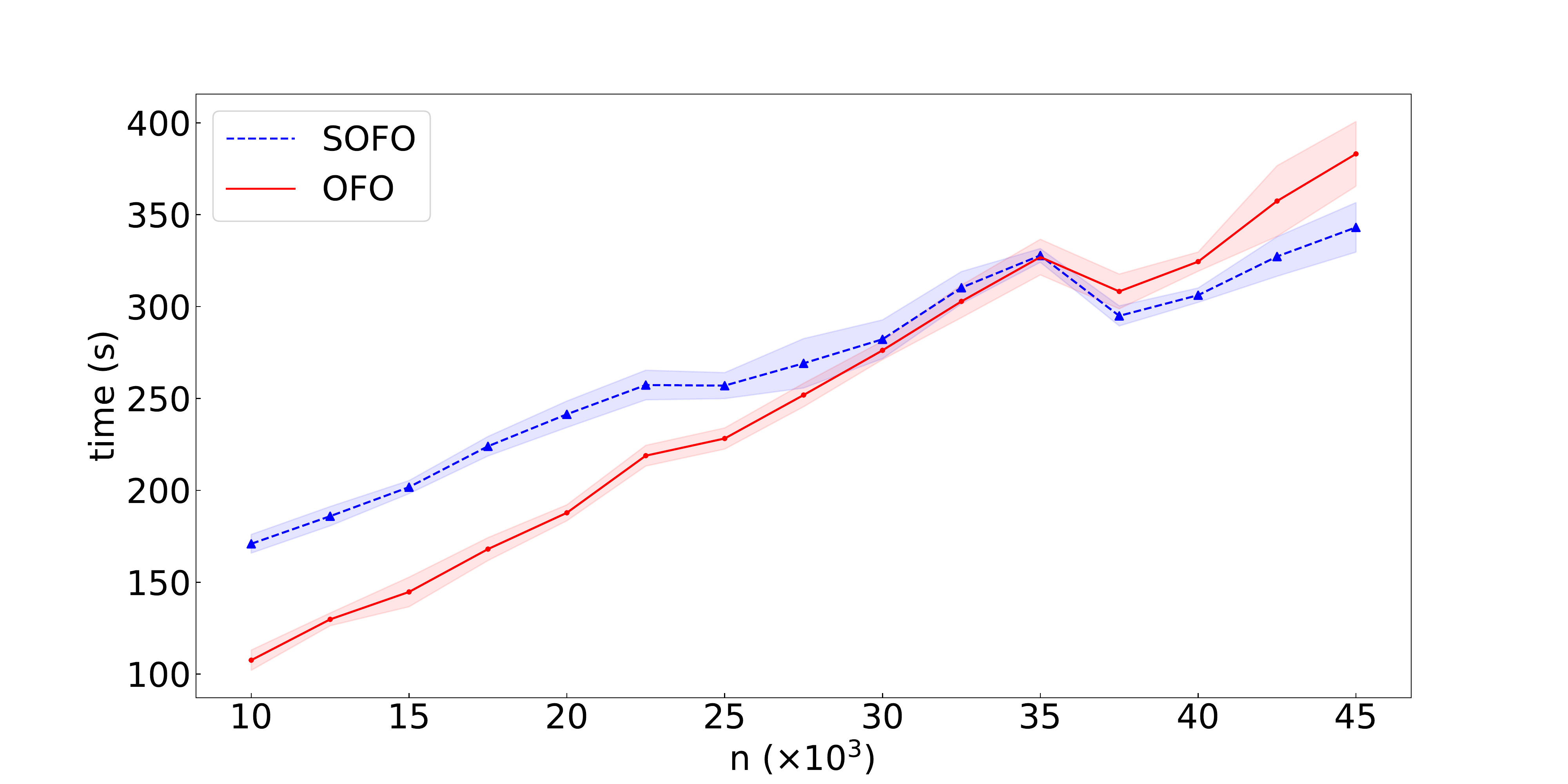} }}%
    \qquad
    \subfloat[$d =300$]{{\includegraphics[width=7.5cm]{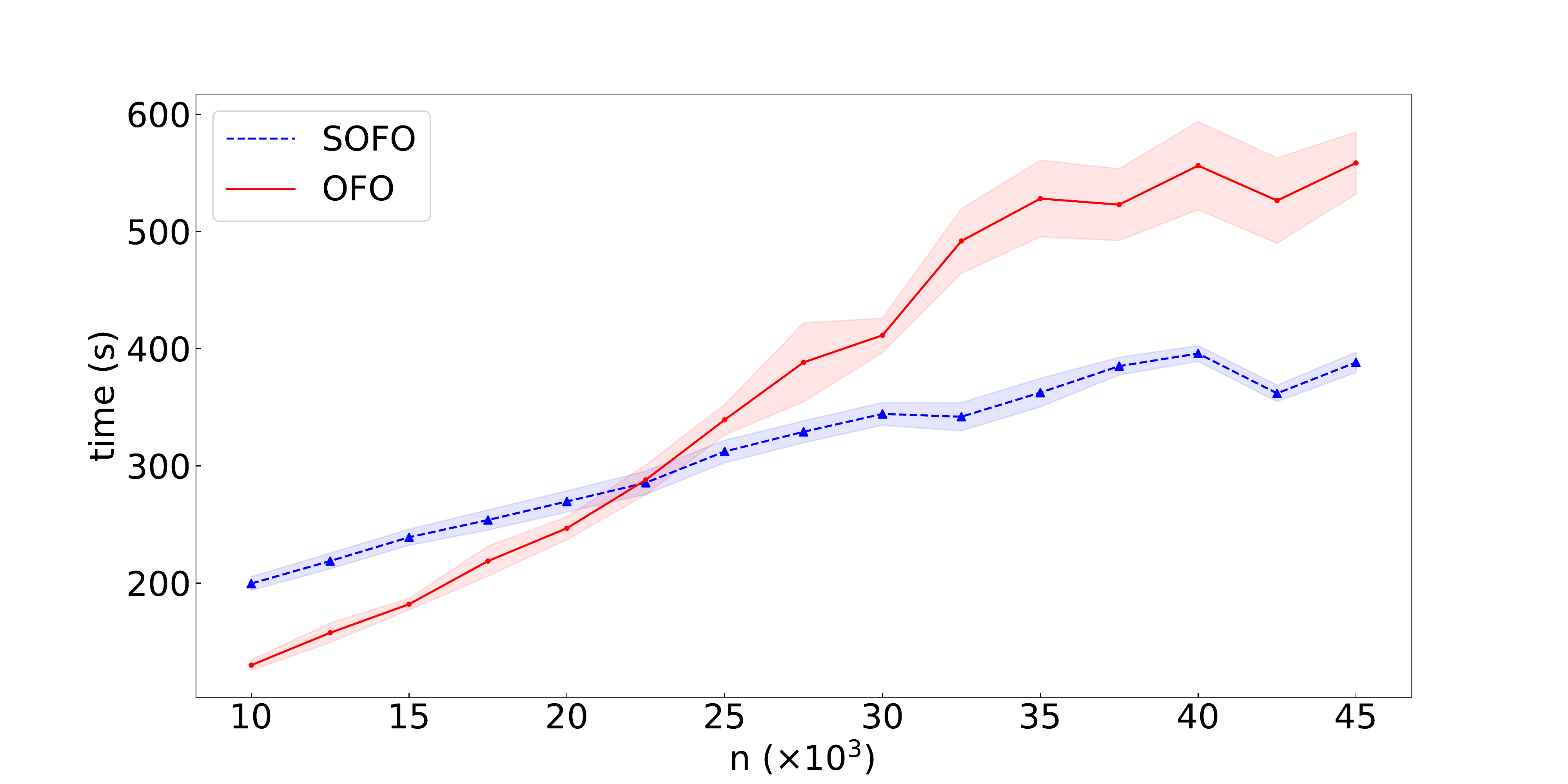}}}%
    \caption{Average solve times (in seconds) to find a solution for SOFO across different values of $n$. The shaded area represents the 95\% confidence intervals.}%
    \label{fig:n num test}%
\end{figure}

Figure \ref{fig:n num test} plots the average solve times in seconds for SOFO-based and OFO-based approaches across different values of $n$. When $n<45000$, we randomly sample $n$ samples from the dataset for every simulation. Since the OFO-based approach does not use the early termination criterion (by computing the SP gap), we also exclude this feature from our SOFO-based approach; although we expect that including the early termination criterion would improve our overall computation times with proper parameter settings. As our theory suggests, SOFO scales better than OFO in $n$ for both $d=174$ and $d=300$. Additionally, we observe that SOFO outperforms OFO as the value of $\costF$ increases when $d$ increases from 174 to 300. However, we see that the SOFO-based approach still has a modest linear grown trend in $n$ in both cases. This comes from the update of the cumulative distribution of $p_t^i$, which we use when random sampling an index from distribution parameter $p_t^i$. Although it is possible to resolve this issue in theory by using a tree structure to update the $p_t^i$ variables, for the problems we consider, using the tree structure in our implementation is inefficient compared to updating the cumulative distribution of the $p_t^i$ variables using well-optimized libraries.
Thus, instead we use the efficient update scheme proposed in Appendix \ref{appn:random sampling}. 
\par 
 Table \ref{table:n test detail1} and Table \ref{table:n test detail2} display the total iterations, average solve time, and final saddle point gap for $d=174$ and $d=300$, respectively. Firstly, observe that the total iteration count is nearly unaffected by the dimensions of $x$ and $p$, which aligns with our theoretical findings. Moreover, by comparing the seconds per-iteration for both SOFO-based and OFO-based approaches, we can appreciate the scalability of the SOFO method with respect to both $n$ and $d$. For example, in both tables, OFO achieves a smaller overall solve time when $n = 20000$ whereas SOFO achieves a smaller overall solve time when $n = 45000$. Notably, the most significant difference in overall solve time occurs in the largest problem size considered, $d = 300$ and $n = 45000$, in which case OFO requires about one and half times more solve time than SOFO.
 While both approaches achieve a final SP gap below the desired level of approximate solution ($\epsilon=0.02$), we can observe that the OFO-based approach achieves a smaller gap. This result can likely be attributed to the reduced variance inherent in deterministic algorithms, such as the OFO-based approach. On the other hand, it is again notable that our SOFO-based approach excels in overall solve time, primarily due to its cheaper iterations and further highlighting the efficiency and scalability of our proposed method.
 
\begin{table*}[t]
    \centering
    \caption{Detailed comparison between SOFO and OFO ($d=174$).}
    \begin{tabular}{|c|c|c|c|c|c|}
    \hline
        &  $n$ & Total Iterations & Solve Time (s) & Seconds Per-Iteration (s) & Final SP Gap \\ \hline
        SOFO & 20000 & 121853 & 241.4068  & 0.00198  & 0.00647  \\ \hline
        OFO & 20000 & 22664 & 187.8148 & 0.00829  & 0.00324  \\ \hline
        SOFO & 45000 & 120489 & 343.0915  & 0.00284  & 0.00547  \\ \hline
        OFO & 45000 & 22611 & 383.0715  & 0.01694  & 0.00312 \\ \hline
    \end{tabular}
    \label{table:n test detail1}
\end{table*}

\begin{table}[t]
    \centering
    \caption{Detailed comparison between SOFO and OFO ($d=300$).}
    \begin{tabular}{|c|c|c|c|c|c|}
    \hline
        &  $n$ & Total Iterations & Solve Time (s) & Seconds Per-Iteration (s) & Final SP Gap \\ \hline
        SOFO & 20000 & 129005 & 269.68359  & 0.00209  & 0.00671  \\ \hline
        OFO & 20000 & 25278 & 246.83508  & 0.00976  & 0.00330  \\ \hline
        SOFO & 45000 & 127561 & 388.28495  & 0.00304  & 0.00536  \\ \hline
        OFO & 45000 & 25224 & 558.37625  & 0.02214  & 0.00282 \\ \hline
    \end{tabular}
    \label{table:n test detail2}
\end{table}

\subsection{Personalized Parameter Selection in Large-Scale Social Network}\label{appn:social network}
In this section, we present another experimental setting for personalized parameter selection in a large-scale social network.
This experiment is adopted from \cite{basu2019optimal}. Social network companies may set certain ``parameters,'' such as the gap between different ad displays on the news feed, for different users in order to personalize the experience for each user and subsequently achieve goals with respect to performance metrics. Clearly, choosing a global parameter for all the customers might not be the optimal solution. Instead, we consider dividing our customer base into $J$ different cohorts and estimate the optimal parameter for each cohort. Let us assume that the parameter can take $L$ possible values, and there are $m$ metrics that we consider. Among these $m$ metrics, revenue is the primary metric that we would like to maximize, and the other $m-1$ metrics, such as click-through rates (CTR), will appear as constraints in our problem. Let our decision variable $x_{jl}$ denote the probability of assigning the $l$-th treatment to the $j-$th cohort. Also, we define $d$ as $d = J \times L$, which is a dimension of our decision variable $x$, and define a random vector $\bu^i \in \bbR^d$ as $\bu^i = (u_{11}^i, \dots, u_{JL}^i)$, where $u_{jl}^i$ represents the casual effect when treatment $l \in L$ is applied to a cohort $j \in [J]$ for a metric $i \in [m]$. We assume that $\bu^i \sim p^i$ with some distribution $p^i$. Suppose we have $n$ samples of $\bu^i$, where $r-$th sample is denoted as $\bu_r^i$, for all $i \in [m]$. Then $F_r^i(x)$ can be expressed as $F_r^i(x) = (\bu_r^i)^\top x$. Our goal is to find an optimal allocation $x^*$ that maximizes revenue under distributionally robust constraints, hence the problem formulation is:
\begin{equation}
\begin{aligned}
    &\underset{x}{\text{max}} &&\underset{p^1 \in \calP^1}{\inf} \  p^{1^\top}\bF^1(x)\\
    &\text{s.t.} \ &&\underset{p^i \in \calP^i}{\sup} (p^i)^\top\bF^i(x) \leq c^i, \quad i = 2,\dots,m\\
    & && 0 \leq x_{jl} \leq 1, \forall j \in [J], l \in [L]\\
    & &&\sum\limits_l x_{jl} = 1 \quad \forall j \in [J].
\end{aligned} \label{eqn:social network}
\end{equation}
In our experiment, we solve the feasibility version of \eqref{eqn:social network}, where the objective function is changed to an inequality $\underset{p^1 \in \calP^1}{\inf} \ p^{1^\top}\bF^1(x) \geq c^1$, rather than solving the decision problem itself.
We use synthetically generated data for our experiment, where the generation process is described in Appendix \ref{appn:ps datagen}. Unless stated otherwise, we set the tolerance level to $\epsilon = 0.05$. Details of other parameter settings can be found in Appendix \ref{appn:social param}
\par 
Figure \ref{fig:K time test-appn} plots the average log duality gap against the CPU-time. Parameter settings are the same as in the previous experiment, except that we fixed the per-iteration sample size $C_K$ to 100, and we set the tolerance level to $\epsilon = 0.001$ in this experiment. Similarly to Figure \ref{fig:K time test}, we see that SOFO arrives at the small SP gap solution faster than OFO but, eventually, OFO catches up when $n=1000$. However, when $n=5000$, there is no intersection between the SOFO and OFO curves even after more than 900 seconds. 
\par 
\begin{figure}[ht]
    \centering
    \subfloat[$J = 10, L = 15$]{{\includegraphics[width=7.5cm]{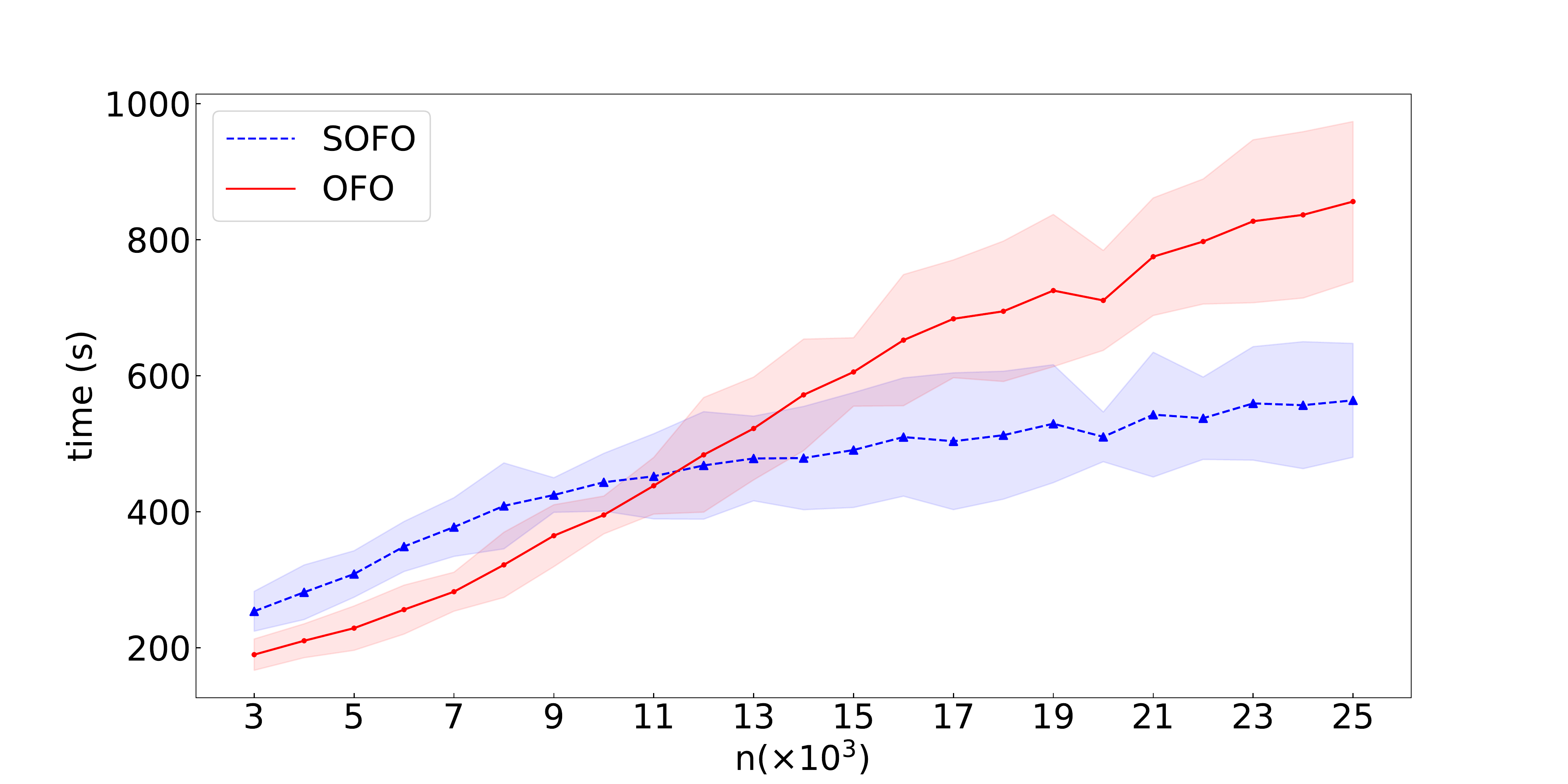} }}%
    \qquad
    \subfloat[$J = 10, L = 25$]{{\includegraphics[width=7.5cm]{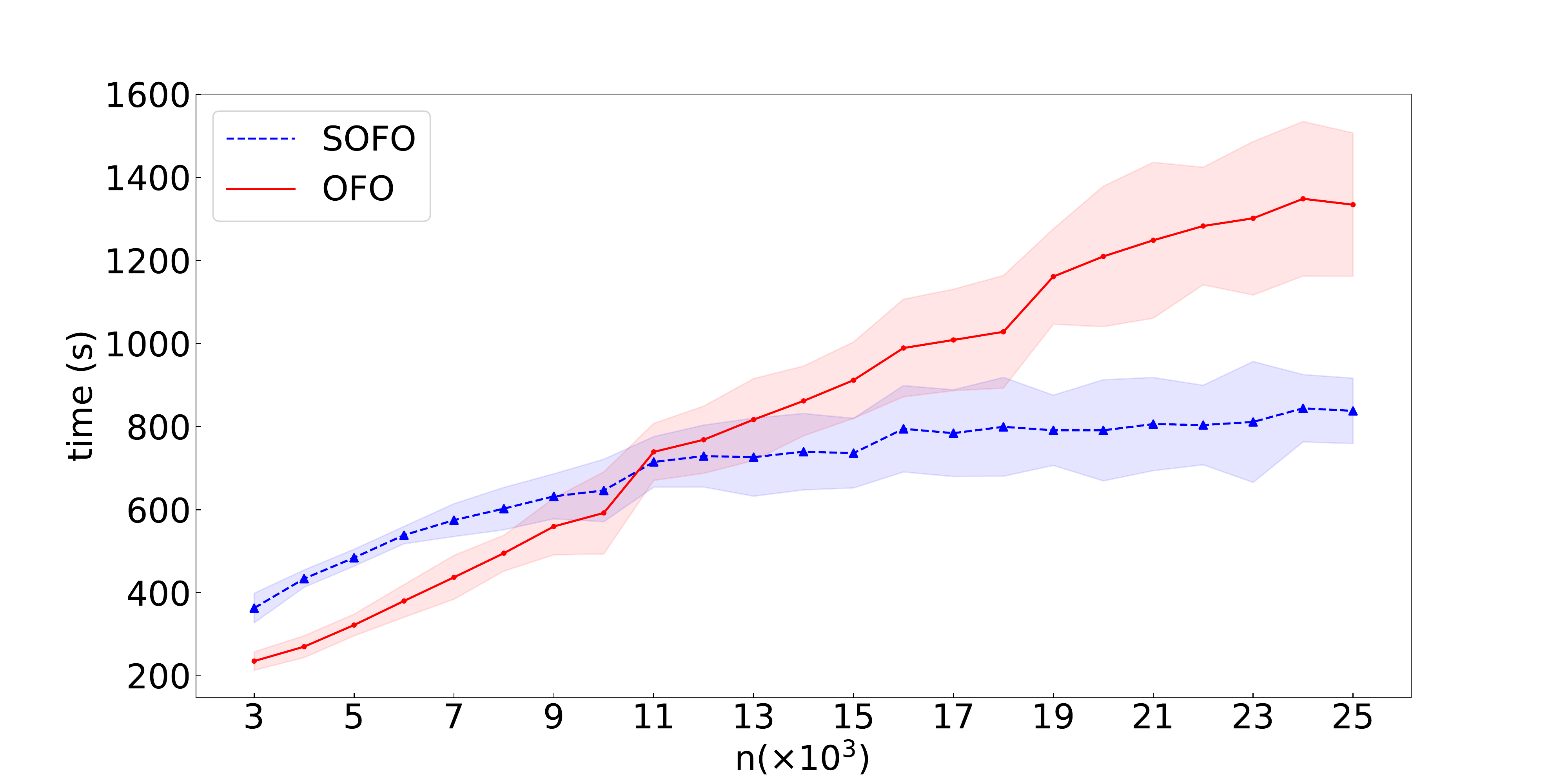} }}%
    \caption{Average solve times (in seconds) for different $n$ ($m = 5, K= 100$). The shaded area represents the 95\% confidence intervals.}%
    \label{fig:n num test-appn}%
\end{figure}

Figure \ref{fig:n num test-appn} plots the average solve times in seconds for the SOFO-based approach and the OFO-based approach for different $n$. We set $c^i = 1.1 \times p_0^{i^\top} \bF^i(x_0)$ for $i \neq 1$ to have a larger feasibility region and $c^1 = 1$. The early termination condition is used to maximize the efficiency of Algorithm \ref{alg:SOFO}. As our theory suggests, in the large $n$ regime ($n \geq 12000$), the SOFO curve becomes flat, while the solve time of the OFO-based approach increases linearly in $n$. Lastly, Table \ref{table:average total iter and per iter} shows the average total iterations and seconds per iteration for the OFO-based approach and the SOFO-based approaches. We can see that the average iterations of the SOFO-based approach is four times higher than that of the OFO-based approach. Also, we see that the seconds per-iteration scales better in the SOFO-based approach than in the OFO-based approach with respect to $n$.
\begin{table}[ht]
\begin{center}
\caption{Average total iterations and seconds per-iteration for OFO and SOFO ($J = 10, L =25, m = 5,\sigma^2 = 0.1, K = 100$).}
    \begin{tabular}{|c|c|c|c|}
    \toprule
         & $n$ & Total Iterations & Seconds Per-Iteration \\
         \hline
    SOFO & 5000 & 208445
 & 0.0023
 \\ 
        \hline
    OFO & 5000 & 55996
 &  0.0058
\\
        \hline
    SOFO & 25000 & 180120
 &0.0046
 \\
        \hline
    OFO & 25000 & 57219
 & 0.023
\\
    \bottomrule
    \end{tabular}
    
    \label{table:average total iter and per iter}
\end{center}
\end{table}

\subsection{Multi-item Newsvendor with Conditional Value at Risk Constraint} \label{subsec:MNVI}
Finally, we consider a multi-item newsvendor problem with conditional value at risk (CVaR) constraint. In the multi-item newsvendor problem, the demand for multiple perishable products is uncertain. Any sold item produces revenue, and any inventory replenishment incurs an ordering/production cost. Also, any unsatisfied demand incurs back-order costs, and we can sell any left-over items at their salvage prices. Our goal is to decide the optimal order quantity of items to minimize the loss while satisfying CVaR and budget constraints. Suppose there are $d$ items and vectors $c,r,s,b \in \bbR^d$ respectively represent production cost, retail price, salvage price and back-order cost of the items and $\xi \in \bbR^d$ denote a random demand vector. It is natural to assume that $c<r$ and $s<r$, and this implies that the optimal order quantities of the items are equal to $\min \{x,\xi \}$, where the minimization function is applied element-wise. Under this notation, we can express our loss function as 
\begin{align*}
    L(x,\xi) &= c^\top x - r^\top \min \{x,\xi \} - s^\top (x - \min \{x,\xi \})+ b^\top (\xi- \min \{x,\xi \}) \\
    & = (c-s)^\top x  - (b + r - s)^\top \min \{x,\xi \} + b^\top \xi
\end{align*}
As we assumed $r>s$, $L(x,\xi)$ is a convex, and in fact piece-wise linear, function of $x$.
Conditional Value at Risk (CVaR) at level $\beta$ is a widely used risk measure interpreted as the expectation of the 100 $\times \beta \%$ worst outcomes of the loss distribution. We assume that the seller is risk-averse so that the seller wants CVaR at level $\beta$ lower than some specified value $\alpha$. We define CVaR at level $\beta$ with respect to distribution $p$ as
\begin{align*}
    \text{CVaR}_{p,\beta}(L(x,\xi)) = \underset{\tau}{\min} \big\{ \tau + \frac{1}{\beta}\mathbb{E}_{\xi \sim p}\big[ (L(x,\xi) - \tau)^+ \big] \big\}.
\end{align*}
Incorporating this into the DRO formulation, we require the order portfolio $x$ to satisfy $$\underset{p \in \calP}{\sup} \ \mathbb{E}_{\xi \sim p}\big[ \text{CVaR}_{p,\beta}(L(x,\xi)) \big]\leq \alpha.$$
Then, the DRO formulation of this problem is:
\begin{equation}
    \begin{aligned}
    &\underset{\tau, x}{\min} \ &&\underset{p^0 \sim \calP}{\sup} \mathbb{E}_{\xi^0 \sim p^0}\big[ L(x,\xi^0) \big] \\
    &\text{s.t.} &&\underset{p^1 \sim \calP}{\sup} \tau + \frac{1}{\beta}\mathbb{E}_{\xi^1 \sim p^1}\big[ (L(x,\xi) - \tau)^+ \big] \leq \alpha\\
    & && x \in S:= \{x| \mathbbm{1}^\top x \leq C, x \geq 0\}.
    \end{aligned}
    \label{form:MINV}
\end{equation}
The set $S$ represents the budget limit in the order portfolio.
A similar model has been studied by \citet{hanasusanto2015distributionally}, where they instead considered a Lagrangian formulation that moves the CVaR constraint into the objective function. We instead consider directly solving the constrained form as given in \eqref{form:MINV}.

In this experiment, we use a synthetically generated data set. We provide the data generation procedure in Appendix \ref{appn:data gen}. In contrast to the previous experiments, we solve the optimization problem \eqref{form:MINV} itself by solving $O(\frac{1}{\ln (\epsilon)})$ feasibility problems and also utilize the early stopping termination by calculating the SP gap. To speed up both SOFO and OFO-based approaches, we use a warm startup technique, which uses the output $\bar{x}$ and $\bar{p}$ of the previous feasibility iteration to initialize $x_0$ and $p_0$ for the current feasibility problem. Later in this section, we present results that show the effectiveness of the warm startup. \par 
We set $T_s$ to be $\frac{T}{5}$ for SOFO and $T_s$ to be $\frac{T}{10}$ for OFO in the first feasibility iteration. To maximize the efficiency of the warm startup strategy, we set $T_s$ to $\frac{T}{10}$ for SOFO and $\frac{T}{20}$ for OFO after the first feasibility iteration. Lastly, we set our tolerance level $\epsilon$ to 0.03 and run 20 instances for each setting. Details of other parameters can be found in Appendix \ref{appn:MINV param}.
After running a similar test as in the previous experiments, we choose $K$ to be 50. 
Figure \ref{fig:n num test, MINV} plots the average solve times in seconds against different $n$. Similar to the previous experiments, we see that OFO is more sensitive than SOFO to the increasing scale of $n$.
\begin{figure}[H]
    \centering 
    {\includegraphics[width=10cm]{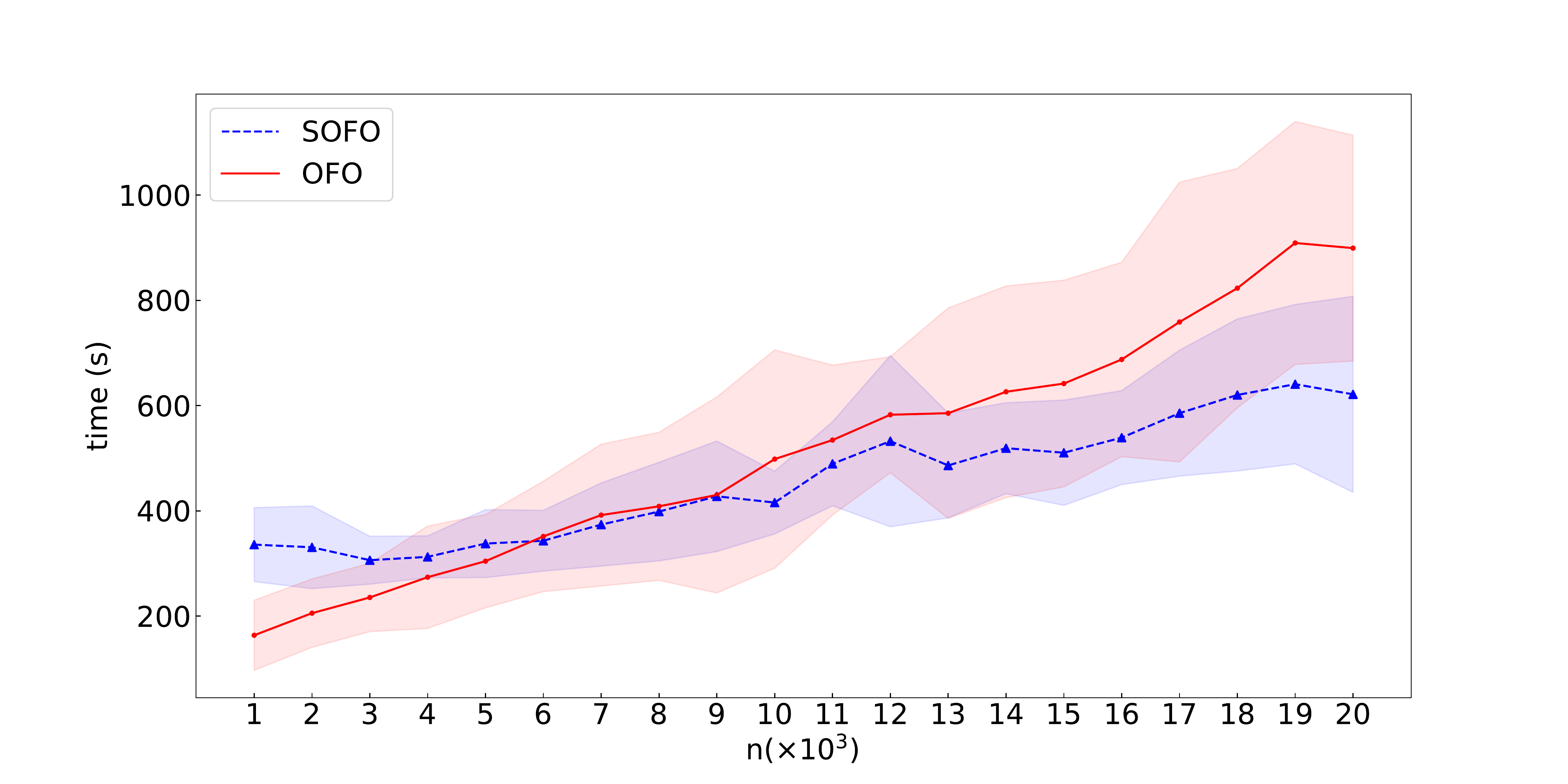} }%
    \caption{Average solve time (in seconds) for different $n$.}%
    \label{fig:n num test, MINV}%
\end{figure}
Lastly, we conclude this experiment by showcasing the effectiveness of the warm start strategy. Table \ref{table:warm startup, MINV} displays the solve times for SOFO and OFO with and without warm startup, where WS stands for warm starting.
\begin{table}[h]
\begin{center}
\begin{tabular}{| c | c | c| c | c|}
\toprule
  $n$ & SOFO w/ WS & SOFO w/o WS & OFO w/ WS & OFO w/o WS\\ 
   \hline
 5000& 337 & 565 & 304 & 422\\
 \hline 
 10000 & 415 & 769 & 498 & 654 \\
 \bottomrule
\end{tabular}
\caption{Average solve times (in seconds) with and without warmstarts (WS).}
\label{table:warm startup, MINV}
\end{center}
\end{table}
Observe that the solve times for both SOFO and OFO are reduced when using warm starting. In particular, it is apparent that SOFO has a slight advantage in relative speedup due to warm starting as compared to OFO.

%% file: Conclusion.tex
We provided a new stochastic first-order meta-algorithm for constrained distributionally robust feasibility problems. To address the issue of calculating stochsatic gradients in the presence of multiple constraints, we propose a new stochastic first-order method called $\epsilon$-SMD. We demonstrate the scalability of $\epsilon$-SMD, and we also develop specialized results for $\chi^2$-divergence set to improve the scalability even further. We numerically validate the improved performance of our meta-algorithm by comparing it with its deterministic counterpart \citep{ho2018online}, and we observe that, especially for large $n$, our stochastic meta-algorithm obtains a solution with a small optimality gap in a reasonable time. 
There are several directions for future research, including developing specialized methods and results for ambiguity sets beyond the $\chi^2$-divergence set and considering problems with more intricate stochastic constraints (e.g., chance constraints).

%% file: Appendix.tex
\section{Additional Results and Discussion}
In this appendix, we provide some additional results and provide complete discussion of several points that were mentioned in the main text.

\subsection{Generalization to Dependent Case}
\label{appn:dependent RV}
In this section, we show that we can reduce an $\epsilon$-feasibility problem with dependent random variables $z^1, \cdots, z^m$, to a $\epsilon$-feasibility problem with independent random variables. Let $q \in \mathbb{R}^l$ be a joint probability distribution of random variables $z^1, \dots, z^m$  and define an uncertainty set of the joint distribution $q$ as $\calQ$. Then, we can interpret $p^i$ as a marginal distribution of $z^i$. Let $A^i$ denote a projection matrix that projects the joint distribution $q$ to a marginal distribution $p^i$ for each $i \in [m]$:    \begin{align*}
        A^i q = p^i, \quad \forall i \in [m]. \label{eqn:projection of q}
    \end{align*}
This setting allows us to capture the dependency between random variables $z^1, \dots, z^m$. If $z^i$ for all $i \in [m]$ are mutually independent, then we set $q$ as $q \in \mathbb{R}^{l}$ with $l = m\times n$ and $A^i \in \mathbb{R}^{n \times l}$, where the first $n$ elements represent $p^1$ and the next $n$ elements represent $p^2$ and so on. Let $\mathbf{0}_{n}$ denote $n \times n$ matrix that all elements are 0 and $I_n$ denote $n-$dimensional identity matrix. Then, in the independent case, we set $A^i$ to be
\begin{align*}
 A^i = \begin{bmatrix} \mathbf{0}_{n} & \cdots & \mathbf{0}_{n} & I_n & \mathbf{0}_{n} & \cdots & \mathbf{0}_{n}  \end{bmatrix}
\end{align*}
where only $i-$th $n \times n$ matrix block is equal to $I_n$, otherwise $\mathbf{0}_{n}$. Let us return to the dependent setting and define a vector-valued function $\mathbf{G}^i: X \rightarrow \mathbb{R}^l $ as $\mathbf{G}^i(x) = A^{i^\top} \bF^i(x)$ for all $i \in [m]$. Under our newly introduced notations, we have
\begin{align*}
\mathbb{E}_{z^i \sim p^i}[F^i (x,z^i)] = (p^i)^\top\bF^i(x) = q^\top \bG^i(x).
\end{align*}
Then we can express $\epsilon$-approximate version of distributionally robust feasibility problem as:
\begin{equation}
    \begin{cases}
    &\epsilon -\text{feasible:} \quad \underset{x \in X}{\text{inf}} \  \underset{q \in \calQ}{\text{sup}} \  \underset{i \in [m]}{\max} \ q^\top \bG^i(x)\leq \epsilon \\
    &\text{infeasible:} \quad \underset{x \in X}{\text{inf}} \    \underset{q \in \calQ}{\text{sup}} \  \underset{i \in [m]}{\max} \ q^\top \bG^i(x) >0 \label{form:DROfeas-G}
\end{cases}    
\end{equation} 
The main structural difference between $\epsilon$-feasibility problem and \eqref{form:DROfeas-G} is that there is a common variable $q$ inside the max function over $i \in [m]$ in \eqref{form:DROfeas-G}, while there are separate variables $\{p^i \}_{i=1}^m$ inside the max function over $i \in [m]$ in $\epsilon$-feasibility problem. However, we can reduce the dependent case to the independent case by generating independent copies of $q$. Suppose $p^i$ be an independent copy of $q$, and $\calP^i$ be a copy of uncertainty set $\calQ$, for all $i \in [m]$. Then, we have 
\begin{align}
    \underset{q \in \calQ}{\sup} \ \underset{i \in [m]}{\max} \ q^\top \bG^i(x) =  \underset{i \in [m]}{\max} \   \underset{q \in \calQ}{\sup}  \ q^\top \bG^i(x) = \underset{i \in [m]}{\max} \   \underset{p^i \in \calP^i}{\sup}  \ (p^i)^\top \bG^i(x) =  \underset{p \in \calP}{\sup} \ \underset{i \in [m]}{\max} \ (p^i)^\top \bG^i(x) \label{eqn:dep to indep}
\end{align}
Using \eqref{eqn:dep to indep}, we can reduce the dependent case to the independent case. 
\begin{remark}
We can easily convert the case where $z^i$ for all $i \in [m]$ are not independent to the independent case by defining a new random variable $\tilde{z}^i$, which is an independent copy of a random variable of a joint distribution of $z^1, \cdots, z^m$.
\end{remark}
\par 
It is natural to ask why we need to work on the independent structure instead of the dependent structure, as converting the dependent case to the independent case increases the distribution parameter's dimension at most $m$ times.
Let us assume dependence between the random variables $z^1, \cdots, z^m$ and consider solving \eqref{form:DROfeas-G} to answer this question.
Notice that we can still apply Theorem \ref{thm:Nam thm 3.1} for the dependent case. However, it is not well studied how to generate a sequence $\{x_t,p_t,\theta_t\}_{t=1}^T$ that satisfies termination condition 
\begin{align*}
    \espbullet + \espcirc \leq \frac{\epsilon}{2}.
\end{align*}
with the convex-non-concave function $\underset{i \in [m]}{\max} \ q^\top \bG^i(x)$. Also, even if we know the way to generate such a sequence by using non-convex optimization methods, usually those methods requires strong assumption on the function $F_r^i(x)$ such as strong convexity to guarantee convergence. However, by assuming the independence between the random variables, we can use online convex optimization (OCO) tools to solve $\epsilon$-feasibility problem with a mild assumption on the function $F_r^i(x)$.

\subsection{Saddle Point Gap Computation}\label{appn:SP gap}
In this section, we explain methods to compute the SP gap as defined in \eqref{eqn:SP gap}. In order to compute the SP gap, we need to solve $\underset{p \in \calP}{\sup} \ \phi(\bar{x},p)$ and $\underset{x \in X}{\inf} \ \phi(x, \bar{p})$. In general, we have 
\begin{equation}
    \begin{aligned}
        \underset{p \in \calP}{\sup} \ \phi(\bar{x},p) = \underset{i \in [m]}{\max} \ \underset{p^i \in \calP^i}{\sup} (p^i)^\top \bF^i(\bar{x}),
    \end{aligned}
    \label{eqn:sup pi}
\end{equation}
and we have to solve $m$ independent convex optimization with linear objective to solve \eqref{eqn:sup pi}. The complexity of solving each convex optimization largely depends on the structure of the convex set $\calP$. However, in some cases such as $\calP^i$ being $f$-divergence uncertainty set, \citep{namkoong2016stochastic} showed that we are able to obtain the explicit solution of \eqref{eqn:sup pi} in $O(n)$.

\par
To calculate $\underset{p \in \calP}{\sup} \ \phi(\bar{x},p)$, we use KKT conditions to derive closed-form solution for $\underset{p \in \calP}{\sup} \ \phi(\bar{x},p)$. As $\calP^i$ shares the same stucture for all $i \in [m]$, we only need to consider solving $\underset{p^i \in \calP^i}{\sup} (p^i)^\top \bF^i(\bar{x})$. For simplicity, we omit the constraint index $i$ here. We change $\underset{p \in \calP}{\sup} \  p^\top \bF(\bar{x})$ to $-\underset{p \in \calP}{\inf} \  p^\top \bF(\bar{x})$ and derive its closed-form solution. Notice that $F_r(\bar{x})$ is a $r-$th element of a vector $\bF(\bar{x})$. We define Lagrangian function $\mathcal{L}(p; \lambda, \theta)$ as
\begin{equation}
    \begin{aligned}
        \mathcal{L}(p; \lambda, \theta) = -p^\top \bF(\bar{x}) - \frac{\lambda}{n^2} \Big( \rho  - \sum_{r=1}^n \frac{(np^{(r)}-1)^2}{2} \Big) - \theta^{\top}(p - \frac{\delta}{n}\mathbbm{1}_n) \label{eqn:Lagrangian2}
    \end{aligned}
\end{equation}
where $\theta \in \mathbb{R}_+^{n}$ and $\lambda \geq 0 $. Using KKT conditions and strict complementary slackness, we get 
\begin{equation}
    \begin{aligned}
        p(\lambda) = \Big( \frac{\bF(\bar{x})}{\lambda} + \frac{1}{n}\mathbbm{1}_n - \frac{\delta}{n}\mathbbm{1}_n \Big)_+ + \frac{\delta}{n}\mathbbm{1}_n
    \end{aligned}
    \label{eqn:p update2}
\end{equation}
where $p(\lambda) = \text{argmin}_{p \in \pset} \text{sup}_{\theta \in \mathbb{R}_+^n} \mathcal{L}(p; \lambda, \theta)$. Substituting $p(\lambda)$ into equation \eqref{eqn:Lagrangian2}, we have 
    \begin{align*}
        g(\lambda) = -p^\top \bF(\bar{x}) - \frac{\lambda}{n^2}\Big(\rho -  \sum_{r=1}^n \frac{(np(\lambda)^{(r)}-1)^2}{2} \Big).
    \end{align*}
As $g(\lambda)$ is a concave function of $\lambda$, $g'(\lambda)$ is a non-increasing function of $\lambda$. We define a set $I(\lambda) = \{ r \in [n]| \frac{\bF(\bar{x})}{\lambda} + \frac{1}{n}\mathbbm{1}_n \geq \frac{\delta}{n}\} $. Then, $g'(\lambda)$ is like: 
    \begin{align*}
         \frac{\partial}{\partial \lambda} g(\lambda) = \frac{1}{2} \sum_{r \in I(\lambda)} \frac{F_r(\bar{x})^2}{\lambda^2} - \frac{\rho}{n^2} + \frac{(1-\delta)^2}{2n^2}(n- |I(\lambda)|)
    \end{align*}
We use bisection search to find $\lambda^*$ such that $\lambda^* \in \underset{\lambda \geq 0}{\text{argmax}} \ g(\lambda)$. We input $\lambda^*$ to \eqref{eqn:p update2} and get $\underset{p \in \calP}{\sup} \  p^\top \bF(\bar{x})$.

\par 
To solve $\underset{x \in X}{\inf} \ \phi(x, \bar{p})$, we need to solve the following convex optimization:
\begin{equation}
    \begin{aligned}
      &\underset{x \in X, t\in \bbR}{\inf} &&t \\
      &\text{s.t.} && t \geq \bar{p}^{i^\top}\bF^i(x), \ \forall i \in [m].
    \end{aligned}
    \label{eqn:inf pi}
\end{equation}
The complexity of solving \eqref{eqn:inf pi} largely depends on the structure of function $F_r^i(x)$. However, in most of the applications, $F_r^i(x)$ would be in a tractable form such as a linear or quadratic function of $x$, and we can solve \eqref{eqn:inf pi} by using state-of-the-art convex optimization solvers such as Gurobi or Mosek.

\subsection{High-Probability Convergence of SMD}
In this section, we introduce the high-probability convergence of online SMD. This is a minor extension of \citet{nemirovski2009robust}. Let us consider the following online setting. We have a finite time horizon $T$, compact and convex domain $X$ and in each time period $t \in [T]$ we have a convex loss function $f_t: X \rightarrow \mathbb{R}$ and weight $\theta_t$ such that $\sum_{t=1}^T \theta_t = 1$ . At every iteration $t \in [T]$, the player must make a decision $x_t \in X$ and suffer a loss $f_t(x_t)$, where $f_t$ is unknown to the player beforehand. We want to minimize the weighted regret
    \begin{align*}
        \sum_{t=1}^T \theta_t f_t(x_t) - \underset{x \in X}{\inf} \ \sum_{t=1}^T \theta_t f_t(x).
    \end{align*}
by using SMD for an update of $x_t$. We assume that the set $X$ follows the mirror descent setup (Assumption \ref{assm:mds}) and our problem satisfies the following light tail assumption. 
\begin{assumption}[light tail]\label{assm:light tail}
Let $x \in X$ be given and $g_t$ be a given stochastic subgradient of $f_t(x)$ that we can obtain from stochastic subgradient oracle. There exists $G_l>0$ such that:
\begin{align}
    \mathbb{E}\Big[ \exp\Big(\frac{\|g_t\|_{x,*}^2}{G_l^2}   \Big) | x \Big] \leq \exp(1), \quad \forall t \in [T]. \label{eqn:light tail}
\end{align}
\end{assumption}
Then, the online SMD has a high-probability convergence guarantee. 
\begin{proposition}[Extension of proposition 2.2 of \citet{nemirovski2009robust}]\label{prop:Nemi high prob conv}
Let $\{x_t\}_{t=1}^T$ be a sequence generated by online stochastic mirror descent (online SMD) with $g_t$ as a stochastic subgradient of $f_t(x)$ at $x_t$ and with diminishing step size
\begin{equation}
    \begin{aligned}
        \alpha_t = \frac{c_x}{\sqrt{t}}, \quad c_x = \frac{1}{G_l}\sqrt{\frac{D_x}{\Omega}}. \label{eqn:high prob ss}
    \end{aligned}
\end{equation}
 Suppose, the light tail assumption \eqref{eqn:light tail} holds with constant $G_l$ and $\theta \in \Delta_T$ be the normalization of $\{\alpha_t\}_{t=1}^T$. For any $\Omega \geq 1$, we have
\begin{align*}
    &\mathbb{P} \Big\{ \sum_{t=1}^T \theta_t f_t(x_t) -\underset{x \in X}{\inf} \sum_{t=1}^T \theta_t f_t(x)> \frac{w(\Omega)\cdot G_l \cdot \ln (T) \cdot \sqrt{2D_x}}{\sqrt{T}}\Big\} \leq 2\exp (-\Omega),
\end{align*}
where $w(\Omega) = 3\sqrt{\Omega}$.
\end{proposition}
\proof{Proof.}
Let us define $h_t$ as $h_t := \mathbb{E}\big[g_t |x_t \big]$ and $e_t$ as $e_t = g_t - h_t$. As $h_t \in \partial f_t(x_t)$, for any $x \in X$,  we have 
\begin{align}
    \alpha_t(f_t(x_t)  - f_t(x)) &\leq \alpha_t h_t^\top (x_t - x) \nonumber\\
    &= \alpha_t(h_t - g_t)^\top (x_t - x) + \alpha_t g_t^\top (x_t - x)\nonumber\\
    &\leq \alpha_te_t^\top (x- x_t) + B_{\psi_x}(x,x_t) - B_{\psi_x}(x,x_{t+1}) +\frac{\alpha_t^2\|g_t\|_*^2}{2}.
    \label{eqn:prop nemi eqn1}
\end{align}
Summing from $t=1$ to $T$ and dropping a negative term $B_{\psi_x}(x,x_{T+1})$, we have
\begin{equation}
    \begin{aligned}
        \sum_{t=1}^T \alpha_t (f_t(x_t) - f_t(x)) \leq B_{\psi_x}(x, x_1) + \underbrace{\frac{1}{2}\sum_{t=1}^T \alpha_t^2 \|g_t\|_*^2}_{A_T} + \underbrace{\sum_{t=1}^T \alpha_t e_t^\top (x-x_t)}_{B_T}. \label{eqn:combine here}
    \end{aligned}
\end{equation}
As $\mathbb{E}\big[ \exp (\|g_t\|_*^2/G_l^2) \big] \leq \exp (1)$ implies $\mathbb{E}\big[ \|g_t\|_*^2 |x_t \big] \leq G_l^2$, we have 
\begin{align}
    &\|h_t\|_* \leq \|\mathbb{E} \big[ g_t | x_t \big]\|_* \leq \sqrt{\mathbb{E}\big[ \|g_t\|_*^2| x_t \big]} \leq G_l, \nonumber \\
    &\|e_t\|_*^2 = \|g_t - h_t \|_*^2 \leq 4 \max \{ \| g_t \|_*^2, \|h_t \|_*^2 \} \label{eqn:bound of et}.
\end{align}
By \eqref{eqn:bound of et}, we have 
\begin{align}
    &\mathbb{E}[\exp (\|e_t\|_*^2/(2G_l)^2)] \leq \exp(1), \label{eqn:et light tail}\\
    &\mathbb{E}[\exp(A_T/ \sigma_A)] \leq \exp(1), \quad \sigma_A = \frac{1}{2}G_l^2 \sum_{t=1}^T \alpha_t^2 \nonumber
\end{align} 
Let us define $\sigma_A$ as $\sigma_A := \frac{1}{2}G_l^2 \sum_{t=1}^T \alpha_t^2$. Then, by the convexity of the exponential function, we have 
\begin{align}
    \mathbb{E}[\exp(A_T/ \sigma_A)] \leq \frac{1}{2} \sum_{t=1}^T \frac{G_l^2 \alpha_t^2}{\sigma_A} \mathbb{E} \big[ \exp (\|g_t\|_*^2/G_l^2) \big] \leq \frac{1}{2} \sum_{t=1}^T \frac{G_l^2 \alpha_t^2}{\sigma_A} \exp (1) = \exp(1). \nonumber
\end{align}
Using Markov inequality, we have 
\begin{align*}
    \mathbb{P} \Big\{ A_T \geq (1+\Omega) \sigma_A \Big\} \leq \exp(-\Omega). \label{eqn:AT}
\end{align*}

Now, we work on $B_T$. Let us define $\zeta_t = \alpha_t (x - x_t)^\top e_t$. Then, we have
\begin{align*}
    |\zeta_t| \leq \alpha_t \|x_t - x\| \|e_t\|_* \leq \alpha_t \|e_t\|_*\sqrt{2B_{\psi_x}(x_t,x)} \leq \alpha_t\|e_t\|_* \sqrt{2D}.
\end{align*}
By \eqref{eqn:et light tail}, we have
\begin{align}
    \mathbb{E}[\exp (|\zeta_t|^2/(8D\alpha_t^2G_l^2))] \leq \exp(1).
\end{align}
Using martingale technique (see Appendix of \cite{nemirovski2009robust}), we have 
\begin{align}
    \mathbb{P} \Big\{ B_T \geq 2G_l\Omega \sqrt{2\sum_{t=1}^T \alpha_t^2} \Big\} \leq \exp(-\frac{\Omega^2}{4}). \label{eqn:BT}
\end{align}
Combining \eqref{eqn:combine here}, \eqref{eqn:AT}, and \eqref{eqn:BT}, we have 
\begin{align*}
    \mathbb{P} \big\{\sum_{t=1}^T \alpha_t (f_t(x_t) - f_t(x)) \geq D + \frac{1+\Omega}{2}G_l^2 \sum_{t=1}^T \alpha_t^2 + 2G_l \Theta \sqrt{2D \sum_{t=1}^T \alpha_t^2 } \big\} \leq \exp(-\Omega) + \exp(-\Theta^2/4).
\end{align*}
For $\Omega \geq 1$, we have
\begin{align*}
    \mathbb{P} \big\{\sum_{t=1}^T \alpha_t (f_t(x_t) - f_t(x)) \geq D + \Omega G_l^2 \sum_{t=1}^T \alpha_t^2 + 2G_l \Theta \sqrt{2D \sum_{t=1}^T \alpha_t^2 } \big\} \leq \exp(-\Omega) + \exp(-\Theta^2/4).
\end{align*}
Setting $\Theta = 2\sqrt{\Omega}$, inputting diminishing step-size \eqref{eqn:high prob ss}, and dividing by $\sum_{t=1}^T \alpha_t$, for $\Omega \geq 1$, we have 
\begin{align*}
    \mathbb{P} \Big\{\sum_{t=1}^T \theta_t (f_t(x_t) - f_t(x)) \geq \frac{3G_l \cdot \ln (T) \cdot \sqrt{2D_x\Omega}}{\sqrt{T}} \Big\} \leq 2\exp (-\Omega).
\end{align*}
Setting $x$ to $x \in  \underset{x\in X}{\text{argmin}} \ \sum_{t=1}^T \theta_t f_t(x)$, we obtain the desired result. \hfill \Halmos
\endproof

\subsection{Lower Bound of Per-Iteration Sample Size}\label{appn:lower bound}
In Section \ref{subsec:epsilon SMD}, we mentioned that if a per-iteration sample size $K$ has a low dependency on $n$ and satisfies $K \ll n$, for all $t \in [T]$, then Algorithm \ref{alg:Stoc Ax} has much lower per-iteration cost than SMD for solving $\epsilon$-feasibility problem. Indeed, $K$ is a key factor that determines the complexity of $\epsilon$-SMD and the parameters of the approximation condition $C_K$ and $\nu_1$. Our goal in this section is to find $K$ that satisfies the approximation condition with $C_K$ and $\nu_1$, and show that such $K$ is almost independent of $n$. In particular, we show that $K$ has a low dependency on $n$ but is more related to the distribution of random variable $F^i(x,z^i)$, where $z^i \sim p^i$. 
\begin{proposition}\label{prop:K bound Hoeff}
For given $C_K, \nu_1 \in (0,1)$, if we have $K$ such that $K \geq \frac{8M^2}{C_K^2\epsilon^2} \log \frac{2mT}{\nu_1}$, then $K$ satisfies the approximation condition with $C_K$ and $\nu_1$ 
\end{proposition}
\proof{Proof.}
This can be shown using Hoeffding's inequality \citep{wainwright2019high}. Under our assumption, there exists $M^i>0$ such that $|F_r^i(x)| \leq M^i$ for all $r \in [N]$ and $x \in X$ for each constraint $i \in [m]$. Let $M = \underset{i \in [m]}{\max}\ M^i$. Then, by Hoeffding's inequality, we have 
\begin{align}
 \mathbb{P}\big\{|\hat{f}_{t}^i(x_t)-f_t^{i}(x_t)|\geq \frac{C_K \epsilon}{2}\big\} \leq 2\exp^{-\frac{KC_K^2 \epsilon^2}{8M^2}}, \ \forall i \in [m], \ \forall t \in [T] \label{eqn:Hoeffding}
\end{align}
Using union bound and \eqref{eqn:Hoeffding}, we have  
    \begin{align}
    \mathbb{P}\big\{\exists (i,t) \in [m] \times  [T] \text{ s.t. } |\hat{f}_{t}^i(x_t)-f_t^{i}(x_t)|\geq \frac{C_K \epsilon}{2} \big\} \leq \sum_{t=1}^T \sum_{i=1}^m \mathbb{P}\big\{|\hat{f}_{t}^{i}(x_t)-f_t^{i}(x_t)|\geq \frac{C_K \epsilon}{2}\big\}
    \leq 2mT\exp^{-\frac{KC_K^2 \epsilon^2}{8M^2}}.
    \label{eqn:hoeffding temp}
    \end{align}
By \eqref{eqn:hoeffding temp}, for $K$ to satisfy the approximation condition with $C_K$ and $\nu_1$, $K$ must satisfy $2mT\exp^{-\frac{KC_K^2 \epsilon^2}{8M^2}} \leq \nu_1$ and this is equivalent to $ K \geq \frac{8M^2}{C_K^2\epsilon^2} \log \frac{2mT}{\nu_1}$.
\hfill \Halmos
\endproof
For given $C_K \in (0,\frac{1}{2})$ and $ \nu_1 \in (0,1)$, let us define $K_h:=\frac{8M^2}{C_K^2\epsilon^2} \log \frac{2mT}{\nu_1}$. Then $K_h$ satisfies the approximation condition with $C_K$ and $\nu_1$ by Proposition \ref{prop:K bound Hoeff}. Also, we observe that $K_h$ has low dependency on our total iteration $T$, the number of constraints $m$, and probability $\nu_1$ as all these parameters are inside the log term. 

\par 
We derive another condition for $K$ to satisfy the approximation condition using Bennet's inequality \citep{wainwright2019high}.
\begin{proposition}\label{prop:Bennet}
Let $p_t^i \in \calP^i$ be given and $I_t^i$ be an index random sampled from the distribution $p_t^i$ for all $i \in [m]$. Also, we assume that for all $x \in X$, we have 
have,
\begin{align*}
    Var\big(F_{I_t}^i(x)\big) \leq \sigma^2 \, \quad \forall x \in X, i \in [m].
\end{align*}
Then, for given $C_K, \nu_1 \in (0,1)$, if we have $K$ such that $$K \geq \frac{4M^2}{\sigma^2}\frac{1}{h(\frac{2MC_K\epsilon}{\sigma^2})} \log \frac{2mT}{\nu_1},$$ where $h(t) := (1+t) \log (1+t) - t$, for $t \geq 0$, then $K$ satisfies the approximation condition with $C_K$ and $\nu_1$.
\end{proposition}
We omit the proof as it is almost same as the proof of Proposition \ref{prop:K bound Hoeff}. The only difference is that we are using different concentration inequality. We define $K_b$ as $K_b := \frac{4M^2}{\sigma^2}\frac{1}{h(\frac{2MC_K\epsilon}{\sigma^2})}\log \frac{2mT}{\nu_1} $. Notice that $K_b$ is also independent of $n$ and $K_b$ has lower dependency on $\epsilon$ than $K_h$. We want to point out that $K_h$ and $K_b$ are independent of $n$ and more related to the distribution information of the random variable. Also, there is a high chance of $K_b$ being smaller than $n$ if the function $F_r^i(x)$ has a low variance and the $\epsilon$ is set to a moderate level (around 0.01 to 0.001). Although Proposition \ref{prop:K bound Hoeff} and Proposition \ref{prop:Bennet} give us theoretical guarantees of $K$ to satisfy the approximation condition, it is generally hard to find $K$ that satisfies the approximation condition in practice. So, in our experiments, we relax the condition of the approximation condition and use $K$ that satisfies $|\hat{f}_{t}^{i}(x_t)-f_t^{i}(x_T)|\leq \frac{C_K \epsilon}{2}$ for 60\% of the total iterations. Our experimental results (section \ref{sec:experiment} show that we are still able to obtain a solution with a small optimality gap with such $K$.

\section{Proofs for Sections \ref{sec:Generic} and \ref{sec:SOFO}}
\subsection{Proof for Theorem \ref{thm:high prob m>1}}
\proof{Proof.}
For simplicity, we use $\epsilon^{\bullet}$ to denote $\epsilon^{\bullet}(\{x_t,p_t,\theta_t\}_{t=1}^T)$ and use $\epsilon^{i}$ to denote $\epsilon^{i}(\{x_t,p_t^i,\theta_t\}_{t=1}^T)$ for all $i \in [m].$ By Assumption \ref{assm:hpc}, we have 
\begin{align*} 
    \mathbb{P}\Big\{\epsilon^{\bullet}+ \epsilon^{i} \geq w(\Omega)\left(\sR_x(T)+ \sR_p^i(T)\right)+\bar{\epsilon}\Big\} &\leq \mathbb{P}\Big\{  \epsilon^{\bullet} \geq w(\Omega)\sR_x(T) +\bar{\epsilon} \Big\}+ \mathbb{P}\Big\{   \epsilon^{i} \geq w(\Omega)\sR_p^i(T) \Big\}\nonumber\\
    &\leq 4 \exp(-\Omega), \quad  \forall i \in [m].
\end{align*}
Then we have:
\begin{align*}
    \mathbb{P}\Big\{\text{There exists at least one $i \in [m]$ s.t. } \epsilon^{\bullet}+ \epsilon^{i} &\geq w(\Omega)\left(\sR_x(T)+ \sR_p^i(T)\right)+\bar{\epsilon}\Big\} \leq 4m \exp(-\Omega).
\end{align*}
\begin{align}
    \Rightarrow 1- 4m \exp(-\Omega) &\leq \mathbb{P}\Big\{\epsilon^{\bullet}+ \epsilon^{i} \leq w(\Omega)\left(\sR_x(T)+ \sR_p^i(T)\right) + \bar{\epsilon}, \forall i \in [m] \Big\} \nonumber\\
    &\leq \mathbb{P}\Big\{\epsilon^{\bullet}+ \epsilon^{i} \leq w(\Omega)\left(\sR_x(T)+ \underset{i \in [m]}{\max} \ \sR_p^i(T)\right)+\bar{\epsilon}, \forall i \in [m] \Big\} \nonumber \\
    & = \mathbb{P}\Big\{\epsilon^{\bullet}+\underset{i \in [m]}{\max} \ \epsilon^{i} \leq w(\Omega)\left(\sR_x(T)+ \underset{i \in [m]}{\max} \ \sR_p^i(T)\right)+\bar{\epsilon} \Big\}. \nonumber \\
    & = \mathbb{P}\Big\{\epsilon^{\bullet}+\espcirc \leq w(\Omega)\left(\sR_x(T)+ \underset{i \in [m]}{\max} \ \sR_p^i(T)\right) + \bar{\epsilon} \Big\}. \nonumber \hfill \Halmos
\end{align}
\endproof

\subsection{Proof for Lemma \ref{lem:main lemma}}
\proof
{Proof.}
For given $x \in X$, let us define $i^*$ as $i^* \in \underset{i \in [m]}{\text{argmax}} \ f_t^i(x)$. Suppose $\hat{i} \neq i^*.$ Otherwise, we have $h_t \in \partial \phi_t(x) \subseteq \partial_{\bar{\epsilon}} \phi_t(x)$. Because $|\hat{f}_{t}^{i}(x)-f_t^{i}(x)|\leq \frac{\bar{\epsilon}}{2}$ holds for all $i \in [m]$, we have 
\begin{align*}
|f_t^{i^*}(x)-f_t^{\hat{i}}(x)|&\leq |f_t^{i^*}(x) - \hat{f}_{t}^{\hat{i}}(x)| + |f_t^{\hat{i}}(x) - \hat{f}_{t}^{\hat{i}}(x)| \\
&\leq | \underset{i \in [m]}{\text{max}} f_t^i(x) - \underset{i \in [m]}{\text{max}} \hat{f}_{t}^{i}(x)| + \frac{\bar{\epsilon}}{2}\\
& \leq \underset{i \in [m]}{\text{max}} \{|f_t^i(x) - \hat{f}_{t}^{i}(x)|\}+\frac{\bar{\epsilon}}{2} \leq \bar{\epsilon}
\end{align*}
The second inequality comes from the definition of $i^*$ and $\hat{i}$ and our condition $|\hat{f}_{t}^{i}(x)-f_t^{i}(x)|\leq \frac{\bar{\epsilon}}{2}$, for all $i \in [m]$. The third inequality can be derived by $\underset{i \in I}{\max}\ \{a_i\} -\underset{i \in I}{\max} \ \{b_i\} \leq \underset{i \in I}{\max}\{a_i-b_i\}$ for given $\{a_i\}_{i\in I}$ and $\{b_i\}_{i\in I}$. The last inequality comes from our condition $|\hat{f}_{t}^{i}(x)-f_t^{i}(x)|\leq \frac{\bar{\epsilon}}{2}$, for all $i \in [m]$.
Eventually we get 
\begin{align}
    |f_t^{i^*}(x)-f_t^{\hat{i}}(x)| \leq \bar{\epsilon} .\label{eqn:lem1-1}
\end{align} Suppose $h_t$ such that $h_t \in \partial f_t^{\hat{i}}(x)$ is given. Then we have   
\begin{align*}
    f_t^{\hat{i}}(y) \geq f_t^{\hat{i}}(x) + (h_t^{\hat{i}})^\top(y-x), \ \forall y \in X
\end{align*}
By \eqref{eqn:lem1-1}, we have $f_t^{\hat{i}}(x) \geq f_t^{i^*}(x)-\bar{\epsilon}$. Then we have
\begin{align*}
    f_t^{\hat{i}}(y) \geq f_t^{\hat{i}}(x) + (h_t^{\hat{i}})^\top(y-x), \ \forall y \in X \Rightarrow f_t^{\hat{i}}(y) \geq f_t^{i^*}(x) + (h_t^{\hat{i}})^\top(y-x)- \bar{\epsilon}, \ \forall y \in X
\end{align*}
Finally using  $\phi_t(y)=\underset{i \in [m]}{\max} \ f_t^i(y) \geq f_t^{\hat{i}}(y)$ and $\phi_t(x) = f_t^{i^*}(x)$, we have 
\begin{align*}
    \phi_t(y) \geq \phi_t(x) + (h_t^{\hat{i}})^\top(y-x)-\bar{\epsilon}, \ \forall y \in X
\end{align*}
This implies $h_t \in \partial_{\bar{\epsilon}} \phi_t(x)$. \hfill \Halmos
\endproof

\subsection{Bounds related to $p \in \pset$}\label{appn:Dp bound}
Before we prove Theorem \ref{thm:x high prob conv} and Proposition \ref{prop:Rp with entropy}, we drive some bounds that is related to $p \in \pset$.
\par 
Let $p \in \pset$ be given. Then, we have 
\begin{equation}
    \begin{aligned}
        \sum_{r=1}^n p^{(r)} = 1 + \sum_{r=1}^n \big(p^{(r)} - \frac{1}{n} \big) \leq 1 + \sqrt{2n \sum_{r=1}^n  \frac{(p^{(r)}-\frac{1}{n})^2}{2}} = 1 + \sqrt{\frac{2}{n} \sum_{r=1}^n  \frac{(np^{(r)}-1)^2}{2}} \leq 1 + \sqrt{\frac{2\rho}{n}}. \nonumber
    \end{aligned}
\end{equation}
The first inequality comes from Cauchy-Schwarz inequality and the second inequality comes from $ D_{\chi^2}(p\| \frac{\mathbbm{1}_n}{n} )\leq \frac{\rho}{n}$. Let us define $C_g$ as $C_g = 1 + \sqrt{\frac{2\rho}{n}}$. By assuming $\rho \leq \frac{n}{2}$, we have $C_g \leq 2$. In short, we have
\begin{equation}
    \begin{aligned}
        \mathbbm{1}_n^\top p \leq C_g \leq 2.
    \end{aligned}
    \label{eqn:bound of sum of p}
\end{equation}
\par
Also, we derive an upper bound of $D_p: = \underset{p,q \in \pset}{\max} \ B_{\psi_p}(p,q)$. 
First of all, we have $B_{\psi_p}(p,q) = \frac{1}{2}\|p-q\|_2^2$ and for any $p \in \pset$, we have 
    \begin{align}
        B_{\psi_p}(p,\frac{1}{n}\mathbbm{1}_n) = \frac{1}{2}\|p - \frac{1}{n}\mathbbm{1}_n\|_2^2 =\frac{1}{2n^2}\|np - \mathbbm{1}_n\|_2^2 = \frac{D_{\chi^2}(p \| \frac{1}{n}\mathbbm{1}_n)}{n} \leq \frac{\rho}{n^2}.
        \label{eqn:bound of breg p}
    \end{align}
Using \eqref{eqn:bound of breg p}, we have 
\begin{equation}
    \begin{aligned}
            D_p &= \big(\underset{p,q \in \pset}{\max} \sqrt{\ B_{\psi_p}(p,q)} \big)^2 = \big( \underset{p,q \in \pset}{\max} \frac{1}{\sqrt{2}}\|p-q\|_2 \big)^2\\
    &\leq \big( \underset{p,q \in \pset}{\max} \Big\{ \frac{1}{\sqrt{2}}\|p-\frac{1}{n}\mathbbm{1}_n\|_2 +\frac{1}{\sqrt{2}}\|q-\frac{1}{n}\mathbbm{1}_n\|_2\Big\} \big)^2\\
    &=\underset{p,q \in \pset}{\max} \big( \sqrt{B_{\psi_p}(p, \frac{1}{n}\mathbbm{1}_n)} + \sqrt{B_{\psi_p}(q, \frac{1}{n}\mathbbm{1}_n)} \big)^2\\
    &\leq 4n^{-2}\rho
    \end{aligned}
    \label{eqn:bound of dp}
\end{equation}

\subsection{Proof for Theorem \ref{thm:x high prob conv}}\label{appn:proof of thm x}
\proof{Proof.}
In this theorem, as we are discussing the update of $x$, we simplify $g_{x,t}$ to $g_t$ for clarity in the proof.
To show the high-probability convergence of $\epsilon$-SMD, 
there are two major challenges. First,
we need to show that the light tail assumption (Assumption \ref{assm:light tail}) holds for the update gradient $g_t$. Secondly, we need to extend Proposition \ref{prop:Nemi high prob conv} to the case where we use $\epsilon$-subgradients.
\par 
The first thing can be done by finding a uniform bound of $\|g_t\|_{x,*}$. That is, if there exists $G_l>0$ such that $\|g_t\|_{x,*} \leq G_l$ for a given random vector $g_t$, then our gradient estimator $g_t$ satisfies the light tail assumption with constant $G_l$.
\par 
Let $i \in [m], r \in [n]$, $x \in X, p^i \in \pset$ and $g_{x}$ such that $g_{x} \in \partial F_r^{i}(x)$ be given.
As $F_r^i(x)$ is a $G$-Lipshcitz function from our assumption, this implies for any $g_x$ such that $g_x \in \partial  F_r^{i}(x)$  satisfies
\begin{align*}
    \|g_x\|_{x,*} \leq G.
\end{align*}
As a result, the light tail assumption holds for $g_t$ with constant $G$ for all $t \in [T]$.
\par 
Now we extend the proof of Proposition \ref{prop:Nemi high prob conv} to when $h_t$ is a $\bar{\epsilon}$-subgradient of $f_t(x)$ at $x_t$ by adding $\alpha_t \bar{\epsilon}$ to the right-hand side of \eqref{eqn:prop nemi eqn1}.
Similary to the proof of Proposition \ref{prop:Nemi high prob conv}, define $h_t$ as $h_t := \mathbb{E}\big[g_t |x_t \big]$ and $e_t$ as $e_t = g_t - h_t$.
By the approximation condition (Definition \ref{def:approx}), we condition on when \eqref{eqn:K cond} holds and let $\bar{\epsilon}$ be defined as $\bar{\epsilon}=C_K\epsilon$. As $h_t \in \partial_{\bar{\epsilon}} \phi_t(x_t)$ by Lemma \ref{lem:main lemma}, for any $x \in X$ and $t \in [T]$,  we have 
\begin{align*}
    \alpha_t(\phi_t(x_t)  - \phi_t(x)) &\leq \alpha_t h_t^\top (x_t - x) + \alpha_t \bar{\epsilon} \nonumber\\
    &= \alpha_t(h_t - g_t)^\top (x_t - x) + \alpha_t g_t^\top (x_t - x) + \alpha_t \bar{\epsilon}\nonumber\\
    &\leq \alpha_te_t^\top (x- x_t) + B_{\psi_x}(x,x_t) - B_{\psi_x}(x,x_{t+1}) +\frac{\alpha_t^2\|g_t\|_*^2}{2} + \alpha_t \bar{\epsilon}.
\end{align*}
The rest of the proof follows the same as the proof of Proposition \ref{prop:Nemi high prob conv}. Eventually, we get 
\begin{align*}
    \mathbb{P} \Big\{ \espbullet \geq \frac{3G \cdot \ln (T) \cdot \sqrt{2D_x\Omega}}{\sqrt{T}} + C_K\epsilon | \eqref{eqn:K cond} \text{ is true} \Big\} \leq 2\exp (-\Omega).
\end{align*}
This implies 
\begin{align*}
    \mathbb{P} \Big\{\espbullet \leq \frac{3G \cdot \ln (T) \cdot \sqrt{2D_x\Omega}}{\sqrt{T}} + C_K\epsilon | \eqref{eqn:K cond} \text{ is true} \Big\} \geq 1 - 2\exp (-\Omega).
\end{align*}
Using the approximation condition, we have 
\begin{align*}
    \mathbb{P} \Big\{\espbullet \geq \frac{3G \cdot \ln (T) \cdot \sqrt{2D_x\Omega}}{\sqrt{T}} + C_K\epsilon \Big\} &\geq (1-2\exp (-\Omega))(1-\nu_1)\\
    &\geq 1-2\exp(-\Omega) -\nu_1.
\end{align*}
The first inequality comes from 
$$\mathbb{P}\big\{ A \big\} \geq \mathbb{P}\big\{ A \cap B \big\} = \mathbb{P}\big\{ A | B \big\} \mathbb{P}\big\{ B \big\} $$
for given events $A$ and $B$. \hfill \Halmos
\endproof

\subsection{Proof for Lemma \ref{lemma:vartheta}}
\proof{Proof.}
Under the approximation condition, we have
 \begin{align*}
     |\vartheta^i - \hvt^i| \leq \frac{C_K\epsilon}{2} \leq \frac{C_K\epsilon}{2}, \quad \forall i \in [m]
 \end{align*}
 with probability at least $1-\nu_0$. Then following the similar technique that we used in the proof of Lemma \ref{lem:main lemma}, we get 
 \begin{align*}
    |\vartheta - \hvt| \leq C_K\epsilon 
\end{align*}
with probability at least $1- \nu_0$. \hfill \Halmos
\endproof

\subsection{Proof for Theorem \ref{thm:feas}}
\proof{Proof.}
Let $\kcirc = \underset{i \in [m]}{\max} \ \sR_p^i(T_\epsilon) /\epsilon$.
According to Proposition \ref{prop:Nam 3.2}, if $\vartheta \leq (1- \kcirc)\epsilon$, then $\bar{x}_T$ becomes $\epsilon$-feasible certificate of $\epsilon$-feasibility problem and if $\vartheta > \kbullet \epsilon$, then we declare that the problem is infeasible. Now, we consider the following cases.
\begin{enumerate}
    \item $\hvt \geq (\kbullet + C_K)\epsilon$:
    This implies $\vartheta \geq \kbullet \epsilon$ by \eqref{eqn:vartheta and hvt diff}, and we declare that the problem is infeasible.
    \item $\hvt \leq (\kbullet + C_K)\epsilon$:
    Then we have $\vartheta \leq (\kbullet + 2C_K) \epsilon \leq (1-\kcirc)\epsilon.$ The second inequality comes from \eqref{eqn:bound sum of kappa}. So $\bar{x}$ is an $\epsilon$-feasible certificate. 
\hfill \Halmos
\end{enumerate}
\endproof

\subsection{Proof for Proposition \ref{prop:Rp with entropy}}
\proof{Proof.}
Similar to the proof of Theorem \ref{thm:x high prob conv}, we first need to find a constant that satisfies the light tail assumption for the gradient estimator $g_{p,t}$ by finding a uniform bound of $\|g_{p,t}^i\|_{p,*}$. Let us assume $x \in X$, $i \in [m]$, $I_p \in [n]$ and $p^i \in \pset$ be given. Also, suppose  $g_p^i \in \mathbb{R}^n$ such that $g_p^{i^{(r)}} = \frac{(\mathbbm{1}_n^{\top}p^i)F_r^i(x)}{p^{i^{(r)}}} \mathbbm{1}(r = I_p)$ for all $r \in [n]$, be given. As $g_p^i$ is a 1-sparse vector and $|F_r^i(x)|\leq M^i$ and $p^i \geq \frac{\delta}{n}\mathbbm{1}_n$, we easily get
\begin{align*}
    \|g_p^i\|_{p,*} \leq \frac{C_gM^in}{\delta}
\end{align*}
by \eqref{eqn:bound of sum of p} and this implies that the gradient estimator $g_{p,t}^i$ satisfy the light tail assumption with constant $\frac{C_gM^in}{\delta}$, for all $i \in [m]$. We also have upper bound of $D_p$ by \eqref{eqn:bound of dp}. We get the desired result by inputting $G_l =  \frac{C_gM^in}{\delta}$ and $D_p = \frac{4\rho}{n^2}$ to Proposition \ref{prop:Nemi high prob conv}. \hfill \Halmos
\endproof

\subsection{Proof for Theorem \ref{thm:alg1 conv}}
\proof{Proof.}
As $\mathbbm{1}_n^\top p^i = 1$ does not hold for $p^i \in \pset$, we need to slightly modify Theorem \ref{thm:x high prob conv}.
The only thing that is changed is the constant $G_l$ that satisfies the light tail assumption (Assumption \ref{assm:light tail}) for gradient estimator $g_{x,t}$ in Algorithm \ref{alg:Stoc Ax}. Let $i \in [m], r \in [n]$, $x \in X, p^i \in \pset$ and $g_{x}$ such that $g_{x} \in \partial (\mathbbm{1}_n^\top p^i) F_r^{i}(x)$ be given.
As $\mathbbm{1}_n^\top p^i = 1$ does not hold for $p^i \in \pset$ anymore, $I_t^i$ is now random sampled from the normalized distribution $\frac{p_t^i}{\mathbbm{1}_n^\top p_t^i}$ and we modify our function $F_r^i(x)$ to $(\mathbbm{1}_n^\top p_t^i)F_r^i(x)$ according to this change.
As $F_r^i(x)$ is a $G$-Lipshcitz function from our assumption, $\mathbbm{1}_n^{\top}p^i \cdot F_r^{i}(x)$ is a $C_gG$-Lipschitz function with $C_g = 1+ \sqrt{\frac{2\rho}{n}}$ by \eqref{eqn:bound of sum of p}.
This implies for any $g_x$ such that $g_x \in \partial (\mathbbm{1}_n^{\top}p^i) F_r^{i}(x)$  satisfies
\begin{align*}
    \|g_x\|_{x,*} \leq C_g G.
\end{align*}
As a result, the light tail assumption holds for $g_{x,t}$ with constant $C_gG$ for all $t \in [T]$. So replacing $G$ with $C_gG$ in Theorem \ref{thm:x high prob conv}, we obtain the high-probability convergence result of $\epsilon$-SMD under $\chi^2$-divergence set. In fact, for any $\Omega \geq 1$, we have
\begin{equation}
    \begin{aligned}
    \mathbb{P}\Big \{\epsilon^{\bullet}(\{x_t,p_t,\theta_t\}_{t=1}^T) \leq \frac{3C_g\cdot G \cdot \ln (T) \cdot \sqrt{2D_x\Omega}}{\sqrt{T}}+ C_K \epsilon \Big\}
     &\geq (1 - 2\exp (- \Omega))(1-\nu_1) \\
     &\geq 1 - 2\exp (-\Omega)-\nu_1.
    \end{aligned}
    \label{eqn:prop3.7 eqn}
\end{equation}
After combining \eqref{eqn:prop3.7 eqn}, Theorem \ref{thm:high prob m>1} and Proposition \ref{prop:Rp with entropy}, we get the desired result. \hfill \Halmos
\endproof

\section{Efficient Updates for Bandit Mirror Descent}\label{appn:efficient update}
In this section, we study the per-iteration complexity of Algorithm \ref{alg:p under chi}. Even our gradient estimator $g_{p,t}^i$ is 1-sparse, complexity of the projection step of Algorithm \ref{alg:p under chi} could be $O(n)$. \citep{namkoong2016stochastic} show that the complexity of Algorithm \ref{alg:p under chi} is $O(\log (n)+ \log(\frac{1}{\epsilon_g}))$ by using a balanced tree structure, where $\epsilon_g$ is a tolerance level that we define shortly. In this section, we present a slightly improved version of their work. We show that we can perform each step of Algorithm \ref{alg:p under chi} in $O(\log n)$.
\par 
First of all, we can implement the index random sampling step (line 3) in $O(\log(n) )$ using the balanced tree structure. We refer our readers to \citep{cormen2009introduction}, \citep{duchi2008efficient}, \citep{namkoong2016stochastic} for the detail implementation of the balanced tree structure. Also, the computation of the gradient estimator $g_{p,t}^i$ and the update of $w_{t+1}$ can be implemented in a constant time because $g_{p,t}^i$ is 1-sparse. In the rest of this section, we show that the complexity of the projection step is $O(\log \frac{1}{\epsilon_g})$.
\par 
As the projection step is identical for all the constraints and iterations, we omit the constraint index $i$ and the iteration index $t$ unless it is confusing. Initially, we define Lagrangian function $\mathcal{L}(p; \lambda, \theta)$ as
\begin{equation}
    \begin{aligned}
        \mathcal{L}(p; \lambda, \theta) = \frac{1}{2}\|p - w \|_2^2 - \frac{\lambda}{n^2} \Big( \rho  - \sum_{r=1}^n \frac{(np^{(r)}-1)^2}{2} \Big) - \theta^{\top}(p - \frac{\delta}{n}\mathbbm{1}_n) \label{eqn:Lagrangian0}
    \end{aligned}
\end{equation}
where $\theta \in \mathbb{R}_+^{n}$ and $\lambda \geq 0 $. The problem $\underset{p\in \mathcal{P}_{\rho,n,\delta}}{\text{min}} \ B_{\psi_p}(p,w_{t+1})$ satisfies Slater condition as $p = \frac{1}{n}\mathbbm{1}_n$ is an interior point of $\pset$. Also, $B_{\psi_p}(p,w) = \frac{1}{2}\|p-w\|_2^2$ is a strongly convex function on $p$. So, using KKT conditions and strict complementary slackness, we get 
\begin{equation}
    \begin{aligned}
        p(\lambda) = \Big( \frac{1}{1+\lambda}w + \frac{\lambda}{1+ \lambda}\frac{1}{n}\mathbbm{1}_n - \frac{\delta}{n}\mathbbm{1}_n \Big)_+ + \frac{\delta}{n}\mathbbm{1}_n
    \end{aligned}
\end{equation}
where $p(\lambda) = \text{argmin}_{p \in \pset} \text{sup}_{\theta \in \mathbb{R}_+^n} \mathcal{L}(p; \lambda, \theta)$. Substituting $p(\lambda)$ into equation \eqref{eqn:Lagrangian0}, we have 
    \begin{align*}
        g(\lambda) = \frac{1}{2}\|p(\lambda) - w\|_2^2 - \frac{\lambda}{n^2}\Big(\rho -  \sum_{r=1}^n \frac{(np(\lambda)^{(r)}-1)^2}{2} \Big).
    \end{align*}
Our primary goal is to find $\lambda^*$ such that $\lambda^* = \underset{\lambda \geq 0}{\text{argmax}} \ g(\lambda)$ and input $\lambda^*$ to $p(\lambda)$ to obtain $p_{t+1}$.
Let us define a set $I(\lambda) = \{ r \in [n]| \frac{1}{1+\lambda}w^{(r)} + \frac{\lambda}{1+\lambda}\frac{1}{n} \geq \frac{\delta}{n}\} $. Then $g'(\lambda)$ is like: 
    \begin{align*}
         \frac{\partial}{\partial \lambda} g(\lambda) = &\frac{1}{2(1+\lambda)^2} \sum_{r \in I(\lambda)} w^{(r)^2} - \frac{1}{n(1+\lambda)^2} \sum_{r \in I(\lambda)} w^{(r)} + \Big(\frac{1}{2n^2(1+\lambda)^2} - \frac{(1-\delta)^2}{2n}\Big)|I(\lambda)|\\
         &+ \frac{(1-\delta)^2}{2n} - \frac{\rho}{n^2}.
    \end{align*}
Suppose $\lambda^*\geq 0$ such that $g'(\lambda^*) = 0$ exists. As $g(\lambda)$ is a concave function of $\lambda$, $g'(\lambda)$ is a non-increasing function of $\lambda$. So, for given $\epsilon_g >0$, we can find $\lambda^*$ such that $|g'(\lambda^*)| \leq \epsilon_g$ by performing $O(\log \frac{1}{\epsilon_g})$ binary search iterations.
We show that we can update the values of $\sum_{r \in I(\lambda)} w^{(r)^k}$ for $k=0,1,2$ in a constant time in the Appendix \ref{appn:update of w} and \ref{appn:update of p}.
As a result, finding $\lambda^*$ takes $O(\log \frac{1}{\epsilon_g})$ steps. However, there may be cases where no $\lambda^* \geq 0$ satisfying $g'(\lambda^*)=0$ exists, and the termination condition, $|g'(\lambda)|\leq \epsilon_g$, is never met. 
We address this issue in the Appendix \ref{subsec:Calculation of Lambda} and show that we can find $\lambda^* \in \underset{\lambda \geq 0}{\text{argmax}} \ g(\lambda)$, in $O(\log \frac{1}{\epsilon_g})$ and even in a constant time for some cases. In summary, the complexity of obtaining $\lambda^*$ is $O(\log n)$ as $O(\log \frac{1}{\epsilon_g}) \approx O(\log n)$.
\par
What remains is to determine the complexity of calculating $p_{t+1}(\lambda^*)$ after obtaining $\lambda^*$. Indeed, we can obtain  $p_{t+1}(\lambda^*)$ in a constant time by using the 1-sparsity of the gradient estimator $g_{p,t}$
Suppose $I_p$ be an index that is randomly chosen to update $p_t$ in Algorithm \ref{alg:p under chi}. By the 1-sparsity of the gradient estimator $g_{p,t}$, $w_{t+1}^{(r)} = p_t^{(r)} \geq \frac{\delta}{n}$ for all $r \in [n] \backslash \{I_p\}$.
Then this implies that $\frac{1}{1+\lambda}w_{t+1} + \frac{\lambda}{1+ \lambda}\frac{1}{n}\mathbbm{1}_n \geq \frac{\delta}{n}\mathbbm{1}_n$, for any $\lambda \geq 0$ and we have $p_{t+1}(\lambda)^{(r)} = \frac{1}{1+\lambda}w_{t+1}^{(r)} + \frac{\lambda}{1+ \lambda}\frac{1}{n}$ for all $r \in [n] \backslash \{I_p\}$.
Notice that all elements of $p_{t+1}(\lambda)$ except the randomly chosen index $I_p$ is updated by the same multiplicative factor $\frac{1}{1+ \lambda}$ and additive factor $\frac{\lambda}{n(1+\lambda)}$.
As a result, we can obtain $p_{t+1}(\lambda^*)$ in a constant time by merely updating the multiplicative and the additive factors and appropriately adjusting the value of $p_{t+1}^{(I_p)}$. 
Putting all the pieces together, we conclude that we can efficiently update $p_t$ whose complexity is almost independent to $n$.
\subsection{Efficient Calculation of $\lambda^*$}\label{subsec:Calculation of Lambda}
Here, we present an efficient method to get $\lambda^*$ such that $\lambda^* \in \underset{\lambda \geq 0}{\text{argmax}} \ g(\lambda)$. After the change of variables $\alpha = \frac{\lambda}{1+ \lambda}$ and $I(\alpha):= \{r \in [n]: (1-\alpha)w^{(r)} + \frac{\alpha}{n}\geq \frac{\delta}{n}\}$, $g(\alpha)$ is still a concave function of $\alpha$ and $g'(\alpha)$ is like:
\begin{equation}
    \begin{aligned}
    \frac{\partial}{\partial \alpha} g(\alpha) = &\frac{1}{2} \sum_{r \in I(\alpha)} w^{(r)^2} - \frac{1}{n}\sum_{r \in I(\alpha)} w^{(r)} + \frac{1}{2n^2(1-\alpha)^2} ((1-\alpha)^2-(1-\delta)^2)|I(\alpha)|\\
    &+\frac{1}{2n^2(1-\alpha)^2} (n(1-\delta)^2 - 2\rho)
    \end{aligned}
    \label{eqn:g'(alpha)}
\end{equation}
Now suppose there exists $\alpha^* \in [0,1]$ such that $g'(\alpha^*) = 0$. Then, for given $\epsilon_g >0$, using bisection search on non-increasing function $g'(\alpha)$, we can find $\alpha^*$ such that $|g'(\alpha^*)| \leq \epsilon_g$ with $O(\log \frac{1}{\epsilon_g})$ bisection search. In later section, we show that we can obtain $\sum_{r \in I(\alpha)} w_{t+1}^{(r)^k}$ for $k = 0, 1, 2$ in a constant time by using the information of $\sum_{r=1}^n p_{t}^{(r)^k}$, $k=0,1,2$. As a result, we can update $p_t$ in $O(\log \frac{1}{\epsilon_g})$, when there exists $\alpha^* \in [0,1]$ such that $g'(\alpha^*) = 0$.
\par 
Let us consider the case when $\alpha^*$ that satisfies $g'(\alpha^*)=0$ does not exist. In this case, we need other termination conditions to find $\alpha^*$ such that $\alpha^* \in \underset{\alpha \in [0,1]}{\text{argmax}} \ g(\alpha)$.
First of all, nonexistence of $\alpha^*$ such that $g'(\alpha^*) = 0$ implies that $g'(\alpha)$ is either non-negative or non-positive for all $\alpha \in [0,1]$.
If $g'(\alpha)<0$ for all $\alpha \in [0,1]$, then $0 = \underset{\alpha \in [0,1]}{\text{argmax}} \ g(\alpha)$. We can simply check whether the condition, $g'(\alpha)<0$ for all $\alpha \in [0,1]$, is satisfied by checking the value of $g'(0)$, as $g'(\alpha)$ is a non-increasing function. 
So, for sufficiently small $\epsilon_{\alpha}>0$, if $\alpha < \epsilon_{\alpha}$ and $g'(\alpha)<0$, we terminate our bisection search and return 0.
Also, we can easily check that $\lim_{\alpha \rightarrow 1} I(\alpha) = [n]$ and $\lim_{\alpha \rightarrow 1} \ g'(\alpha) = -\infty$ and this implies $g'(\alpha)>0$ for all $\alpha \in [0,1]$ does not hold.
Notice that we now have two different termination conditions.
One is to terminate if we find $\alpha$ such that $|g'(\alpha)| < \epsilon_g$ and the other is to terminate if we find $\alpha$ such that $\alpha < \epsilon_{\alpha}$ and $g'(\alpha)<0$.
Based on these termination conditions and the virtue of updates of Algorithm \ref{alg:p under chi}, we present an efficient method to obtain $\alpha^*$ such that $\alpha^* \in \underset{\alpha \in [0,1]}{\text{argmax}} \ g(\alpha)$, even in a constant time in some cases.
\par 
For given index $r \in [n]$, if $w^{(r)}$ satisfies $w^{(r)} \geq \frac{\delta}{n}$, then $r$ belongs to the set $I(\alpha)$ regardless of the value of $\alpha$. As $g_{p}^i$ is a 1-sparse vector and $p^{(r)} \geq \frac{\delta}{n}$, for all $r \in [n]$, by the virtue of the update (Algorithm \ref{alg:p under chi}.5), $p^{(r)} = w^{(r)}\geq \frac{\delta}{n}$ for $r \neq I_{p}$, where $I_{p}$ is a randomly selected index from the distribution $\frac{p}{\mathbbm{1}_n^\top p}$.

This implies $[n]\backslash\{I_{p}\} \subseteq I(\alpha)$ is always satisfied and we only need to check whether $I_{p}$ belongs to the set $I(\alpha)$ or not. If $w^{(I_{p})} = p_t^{(I_{p})} + \alpha_{p}g_{p} \geq \frac{\delta}{n}$, then the set $I(\alpha)$ is always equal to $[n]$ regardless of the value of $\alpha$. Otherwise, $I(\alpha)$ depends on $\alpha$.
If $w^{(I_{p})} \geq \frac{\delta}{n}$, $g'(\alpha)$ can be expressed as 
\begin{equation}
    \begin{aligned}
         \frac{\partial}{\partial \alpha} g(\alpha) = & \frac{1}{2}\sum_{r=1}^n (w^{(r)}- \frac{1}{n})^2 - \frac{\rho}{n^2 (1-\alpha)^2}.
    \end{aligned}
    \label{eqn:appdx alpha eq1}
\end{equation}
Suppose $\frac{1}{2}\sum_{r=1}^n (w^{(r)}- \frac{1}{n})^2 \geq \frac{\rho}{n^2}$ holds. Then $\frac{1}{2}\sum_{r=1}^n (w^{(r)}- \frac{1}{n})^2 \geq \frac{\rho}{n^2}$ implies $g'(0) \geq 0$ and there exists $\alpha^*$ such that $|g'(\alpha^*)| \leq \epsilon_g$ as $g'(1) < 0$. Using \eqref{eqn:appdx alpha eq1}, we have closed form solution for $\alpha^*$ as the following:
\begin{equation}
    \begin{aligned}
     \frac{\partial}{\partial \alpha} g(\alpha^*) = 0 \iff \alpha^* = 1 - \frac{1}{n}\sqrt{\frac{\rho}{(\frac{1}{2}\sum_{r=1}^n (w^{(r)}- \frac{1}{n})^2)}}
    \end{aligned}
    \label{eqn:alpha closed form1}
\end{equation}
In this case, if we can get the value of $\sum_{r=1}^n w_{t+1}^{(r)^k}$, for $k=0,1,2$, in a constant time, then we get the $\alpha^*$ in a constant time without any bisection search by \eqref{eqn:alpha closed form1}. Indeed, if we have the information of $\sum_{r=1}^n p^{(r)^k}$, for $k=0,1,2$, we are able to calculate $\sum_{r=1}^n w^{(r)^k}$, for $k=0,1,2$, in a constant time, as only one element has a different value between $p$ and $w$. On the other hand, if $\frac{1}{2}\sum_{r=1}^n (w^{(r)}- \frac{1}{n})^2 < \frac{\rho}{n^2}$, then this implies $g'(0)<0$ and we return 0. In conclusion, when $I(\alpha) = [n]$, we get $\alpha^*$ in a constant time.
\par 
Now we consider the case when $w^{(I_p)} < \frac{\delta}{n}$. Then, whether $I_p \in I(\alpha)$ or not, depends on $\alpha$. Let us define $\bar{\alpha}$ as $\bar{\alpha} = \frac{\delta/n - w^{(I_p)}}{1/n- w^{(I_p)}}$. If $\alpha \geq \bar{\alpha}$, then $I_p \in I(\alpha)$ and otherwise, $I_p \notin I(\alpha)$. As $g'(\alpha)$ is a non-increasing function, we implement our regular bisection search for $\alpha^*$ and update a search interval for $\alpha$ depending on the $g'(\alpha)$ value. Suppose we denote our current search interval as $[l_{\alpha}, u_{\alpha}]$. If $l_{\alpha}> \bar{\alpha}$, then this implies $I(\alpha) = [n]$ for all $\alpha$ that we will going to search in the future. 
This again means that $g'(\alpha)$ would follow the expression \eqref{eqn:appdx alpha eq1} and we can apply the same procedure as when $w^{(I_p)}\geq \frac{\delta}{n}$. If $u_{\alpha} < \bar{\alpha}$, then this implies $I(\alpha) = [n] \backslash \{I_{p}\}$ for all the $\alpha$ that we are going to search in the future. With $I(\alpha) = [n]\backslash \{I_p\}$, $g'(\alpha)$ can be expressed as 
\begin{align*}
    g'(\alpha) = \frac{1}{2}\sum_{r \in [n] \backslash \{I_{p}\}} w^{(r)^2} - \frac{1}{n}\sum_{r \in [n] \backslash \{I_{p}\}} w^{(r)} + \frac{n-1}{2n^2} - \frac{1}{(1-\alpha)^2}(\frac{\rho}{n^2} - \frac{(1-\delta)^2}{2n^2}).
\end{align*}
If $\frac{1}{2}\sum_{r \in [n] \backslash \{I_{p}\}} w^{(r)^2} - \frac{1}{n}\sum_{r \in [n] \backslash \{I_{p}\}} w^{(r)} + \frac{n-1}{2n^2} \leq \frac{\rho}{n^2} - \frac{(1-\delta)^2}{2n^2}$, then this implies $g'(\alpha)\leq 0$, so we set $\alpha^*$ to 0. Otherwise, we have a closed form solution for $\alpha^*$:
\begin{equation}
    \begin{aligned} \alpha^* = 1- \sqrt{\frac{\frac{\rho}{n^2} - \frac{(1-\delta)^2}{2n^2}}{\frac{1}{2}\sum_{r \in I(\alpha)} w^{(r)^2} - \frac{1}{n}\sum_{r \in I(\alpha)} w^{(r)} + \frac{n-1}{2n^2}}}.
    \end{aligned}
    \label{eqn:appdx alpha eq3}
\end{equation}
As we can calculate $\sum_{r \in I(\alpha)} w_{t+1}^{(r)^k}$, for $k=0,1,2$, in a constant time, \eqref{eqn:appdx alpha eq3} can be calculated in a constant time. Notice that the algorithm for finding $\alpha$ would not terminate if an interval $[l_{\alpha}, u_{\alpha}]$ always includes $\bar{\alpha}$. However, this implies our interval would center around $\bar{\alpha}$ and eventually terminate by the other termination condition, which is $|g'(\alpha)| \leq \epsilon_g$. So, the termination of our algorithm is guaranteed. Algorithm \ref{alg:find alpha} is a pseudo-code for finding $\alpha^*$.
\begin{algorithm}[h]
\begin{algorithmic}[1]
\STATE \textbf{input} value of $\sum_{r=1}^n p_{t}^{(r)^k}$ for $k=0,1,2$, $w^{(I_p)}$
\STATE \textbf{output} $\alpha^*$ such that $\alpha^* \in \underset{\alpha \in [0,1]}{\text{argmax}} \ g(\alpha)$
\IF{$w^{(I_p)}\geq \frac{\delta}{n}$}
    \IF{ $\frac{1}{2}\sum_{r=1}^n (w^{(r)}- \frac{1}{n})^2 \geq \frac{\rho}{n^2}$}
        \STATE \textbf{return} $\alpha^* = 1 - \frac{1}{n}\sqrt{\frac{\rho}{(\frac{1}{2}\sum_{r=1}^n (w^{(r)}- \frac{1}{n})^2)}}$
    \ELSE
        \STATE \textbf{return} $\alpha^* = 0$
    \ENDIF
\ELSE
    \STATE Set $l_{\alpha}=0, u_{\alpha}=1, \bar{\alpha} = \frac{\delta/n - w^{(I_p)}}{1/n- w^{(I_p)}}$.
    \WHILE{True}
    \IF{$l_{\alpha}> \bar{\alpha}$}
        \IF{ $\frac{1}{2}\sum_{r=1}^n (w^{(r)}- \frac{1}{n})^2 \geq \frac{\rho}{n^2}$}
        \STATE \textbf{return} $\alpha^* = 1 - \frac{1}{n}\sqrt{\frac{\rho}{(\frac{1}{2}\sum_{r=1}^n (w^{(r)}- \frac{1}{n})^2)}}$
        \ELSE
        \STATE \textbf{return} $\alpha^* = 0$
        \ENDIF
    \ENDIF
    \STATE Algorithm continues in the next page
            \caption{Computing $\alpha^*$ that maximizes $g(\alpha)$}
      \algstore{SOFO}
      \label{alg:find alpha}
        \end{algorithmic}
        \end{algorithm}
        \begin{algorithm}[h]
        \begin{algorithmic}[1]
        \algrestore{SOFO}
    \IF{$u_{\alpha}<\bar{\alpha}$}
        \IF{$\frac{1}{2}\sum_{r \in [n] \backslash \{I_{p}\}} w^{(r)^2} - \frac{1}{n}\sum_{r \in [n] \backslash \{I_{p}\}} w^{(r)} + \frac{n-1}{2n^2} \leq \frac{\rho}{n^2} - \frac{(1-\delta)^2}{2n^2}$}
        \STATE \textbf{return} $\alpha^* = 0$
        \ELSE
            \STATE \textbf{return} $\alpha^* = 1- \sqrt{\frac{\frac{\rho}{n^2} - \frac{(1-\delta)^2}{2n^2}}{\frac{1}{2}\sum_{r \in I(\alpha)} w^{(r)^2} - \frac{1}{n}\sum_{r \in I(\alpha)} w^{(r)} + \frac{n-1}{2n^2}}}$
        \ENDIF
    \ENDIF

    \STATE Set $\alpha = \frac{l_{\alpha}+ u_{\alpha}}{2}$
    \STATE Compute $g'(\alpha)$
    \IF{$|g'(\alpha)| \leq \epsilon_{g}$}
    \STATE \textbf{return} $\alpha$
    \ENDIF
    \IF {$g'(\alpha)<0$}
          
    \STATE Set $u_{\alpha} = \alpha$
    \ELSE 
    \STATE Set $l_{\alpha} = \alpha$
    \ENDIF
    \ENDWHILE
\ENDIF
\end{algorithmic}
\end{algorithm}

\subsection{Update of $\sum_{r \in I(\alpha)}
w_{t+1}^{(r)^k}$, for $k=0,1,2$}\label{appn:update of w}
Here, we describe a calculation step, which calculates $\sum_{r \in I(\alpha)} w_{t+1}^{(r)^k}$, for $k=0,1,2$ based on the information, $\sum_{r=1}^n p_{t}^{(r)^k}$, $k=0,1,2$. If $w_{t+1}^{(I_{p,t})} \geq \frac{\delta}{n}$, then update of $\sum_{r \in I(\alpha)} w_{t+1}^{(r)^k}$ is trivial because $I(\alpha) = [n]$ for all $\alpha \in [0,1]$. So, we only consider when $w_{t+1}^{(I_{p,t})} < \frac{\delta}{n}$ here. Let us define $\bar{\alpha}$ as $\bar{\alpha} = \frac{\delta/n-w_{t+1}^{(I_{p,t})}}{1/n-w_{t+1}^{(I_{p,t})}}$. If $\alpha > \bar{\alpha}$, then $I_{p,t} \in I(\alpha)$ and we have $\sum_{r \in I(\alpha)} w_{t+1}^{(r)^k} =\sum_{r=1}^n p_{t}^{(r)^k} + (w_{t+1}^{(I_{p,t})^k} - p_t^{(I_{p,t})^k})$, for $k=0,1,2$. Otherwise, $I_{p,t} \notin I(\alpha)$ and we have $\sum_{r \in I(\alpha)} w_{t+1}^{(r)^k} = \sum_{r=1}^n p_{t}^{(r)^k} - p_t^{(I_{p,t})^k}.$ This shows that we can get the value of $\sum_{r \in I(\alpha)} w_{t+1}^{(r)^k}$ in a constant time.

\subsection{Update of $\sum_{r=1}^n p_{t+1}^{(r)^k}$, for $k=1,2$}\label{appn:update of p}

Here, we describe a calculation step, which calculates $\sum_{r=1}^n p_{t+1}^{(r)^k}$ using $\sum_{r=1}^n p_{t}^{(r)^k}$, for $k=1,2$. When $k=0$, the value of $\sum_{r=1}^n p_{t+1}^{(r)^k}$ is always equal to $n$. We calculate this in a constant time by using  the update relationship  $p_{t+1}^{(r)} = (1-\alpha)p_t^{(r)} + \frac{\alpha}{n}$ for $r \in [n] \backslash I_{p,t}$. For randomly sampled index $I_{p,t}$, we need to check whether $(1-\alpha)w_{t+1}^{(I_{p,t})} + \frac{\alpha}{n} \geq \frac{\delta}{n}$ or not. If $(1-\alpha)w_{t+1}^{(I_{p,t})} + \frac{\alpha}{n} \geq \frac{\delta}{n},$ then this implies that $I_{p,t} \in I(\alpha)$ and we have $p_{t+1}^{(I_{p,t})} = (1-\alpha)p_t^{(I_{p,t})} + \frac{\alpha}{n}$. Otherwise, $p_{t+1}^{(I_{p,t})}$ is set to $\frac{\delta}{n}$. Using all of these relationships, we obtain $\sum_{r=1}^n p_{t+1}^{(r)^k}$, for $k=1,2$ in a constant time. For simplicity, let us denote $\sum_{r=1}^n p_{t}^{(r)}$ as $\beta_t$ and denote $\sum_{r=1}^n p_{t}^{(r)^2}$ as $\gamma_t$. Let us first start with when $(1-\alpha)w_{t+1}^{(I_{p,t})} + \frac{\alpha}{n} \geq \frac{\delta}{n}$.  Then, we have $p_{t+1} = (1-\alpha)p_t + \frac{\alpha}{n}$. However, notice that we need to replace $p_t^{(I_{p,t})}$ to $w_{t+1}^{(I_{p,t})}$. So, we update the values of $\beta_t$ and $\gamma_t$ according to this change ($\beta_t \leftarrow \beta_t + w_{t+1}^{(I_{p,t})} - p_t^{(I_{p,t})}, \gamma_t \leftarrow \gamma_t +w_{t+1}^{(I_{p,t})^2} - p_t^{(I_{p,t})^2}$). As, all the elements share same multiplication factor and additional factor, we have
\begin{align*}
    &\sum_{r=1}^n p_{t+1}^{(r)} = (1-\alpha)\beta_t +\alpha,\\
    &\sum_{r=1}^n p_{t+1}^{(r)^2} = (1-\alpha)^2\gamma_t + 2\alpha(1-\alpha)\frac{\beta_t}{n} + \frac{\alpha^2}{n}.
\end{align*}
Notice that it only requires constant time to get the desired values.
\par 
Now suppose $(1-\alpha)w_{t+1}^{(I_{p,t})} + \frac{\alpha}{n} < \frac{\delta}{n}$. Then, it is same as the above except that now we update $\beta_t$ and $\gamma_t$ to $\beta_t \leftarrow \beta_t + \frac{\delta}{n} - p_t^{(I_{p,t})}, \gamma_t \leftarrow \gamma_t +\big(\frac{\delta}{n}\big)^2 - p_t^{(I_{p,t})^2}$. The rest of the update procedure remains same.

\subsection{Index Random Sampling} \label{appn:random sampling}
Here we propose an index random sampling technique (line 4 of Algorithm \ref{alg:Stoc Ax}, line 4 of Algorithm \ref{alg:p under chi}) that we used in our Python code. We can obtain a random sampled index $I_p$ in $O(\log n)$ complexity by using a balanced tree structure \citep{namkoong2016stochastic}. However, balanced search tree implementation in Python turns out to be not efficient due to the frequent use of for and while loops, which are known to be typically slow in Python. Also, if we compare the computation time of random sampling with the tree structure with other parts of the code that is coded with well-optimized and C-written packages such as NumPy, the random sampling with tree structure becomes the bottleneck, although it shouldn't be. So, we instead iteratively update the cumulative distribution of $\frac{p_t}{\mathbbm{1}_n^\top p_t}$ to do the index random sampling.
For simplicity, we omit $i$ which is used for constraint index.
\par 
By using cumulative distribution of $\frac{p_t}{\mathbbm{1}_n^\top p_t}$, we show that a cumulative distribution of $\frac{p_{t+1}}{\mathbbm{1}_n^\top p_{t+1}}$ can be derived efficiently. Suppose $S_t$ be a vector of cumulative distribution of distribution  $\frac{p_t}{\mathbbm{1}_n^\top p_t}$. For each iteration $t\in [T]$, let us define $\beta_t$ as $\beta_t : = \mathbbm{1}_n^\top p_t$ and $I_{p,t}$ be a randomly choesn index from the distribution $\frac{p_t}{\mathbbm{1}_n^\top p_t}$. Suppose $w_t^{(I_{p,t})} \geq \frac{\delta}{n}$, then $I(\alpha) = [n]$. By the virtue of the update, which is $p_{t+1} = (1-\alpha) p_t + \frac{\alpha}{n}$ , we have 
\begin{align*}
    \beta_{t+1}S_{t+1} = (1-\alpha) \beta_{t} S_t + v_1
\end{align*}
where $v_1: = \frac{\alpha}{n}(1,\dots, n) \in \mathbb{R}^n$. 
\par 
Now suppose $(1-\alpha)w_t^{(I_{p,t})} + \frac{\alpha}{n} \leq \frac{\delta}{n}$, then we set $p_{t+1}^{(I_{p,t})} = \frac{\delta}{n}$. Then update rule for each element $k \in [n]$ of the cumulative distribution $S_{t+1}$ would follow
\begin{enumerate}[(i)]
    \item $k<I_{p,t}$ \par 
    $$ \beta_{t+1}S_{t+1}^{(k)} = (1-\alpha) \beta_t S_{t}^{(k)} + \frac{k\alpha }{n}$$
    \item $k = I_{p,t}$\par
     \begin{align*}
         \beta_{t+1}S_{t+1}^{(I_{p,t})} &= (1-\alpha) (\beta_t S_{t}^{(I_{p,t})} - p_t^{(I_{p,t})}) + \frac{I_{p,t}\alpha}{n} + \frac{\delta - \alpha}{n}\\
         & = (1-\alpha) \beta_t S_t^{(I_{p,t})}+ \frac{I_{p,t}\alpha}{n} + C_{v_2}\\
         \text{where } C_{v_2}&= \frac{\delta - \alpha}{n} - (1-\alpha)p_t^{(I_{p,t})}.
     \end{align*}
    \item $k>I_{p,t}$ \par 
    \begin{align*}
        \beta_{t+1} S_{t+1}^{(I_{p,t})} = (1-\alpha) \beta_t S_{t}^{(I_{p,t})} + \frac{k\alpha }{n} + C_{v_2}.
    \end{align*}
\end{enumerate}
So we can simply express as 
\begin{align*}
    \beta_{t+1}S_{t+1} = (1-\alpha) \beta_t S_t + v_1 + C_{v_2}v_2
\end{align*}
where $v_2 : = (0,\dots,0,1,\dots,1) \in \mathbb{R}^n$ and 1 starts from index $I_{p,t}$. After obtaining $S_{t+1}$, we random sample $u$ from $\mathrm{uniform}(0,1)$ and find an index of interval of $S_{t+1}$ that $u$ falls into. This is equivalent to random sampling an index from distribution $\frac{p_{t+1}}{\mathbbm{1}_n^\top p_{t+1}}$. It turns out that this is 3 to 4 times faster than using a line of code, numpy.random.choice(n,$\frac{p_{t+1}}{\mathbbm{1}_n^\top p_{t+1}}$).
\section{Additional Experimental Details}
In this section, we present experimental settings or results not presented in the main body. The code is available at the following repository: \url{https://github.com/HyungkiIm/SFOM-DRO}

\subsection{Parameter Tuning} \label{appn:param tune}
We start by deciding the sample size parameter $K$ and $C_K$. We set $C_K=0.05$ at the beginning. Then we decide the sample size $K$ based on the information that shows us the quality of $\hat{i}_t$. This includes table \ref{table:approximation of i-appn} and figure \ref{fig:K iter test-appn}. We choose $K$ that is small compared to $n$. If there is no such $K$, we increase the value of $C_K$ to $C_K+ 0.05$ and repeat the same procedure. 
\par 
After determining $C_K$ and $K$, we tune the parameter $C_s$ such that $w(\Omega) \equiv C_s$, so that we control the value of $T$. We need sufficiently large but tight $T$ so that Algorithm \ref{alg:SOFO} terminates with the SP gap termination. We start by $C_s = 1$ and if Algorithm \ref{alg:SOFO} does not terminate before $T$ iterations, we change $C_s$ to $\sqrt{2}C_s$. This results in $T$ two times larger than before. We repeat this procedure until we get sufficiently large $T$.
\par 
Lastly, we decide $T_s$ based on the calculation time of the SP gap. Also, we choose $T_s$ proportional to $T$. That is, $T_s = \frac{T}{Q}$ for some $Q \geq 1$. If the calculation time of the SP gap is large, then we choose $Q$ to be small so that we do not calculate the SP gap often. Otherwise, we choose larger $Q$ to terminate Algorithm \ref{alg:SOFO} early. In our experiments, we choose $Q = 5, 10, 20, 40$.

\subsection{Logistic Regression with Fairness Constraints}

\subsubsection{Parameter Setting.}\label{appn:FL param}
In this section, we explain why we choose $c=0.05$ for our experiment and give details of other parameter settings. As mentioned in section \ref{sec:experiment}, people choose $c$ to satisfy certain fairness criteria.  
Here, we use the $p\%$-rule \citep{biddle2017adverse}, which measures the fairness of a classifier by the ratio between the percentage of positively classified samples with certain sensitive attribute value and the percentage of positively classified samples without that sensitive attribute value. The ratio of a classifier that satisfies the $p\%$-rule must be greater than equal to $p\%$. Usually, $80\%$-rule is used in practice. 
\citep{zafar2017fairness,akhtar2021conservative} and people usually tune $c$ to satisfy $80\%$-rule as it is hard to incorporate the $p\%$-rule directly into the formulation. We choose $c=0.05$ as this is the value that satisfies the $80\%$-rule \citep{akhtar2021conservative}.
For every experiment, we choose $\rho =5$ and $\delta =0.95$ for defining the ambiguity set, $\epsilon = 0.02$, and $C_K =0.05$. We take the feasible region as $\Theta= \{\theta \in \bbR^d|\|\theta\|_2\leq 5\log(d)\}$. Lastly, we set $G=0.25$ and $M^i=0.25$ for all $i \in [m]$.

\subsubsection{Additional Results}
\begin{figure}[t]
    \centering
    \includegraphics[width=11cm]{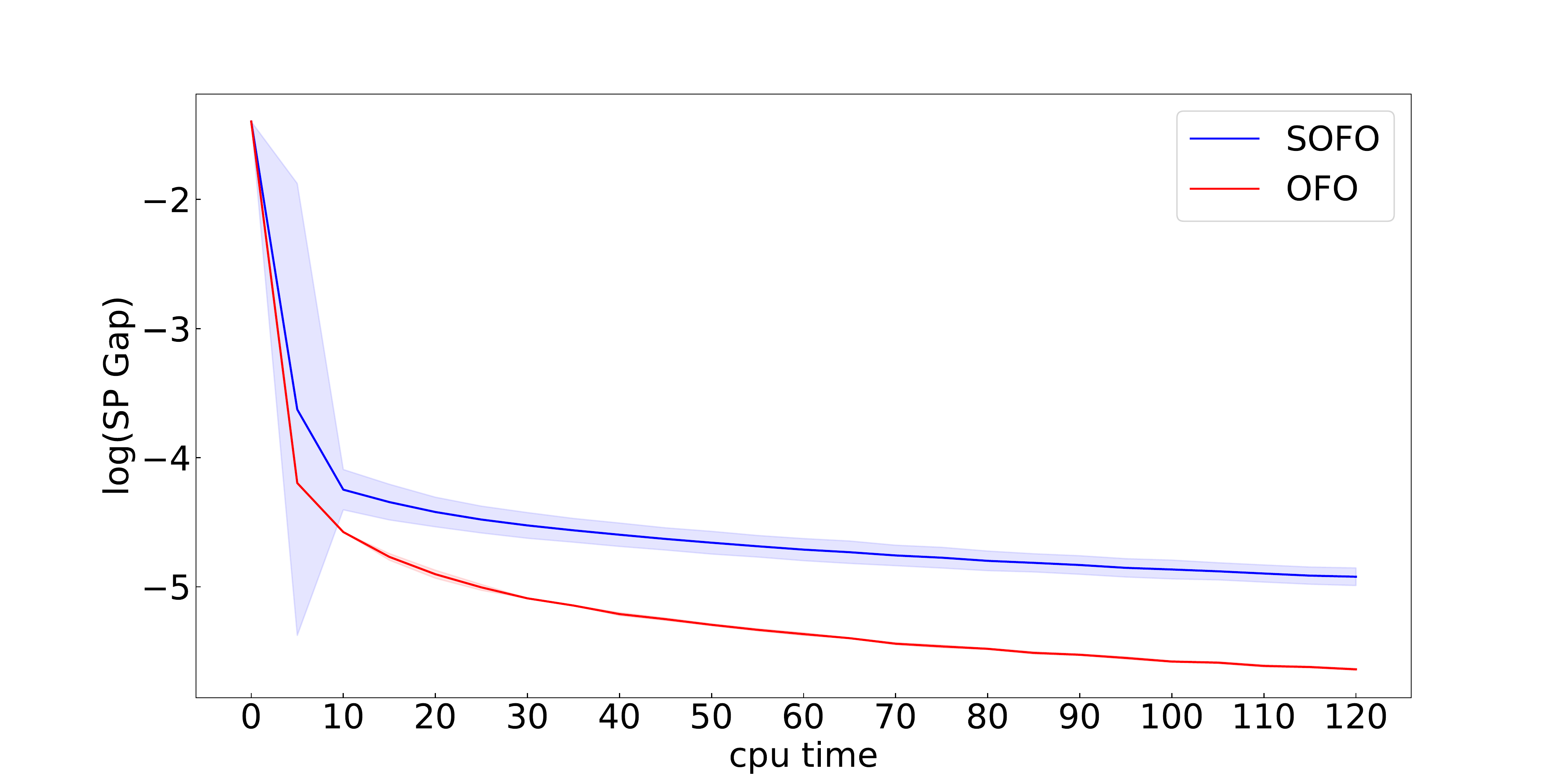}
    \caption{Average log(SP gap) versus cpu time for SOFO and OFO-based approach ($d=300, n=20000$). The shaded area represents the 95\% confidence interval of our test results.}%
    \label{fig:K time test2}%
\end{figure}
In this section, we present some plots and tables that are not presented in the section \ref{sec:experiment}. First of all, figure \ref{fig:K time test2} shows the convergence speed between SOFO and OFO when $n=20000$. Contrary to the result that is presented in figure \ref{fig:K time test}, we can see that the OFO-based approach converges faster than the SOFO-based approach from the beginning of the computation.

\subsection{Parameter Selection in Large-scale Social Network}
\subsubsection{Data Generation}\label{appn:ps datagen}
For each simulation instance, we randomly generate each sample of $\bu^i$ from distribution $\mathcal{N}(\mu^i, \Sigma^i)$.
Also, each element of mean vector $\mu^i \in \bbR^d$ is also randomly generated from distribution $\mathrm{uniform}(0, 1/J)$. The main reason that we are using this uniform distribution is to scale the value of $(p^i)^\top\bF^i(x)$ to $[0,1]$.
If $\mu^i \sim \mathrm{uniform}(0, 1/J)$, for any $x$ that satifies the last two constraints of \eqref{eqn:social network}, we have $\mathbb{E}\big[\bu^{i^\top} x \big] = \frac{1}{2}$. The covariance matrix $\Sigma^i$ is defined as $\Sigma^i = \sigma^2 I_d$, where $I_d$ represents $d \times d$ identity matrix and $\sigma^2$ is a parameter that differs by experiments.

\subsubsection{Parameter Settings}\label{appn:social param}
We use $\delta = 0.9$ and $\rho = 5$ for all the experiments. Also, the constant factors $C_K, C_g, C_s$ are set as $C_K = 0.05,$ $C_g = 1 + \sqrt{\frac{\rho}{n}}$, and $C_s = 1$. Also, in this experiment, we employ the SP gap early termination criterion, and we set the duality gap calculation frequency $T_s$ to $T_s = \frac{T}{40}$. The initializations of $x$ and $p$ are defined as $p_1^i = \frac{1}{n}\mathbbm{1}_n$ and $x_1 = \underbrace{(\frac{1}{L}\mathbbm{1}_L, \dots,\frac{1}{L} \mathbbm{1}_L)}_\text{J times}$ for all the experiments.

\subsubsection{Customized $\epsilon$-SMD with Entropy Function}
This section shows the pseudo-code of $\epsilon$-SMD with entropy distance-generating function and $X = \Delta_d$.
\begin{algorithm}[H]
\caption{$\epsilon$-SMD with Entropy Function}
\begin{algorithmic}[1]
\STATE \textbf{input}: step-size $\alpha_{x,t}$, sample size $K$, $x_{t}$ and $p_t$
\STATE \textbf{output}: $x_{t+1}$

\FOR{$i=1,\dots,m$}
        \STATE Sample $K$ indices from distribution $\frac{p_t^i}{\mathbbm{1}_n^{\top}p_t^i}.$ Let $\mathbf{S}^i$ be an index list of $K$ samples.
\ENDFOR
\STATE compute $\hat{i} \in \underset{i \in [m]}{\text{argmax}} \ \{(\mathbbm{1}_n^\top p_t^i)\hat{f}_{t}^i(x_t) = \frac{\mathbbm{1}_n^{\top}p_t^i}{K_t} \sum\limits_{s \in \mathbf{S}^i} F_s^i(x_t)\}$ \\
\STATE Sample an index $I_{x,t}$ from distribution $\frac{p_t^{\hat{i}}}{\mathbbm{1}_n^{\top}p_t^{\hat{i}}}$. Then get $g_{x,t}$ such that  $g_{x,t} \in \partial \mathbbm{1}_n^{\top}p_t^{\hat{i}} \cdot F_{I_{x,t}}^{\hat{i}}(x_t)$.\\
\STATE Set  $\eta_t^{(r)} = x_t^{(r)}\exp (- \alpha_{x,t} g_{x,t}^{(r)})$, $\forall r \in [d]$
\STATE Set $x_{t+1}^{(r)} = \frac{\eta_t^{(r)}}{\sum_{j=1}^d \eta_t^{(j)}}$, $\forall r \in [d]$
\STATE \textbf{Return} $x_{t+1}$
\end{algorithmic}
\label{alg:Stoc Ax with entropy}
\end{algorithm}
In our experiments of personalized parameter selection in a large-scale social network, the domain of elements of $x$ that correspond to each cohort $j \in [J]$ is $\Delta_L$. So, we apply lines 8-9 of Algorithm \ref{alg:Stoc Ax with entropy} to each corresponding cohort to update the $x_t$ variables in the parameter selection example.

\subsubsection{Experimental Results}
In this section, we present supplementary experimental results that were not included in our main document.
\begin{figure}[t]
    \centering
    \subfloat[n = 5000]{{\includegraphics[width=7.5cm]{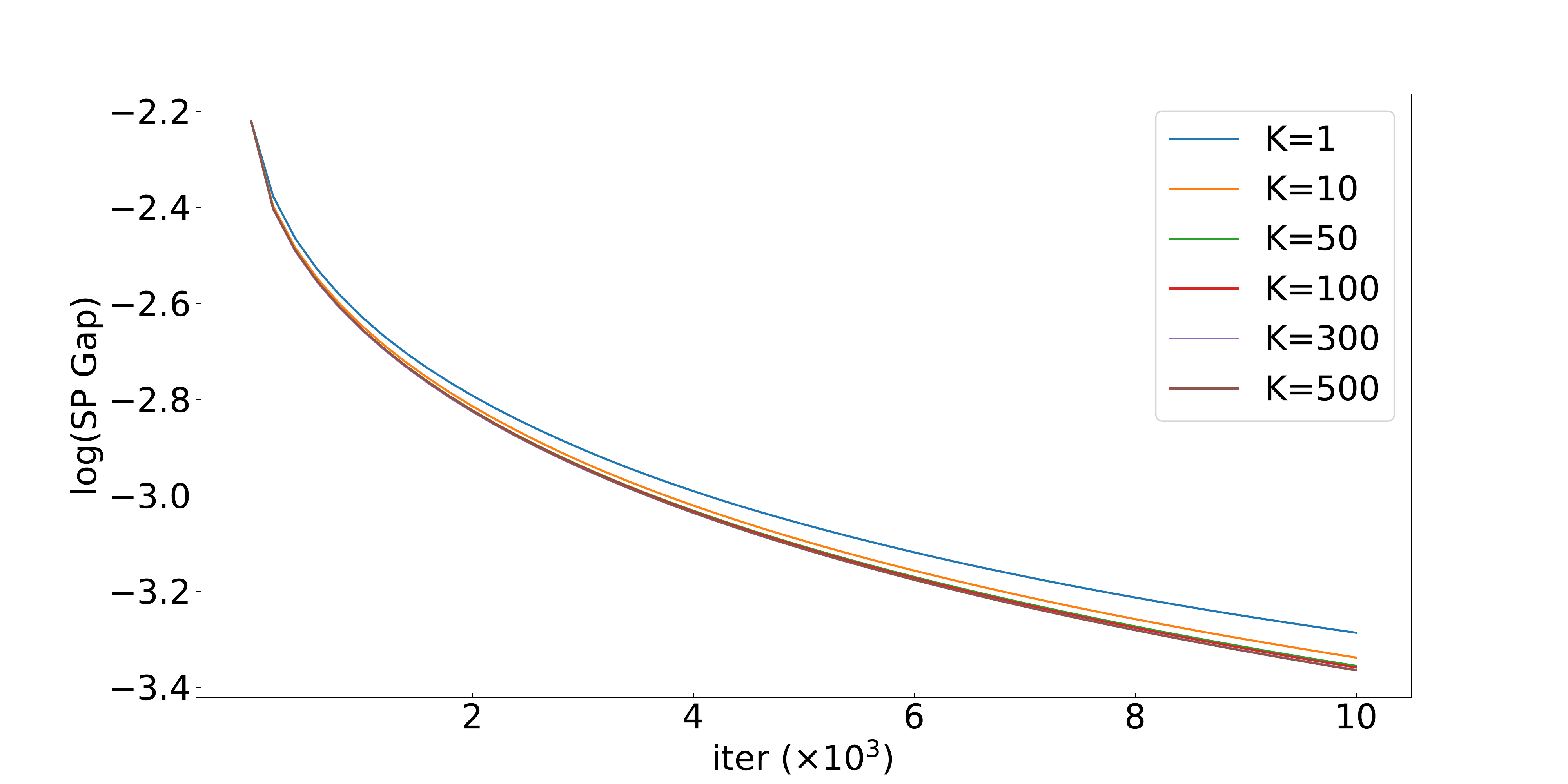} }}%
    \qquad
    \subfloat[n =10000]{{\includegraphics[width=7.5cm]{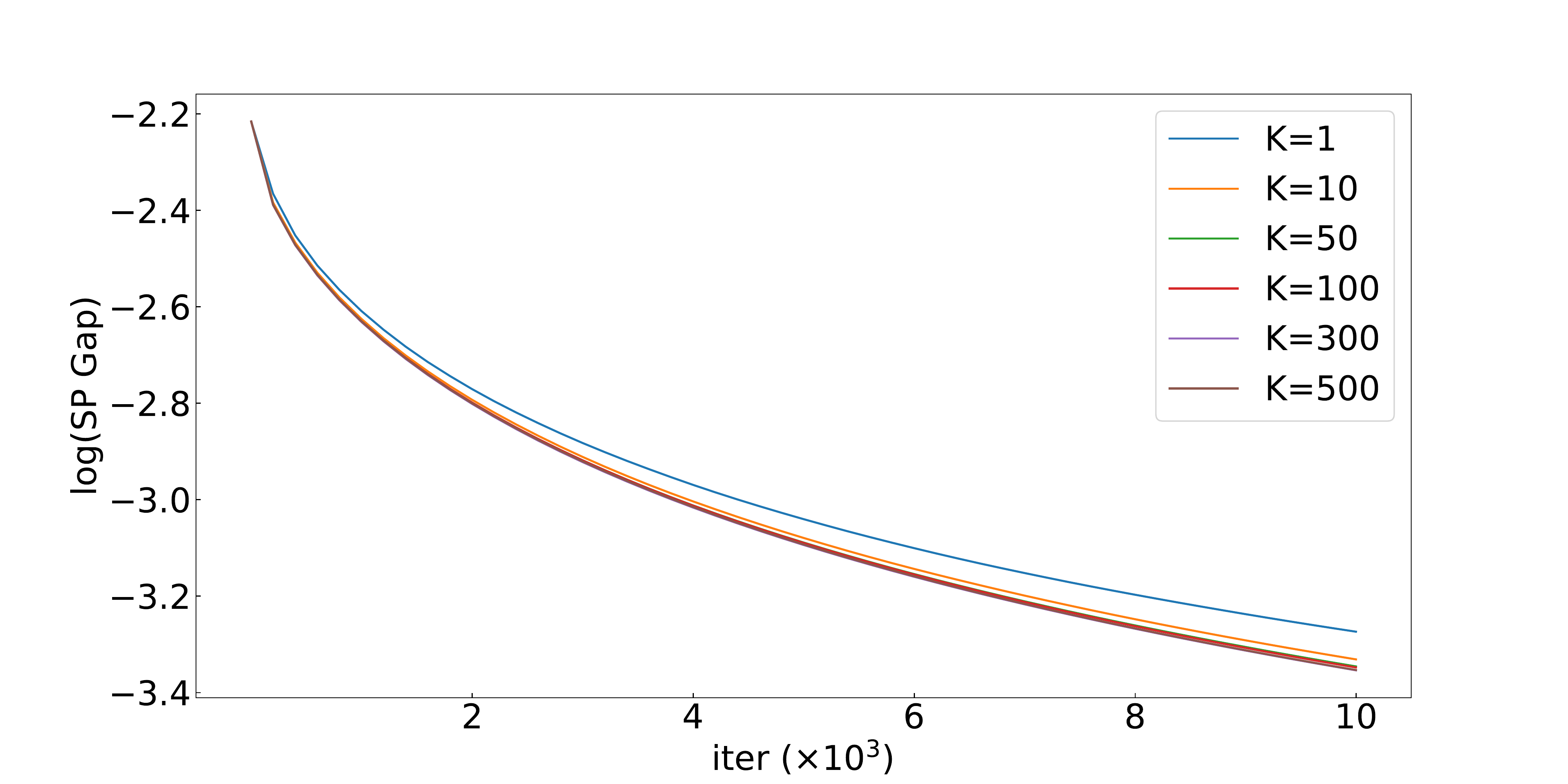} }}%
    \caption{Average log(duality gap) versus number of iterations for different per-iteration sample size SOFO and OFO ($J = 10, L = 25, m = 20, \sigma^2 = 0.2$)}%
    \label{fig:K iter test-appn}%
    \end{figure}
    \par 
As we mentioned earlier, it is essential to choose a moderate per-iteration sample size $K$.
To determine the appropriate per-iteration sample size $K$, we plot the average log duality gap against the number of iterations in figure \ref{fig:K iter test-appn}. We set $c^i = 0.95 \times p_0^{i^\top}\bF^i(x_0)$ to make constraints tight at the start of our algorithm for all $i \in [m]$.
We choose $m$ and $\sigma^2$ to be larger than the other experiments to enlarge an effect of per-iteration sample size. Similar to figure \ref{fig:K Iter}, convergance rate increases as $K$ increases. Moreover, Table \ref{table:approximation of i-appn} shows the average percentage of iterations that satisfy $|f_t^{\hat{i}_t}(x_t) - f_t^{i_t^*}(x_t)| \leq C_K \epsilon$. We can see that when $K\geq 100$, more than 70\% of the iterations satisfied the above inequality, while $\hat{i}_t$ is not approximating $i_t^*$ well when $K<100$. Based on this result, we choose $K=100$ in the following experiment.

\begin{table}[ht]
\begin{center}
\caption{Average percentage of iterations such that $|f_t^{\hat{i}_t}(x_t) - f_t^{i_t^*}(x_t)| \leq C_K \epsilon$ out of 10000 iterations. $(J=10,L=25,m=20,\sigma^2= 0.2)$} 
\begin{tabular}{ |c| c | c | c | c | c | c| }
\toprule
  \backslashbox{$n$}{$K$} & 1 & 10 & 50 & 100 & 300 & 500 \\ 
   \hline
 5000& 0.26 & 6.51 & 43.78  & 72.89 & 98.34 & 99.84\\
 \hline 
 10000 & 0.27 & 6.55 & 43.35 & 72.06 & 98.17 & 99.80\\
 \bottomrule
\end{tabular}
\label{table:approximation of i-appn}
\end{center}
\end{table}
\begin{figure}[ht]
    \centering
    \subfloat[$n = 1000$]{{\includegraphics[width=7.5cm]{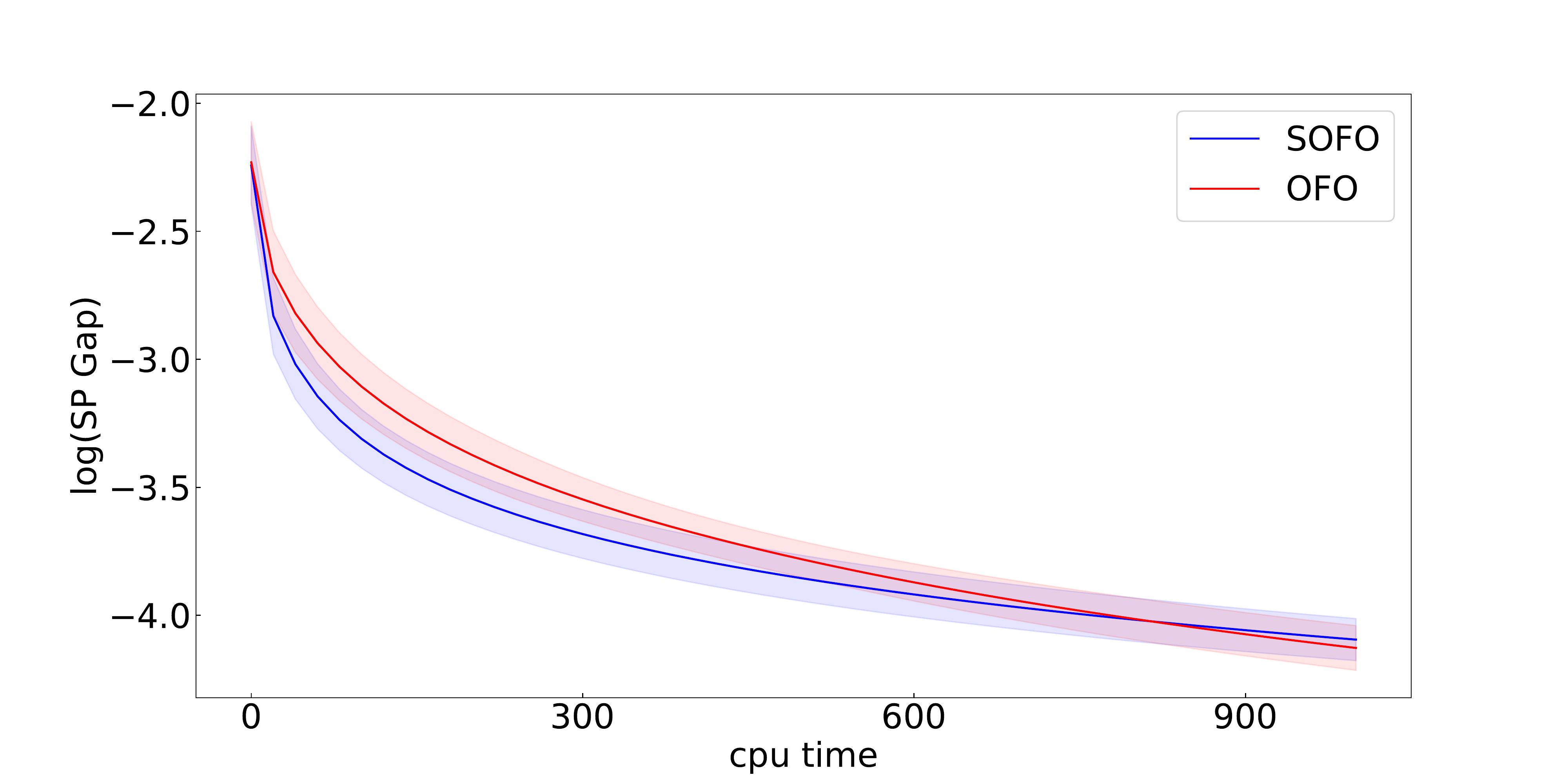} }}%
    \qquad
    \subfloat[$n =5000$]{{\includegraphics[width=7.5cm]{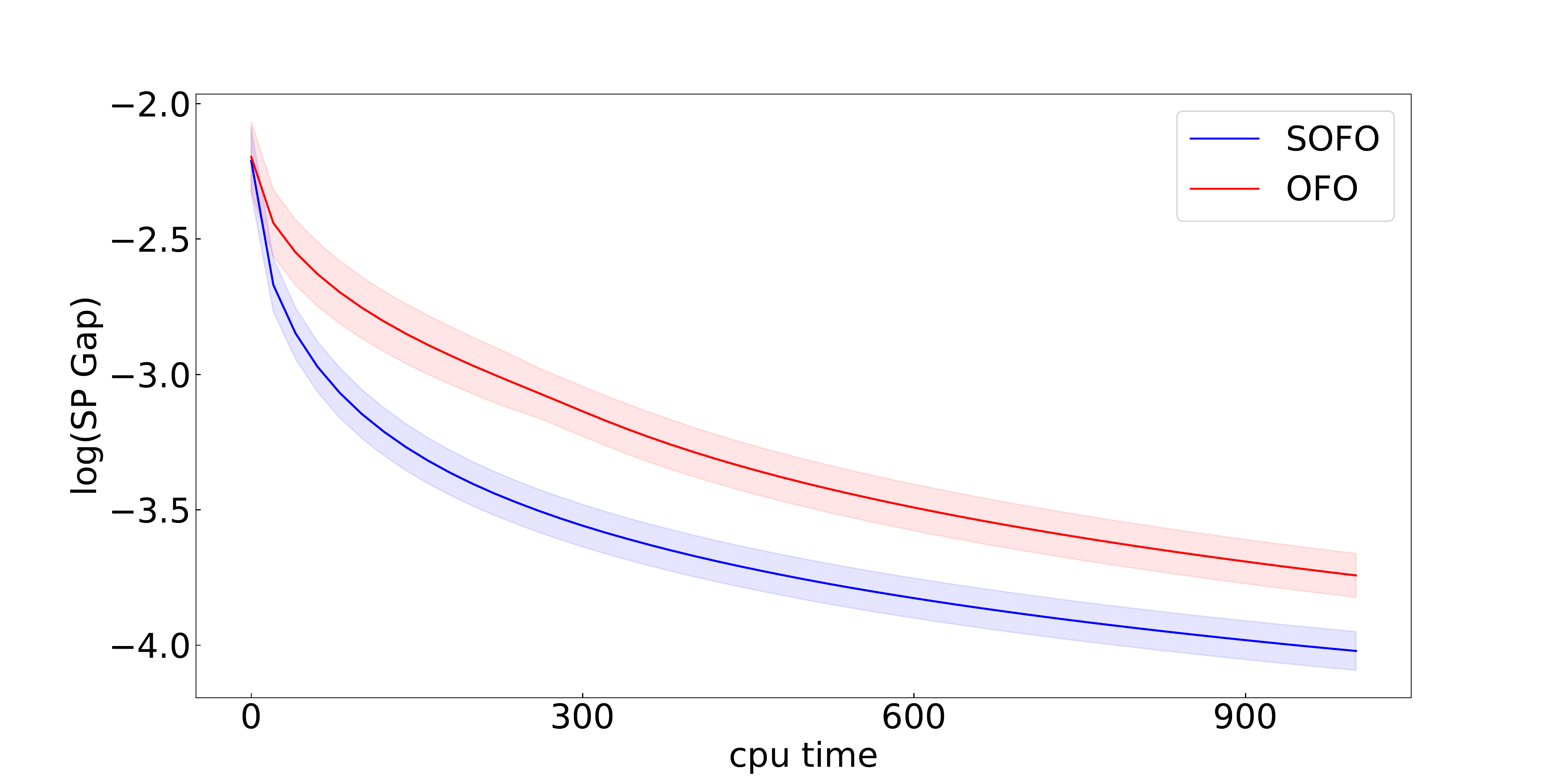} }}%
    \caption{Average log(SP Gap) versus CPU time for SOFO and OFO-based approach. The shaded area represents the 95\% confidence intervals. $(J = 10, L = 25, m = 20, \sigma^2 = 0.2, K= 100)$}%
    \label{fig:K time test-appn}%
    \end{figure} 
\par

\subsection{Multi-item Newsvendor}
\subsubsection{Data Generation.}\label{appn:data gen}
We generate samples of $\xi$ using the following procedure, which largely follows \citep{hanasusanto2015distributionally}.
First, we assume that we have 10 items in our hand and set revenues of all of the items to be 0.5. We fix the salvage price and back-order cost of each items to be $20\%$ and $25\%$ of the corresponding selling prices, respectively. We sample cost vector $c$ from distribution $uniform(0.1,0.25)$ so that our data satisfies $c<r$ and $s<r$. Then we sample mean demand vector $\mu$ from $uniform(0.1,0.2)$ and vector of standard deviation $\sigma$ of the demand from $uniform(0.05 \mu, 0.2 \mu)$. We sample a random matrix $S \in \bbR^{d \times d}$, where each element follows normal distribution, and set $U = S^\top S$ and $u \in \bbR^d$ to be $u = 1/ \sqrt{\mathbf{diag}(U)}$. Then we define a matrix $C$ as $C = \mathrm{diag}(u)U\mathrm{diag}(u)$. We set our covariance matrix $\Sigma$ as $\Sigma  = \mathrm{diag}(\sigma) C \mathrm{diag}(\sigma)$. Finally, for each instance, we generate $n$ samples of $\xi$ from distribution $N(\mu, \Sigma)$.

\subsubsection{Parameter Settings}\label{appn:MINV param}
We choose $\delta=0.9, \rho=5, C_K=0.05,$ and $C_g = 1+\sqrt{\frac{\rho}{n}}$. The initial points of $x$ and $p$ are defined as $x_0 =  \frac{0.1 + 0.2}{2} \mathbbm{1}_d$ and $p_0 = \frac{1}{n} \mathbbm{1}_n$. To ensure we have enough budget, we set $C$ to be $1.2 \times \mathbbm{1}_d^\top \mu$, where $\mu \in \bbR^d$ is the mean demand vector of products. We set $\beta$ to be 0.1 and $\alpha$ to be an approximate upper bound of $\underset{p^1 \sim \calP}{\sup} \tau + \frac{1}{\beta}\mathbb{E}_{\xi^1 \sim p^1}\big[ (L(x,\xi) - \tau)^+ \big]$. We present the details in Appendix \ref{appn:calculation of alpha}.

\subsubsection{Calculating $\underset{x \in X}{\inf} \ \phi(x, \bar{p})$}
In parameter selection example, each element of $\bF^i(x)$ is a linear function of $x$, so \eqref{eqn:inf pi} becomes a linear optimization, which is solvable by using the Gurobi solver.
\par 
For multi-item newsvendor example, let us use $\xi_i$ to denote the $i-$th sample of $\xi$ and let $\bm{\xi}  \in \bbR^{d \times n}$ denote samples of $\xi$, where $\bm{\xi} = \begin{bmatrix} \xi_1 &| \cdots | & \xi_n  \end{bmatrix}$, and also define $\mathbf{L}_{\bm{\xi}}(x): X \rightarrow \bbR^n$ as $\mathbf{L}_{\bm{\xi}}(x) = \begin{bmatrix} L(x,\xi_1)   \\
    \vdots \\ L(x,\xi_n) \end{bmatrix}$. In order to calculate a SP gap, we need to solve $\underset{x \in X}{\inf} \phi(x,\bar{p})$ which is equivalent to:
\begin{equation}
    \begin{aligned}
    &\underset{t, x \in  X}{\inf} \ &&t \\
    &\text{s.t.} && t \geq \bar{p}^{0^\top} \mathbf{L}(x,\bm{\xi}), \\
    & &&t \geq \bar{p}^{1^\top}(\tau \mathbbm{1}_n + \frac{1}{\beta}(\mathbf{L}(x,\bm{\xi}) - \tau \mathbbm{1}_n)_+ - \alpha \mathbbm{1}_n).
    \end{aligned}
\end{equation}
$(\cdot)_+$ function denotes $\max \{\cdot, 0 \}$, where maximum is applied element-wise. Notice that each element of $\mathbf{L}_{\bm{\xi}}(x)$ is like: $L(x,\xi_i) = (c-s)^\top x - (b+r-s)^\top \min \{x, \xi_i \} + b^\top \xi_i.$ In order to solve this problem using Gurobi, we need to replace $ \min \{x, \xi_i \}$ to a vector variable $y \in \bbR^d$ such that $y \geq x$ and $y \geq \xi_i$. As we have $n$ samples, we need to define $d \times n$ variables and the calculation of the duality gap can be costly as $n$ increases.
\par 

\subsubsection{Calculation of $\alpha$} \label{appn:calculation of alpha}
In this section, we study how to set $\alpha$ in multi-item newsvendor model \eqref{form:MINV}. We set $\alpha$ as an approximate upper bound of
\begin{align}
    \underset{p^1 \sim \calP}{\sup} \tau + \frac{1}{\beta}\mathbb{E}_{\xi^1 \sim p^1}\big[ (L(x,\xi) - \tau)^+ \big]. \label{eqn:alpha ub}
\end{align}
To gain an approximate upper boufnd of \eqref{eqn:alpha ub}, we input $x = \mu$, which is the demand vector that we used to generate samples of $\xi$. In fact, the optimal solution $x^*$ might be close to $\mu$. As we can interpret as the average of the $100 \times \beta \%$ worst outcomes of the loss distribution, we input $(1 - \beta)-$th quantile element of a vector $L_{\bm{\xi}}(x)$ to $\tau$. Then, we set $\alpha$ to the average of $\tau + \frac{1}{\beta}\big( (L(x,\xi) - \tau)^+ \big)$, where $\xi$ follows the uniform distribution.

\subsubsection{Choosing $K$ for Multi-item Newsvendor Experiment}\label{appn:result}
Table \ref{table:approximation of i, MINV} shows us the average percentage of iterations that satisfies $|f_t^{\hat{i_t}}(x_t) - f_t^{i_t^*}(x_t)| \leq C_K \epsilon$. 
\begin{table}[h]
\begin{center}
\begin{tabular}{ | c | c | c | c | c| }
\toprule
  \backslashbox{$n$}{$K$} & 50 & 100 & 300 & 500 \\ 
   \hline
 5000& 59.09 & 63.49 & 72.95  & 77.13 \\
 \hline 
 10000 & 58.99 & 63.55 & 72.89 & 77.13\\
 \bottomrule
\end{tabular}
\caption{Average percentage of iterations such that $|f_t^{\hat{i_t}}(x_t) - f_t^{i_t^*}(x_t)| \leq C_K \epsilon$} out of 10000 iterations
\label{table:approximation of i, MINV}
\end{center}
\end{table}
Though $K =50$ only satisfies $|f_t^{\hat{i_t}}(x_t) - f_t^{i_t^*}(x_t)| \leq C_K \epsilon$ about 60\%, the accuracy gain that we achieve by increasing per-iteration sample size $K$ is not significant. Considering we only have two constraints in this example, we set $K$ to be 50.